\documentclass{article}
\usepackage{amsmath}
\usepackage{amssymb}
\usepackage{amsthm}
\usepackage{cite}
\usepackage[arrow,curve,matrix]{xy}
\begin{document}
\def\e#1\e{\begin{equation}#1\end{equation}}
\def\ea#1\ea{\begin{align}#1\end{align}}
\def\eq#1{{\rm(\ref{#1})}}
\theoremstyle{plain}
\newtheorem{thm}{Theorem}[section]
\newtheorem{lem}[thm]{Lemma}
\newtheorem{prop}[thm]{Proposition}
\newtheorem{cor}[thm]{Corollary}
\theoremstyle{definition}
\newtheorem{dfn}[thm]{Definition}
\newtheorem{ass}[thm]{Assumption}
\newtheorem{ex}[thm]{Example}
\def\dim{\mathop{\rm dim}\nolimits}
\def\GL{\mathop{\rm GL}\nolimits}
\def\SL{\mathop{\rm SL}\nolimits}
\def\Sp{\mathop{\rm Sp}\nolimits}
\def\Ker{\mathop{\rm Ker}}
\def\Aut{\mathop{\rm Aut}}
\def\Ad{\mathop{\rm Ad}}
\def\id{{\rm id}}
\def\diag{{\rm diag}}
\def\cha{\mathop{\rm char}}
\def\d{{\rm d}}
\def\na{{\rm na}}
\def\stk{{\rm stk}}
\def\Ho{{\rm Ho}}
\def\Po{{\rm Po}}
\def\uni{{\rm uni}}
\def\rk{{\rm rk}}
\def\re{{\rm re}}
\def\vi{{\rm vi}}
\def\Hom{\mathop{\rm Hom}\nolimits}
\def\End{\mathop{\rm End}}
\def\Iso{\mathop{\rm Iso}\nolimits}
\def\Var{\mathop{\rm Var}\nolimits}
\def\Sch{\mathop{\rm Sch}\nolimits}
\def\Stab{\mathop{\rm Stab}\nolimits}
\def\Spec{\mathop{\rm Spec}\nolimits}
\def\CF{\mathop{\rm CF}\nolimits}
\def\LCF{\mathop{\rm LCF}\nolimits}
\def\SF{\mathop{\rm SF}\nolimits}
\def\LSF{\mathop{\rm LSF}\nolimits}
\def\uSF{\mathop{\underline{\rm SF\!}\,}\nolimits}
\def\oSF{\mathop{\bar{\rm SF}}\nolimits}
\def\uoSF{\mathop{\bar{\underline{\rm SF\!}\,}}\nolimits}
\def\hSF{\mathop{\hat{\rm SF}}\nolimits}
\def\uhSF{\mathop{\hat{\underline{\rm SF\!}\,}}\nolimits}
\def\uLSF{\mathop{\underline{\rm LSF\!}\,}\nolimits}
\def\ge{\geqslant}
\def\le{\leqslant}
\def\bA{{\mathbin{\mathbb A}}}
\def\bF{{\mathbin{\mathbb F}}}
\def\bG{{\mathbin{\mathbb G}}}
\def\Z{{\mathbin{\mathbb Z}}}
\def\Q{{\mathbin{\mathbb Q}}}
\def\C{{\mathbin{\mathbb C}}}
\def\K{{\mathbin{\mathbb K}}}
\def\D{{\mathbin{\mathfrak D}}}
\def\E{{\mathbin{\mathfrak E}}}
\def\F{{\mathbin{\mathfrak F}}}
\def\oF{\mathbin{\overline{\mathfrak F}}}
\def\G{{\mathbin{\mathfrak G}}}
\def\oG{\mathbin{\overline{\mathfrak G}}}
\def\H{{\mathbin{\mathfrak H}}}
\def\cA{{\mathbin{\mathcal A}}}
\def\cP{{\mathbin{\mathcal P}}}
\def\cQ{{\mathbin{\mathcal Q}}}
\def\cR{{\mathbin{\mathcal R}}}
\def\cS{{\mathbin{\mathcal S}}}
\def\R{{\mathbin{\mathfrak R}}}
\def\oR{\mathbin{\overline{\mathfrak R}}}
\def\fS{{\mathbin{\mathfrak S}}}
\def\oS{\mathbin{\overline{\mathfrak S}}}
\def\U{{\mathbin{\mathfrak U\kern .05em}}}
\def\V{{\mathbin{\mathfrak V}}}
\def\al{\alpha}
\def\be{\beta}
\def\ga{\gamma}
\def\de{\delta}
\def\bde{\bar\delta}
\def\ep{\epsilon}
\def\io{\iota}
\def\la{\lambda}
\def\ka{\kappa}
\def\th{\theta}
\def\si{\sigma}
\def\La{\Lambda}
\def\Th{\Theta}
\def\Up{\Upsilon}
\def\Om{\Omega}
\def\ts{\textstyle}
\def\sst{\scriptscriptstyle}
\def\sm{\setminus}
\def\ot{\otimes}
\def\ra{\rightarrow}
\def\longra{\longrightarrow}
\def\t{\times}
\def\ha{{\ts\frac{1}{2}}}
\def\op{\oplus}
\def\ci{\circ}
\def\ti{\tilde}
\def\ab{\allowbreak}
\def\md#1{\vert #1 \vert}
\title{Motivic invariants of Artin stacks \\ and `stack functions'}
\author{Dominic Joyce}
\date{}
\maketitle

\baselineskip 11.65pt plus .25pt

\begin{abstract}
An invariant $\Up$ of quasiprojective $\K$-varieties $X$ with values
in a commutative ring $\La$ is {\it motivic} if
$\Up(X)=\Up(Y)+\Up(X\sm Y)$ for $Y$ closed in $X$, and $\Up(X\t
Y)=\Up(X)\Up(Y)$. Examples include Euler characteristics $\chi$ and
virtual Poincar\'e and Hodge polynomials. We first define a unique
extension $\Up'$ of $\Up$ to finite type Artin $\K$-stacks $\F$,
which is motivic and satisfies $\Up'([X/G])=\Up(X)/\Up(G)$ when $X$
is a $\K$-variety, $G$ a {\it special\/} $\K$-group acting on $X$,
and $[X/G]$ is the quotient stack. This only works if $\Up(G)$ is
invertible in $\La$ for all special $\K$-groups $G$, which excludes
$\Up=\chi$ as $\chi(\bG_m)=0$. But we can extend the construction to
get round this.

Then we develop the theory of {\it stack functions} on Artin stacks.
These are a universal generalization of constructible functions on
Artin stacks. There are several versions of the construction: the
basic one $\SF(\F)$, and variants $\uSF(\F,\Up,\La),\ldots$
`twisted' by motivic invariants. We associate a $\Q$-vector space
$\SF(\F)$ or a $\La$-module $\uSF(\F,\Up,\La)$ to each Artin stack
$\F$, with functorial operations of multiplication, pullbacks
$\phi^*$ and pushforwards $\phi_*$ under 1-morphisms $\phi:\F\ra\G$,
and so on. They will be important tools in the author's series on
`Configurations in abelian categories'.
\end{abstract}

\section{Introduction}
\label{mi1}

Let $\K$ be an algebraically closed field. An invariant $\Up$ of
isomorphism classes $[X]$ of quasiprojective $\K$-varieties $X$
taking values in a commutative ring $\La$ is called {\it motivic} if
whenever $Y\subseteq X$ is a closed subvariety we have
$\Up([X])=\Up([X\sm Y])+\Up([Y])$, and whenever $X,Y$ are varieties
we have $\Up([X\t Y])=\Up([X])\Up([Y])$. The name `motivic' refers
to {\it motives} and {\it motivic integration}, where such
constructions are common. Well-known examples are the {\it Euler
characteristic}, {\it virtual Hodge polynomials} and {\it virtual
Poincar\'e polynomials}.

The first goal of this paper, in \S\ref{mi41}--\S\ref{mi42}, is to
extend such invariants to {\it Artin stacks}. If $\R$ is a finite
type algebraic $\K$-stack with affine geometric stabilizers we
define $\Up'([\R])\in\La$ uniquely with the above motivic
properties, such that if $\R$ is a quotient $[X/G]$ for $X$ a
quasiprojective $\K$-variety and $G$ a {\it special\/} algebraic
$\K$-group, then~$\Up'([\R])=\Up([X])\Up([G])^{-1}$.

Naturally, this is only possible if $\Up([G])$ is invertible in
$\La$ for all special algebraic $\K$-groups $G$. The most important
restriction this imposes is that $\ell-1$ is invertible, where
$\ell=\Up([\bA^1])$. We can arrange this for virtual Hodge and
Poincar\'e polynomials, but not for Euler characteristics, since
then $\ell=1$. Parts of \S\ref{mi5} and \S\ref{mi6} are dedicated to
a version of this construction which allows $\ell=1$, and so defines
a kind of Euler characteristic of Artin stacks.

Very roughly, the idea when $\ell=1$ is that if $T^G$ is a maximal
torus of $G$, then $\Up([G])=\Up([G/T^G])\Up([T^G])$, where
$\Up([G/T^G])$ is invertible in $\La$. So we can write
$\Up'([[X/G]])=\Up([G/T^G])^{-1}\Up'([[X/T^G]])$. Now $[X/T^G]$ is a
finite disjoint union of $\K$-substacks 1-isomorphic to
$Y_i\t[\Spec\K/H_i]$ for quasiprojective $\K$-varieties $Y_i$ and
$\K$-groups $H_i$ of the form $\bG_m^k\t K$ for $K$ finite abelian.
We then define $\Up'([[X/T^G]])=\sum_i\Up([Y_i])[H_i]$, which takes
values in the commutative $\La$-algebra $\bar\La$ with $\La$-basis
isomorphism classes $[H]$ of $\K$-groups $H$ of the form $\bG_m^k\t
K$ for $K$ finite abelian, and products~$[H_1][H_2]=[H_1\t H_2]$.

The above is not quite true: to make $\Up'([[X/G]])$ depend only on
the stack $[X/G]$ and not on $X,G$ we have to introduce in
\S\ref{mi52} the idea of {\it virtual rank}, which treats a
nonabelian $\K$-group $G$ as being a kind of finite $\La$-linear
combination of certain $\K$-subgroups $Q\subseteq T^G$, of which
$T^G$ is the largest. Then $\Up'([[X/G]])$ is a $\La$-linear
combination of $\Up'([[X/Q]])$ over all such~$Q$.

We will apply this in the series \cite{Joyc2,Joyc3,Joyc4,Joyc5}. If
$\cA$ is a $\K$-linear abelian category and $(\tau,T,\le)$ a {\it
stability condition} on $\cA$, we define {\it invariants} of
$\cA,(\tau,T,\le)$ by applying $\Up'$ to the $\K$-stacks ${\rm
Obj}_{\rm ss}^\al(\tau),{\rm Obj}_{\rm st}^\al(\tau)$ of
$\tau$-(semi)stable objects in $\cA$ with class $\al\in K(\cA)$. The
motivic properties of $\Up'$ mean these invariants satisfy
attractive identities and transformation laws, and can be computed
in examples.

The second goal of the paper, in \S\ref{mi3} and
\S\ref{mi43}--\S\ref{mi6}, is to develop the theory of {\it `stack
functions'}. Before discussing this we explain the ideas of
\cite{Joyc1} on {\it constructible functions on stacks}. To each
Artin $\K$-stack $\F$ we associate a $\Q$-algebra $\CF(\F)$ of {\it
constructible functions} on $\F$, spanned by the characteristic
functions of finite type $\K$-substacks $\G\subseteq\F$. If
$\phi:\F\ra\G$ is a 1-morphism we define the {\it pushforward\/}
$\CF^\stk(\phi):\CF(\F)\ra\CF(\G)$ (for $\phi$ representable and
$\cha\K=0$) and the {\it pullback\/} $\phi^*:\CF(\G)\ra\CF(\F)$ (for
$\phi$ of finite type). These have good functorial properties, for
instance $\CF^\stk(\psi\ci\phi)=\CF^\stk(\psi)\ci\CF^\stk(\phi)$,
$(\psi\ci\phi)^*=\phi^*\ci\psi^*$, and pushforwards, pullbacks
commute in Cartesian squares.

Stack functions are a universal generalization of constructible
functions. The basic version, in \S\ref{mi3}, replaces $\CF(\F)$ by
a $\Q$-vector space $\uSF(\F)$ or $\SF(\F)$ spanned by
(representable) 1-morphisms $\rho:\R\ra\F$, for $\R$ of finite type.
These have multiplication and pushforwards and pullbacks along
1-morphisms with the same functoriality properties as constructible
functions, and maps to and from $\CF(\F)$ commuting with
multiplication and pushforwards and pullbacks in various ways. Thus,
stack functions can be used as a substitute for constructible
functions in many problems. But as $\uSF(\F),\SF(\F)$ contain much
more information than $\CF(\F)$ they are a more powerful invariant.
This will be exploited in~\cite{Joyc3,Joyc4,Joyc5}.

For varieties, similar ideas to \S\ref{mi3} can be found in the
subject of {\it motivic integration}. In particular, for a
$\K$-variety $X$, our space $\SF(X)$ agrees with
$K_0(\Var_X)\ot_\Z\Q$, where $K_0(\Var_X)$ is the {\it Grothendieck
group of\/ $X$-varieties} defined by Looijenga \cite[\S 2]{Looi} and
Bittner \cite[\S 5]{Bitt}, and the operations we define on such
$\SF(X)$ agree with operations in \cite[\S 6]{Bitt}. This suggests
our spaces $\uSF,\SF(\F)$ may have applications in the extension of
motivic integration to Artin stacks (see Yasuda \cite{Yasu} for the
extension to Deligne--Mumford stacks).

Sections \ref{mi43}, \ref{mi5} and \ref{mi6} integrate these ideas
with the material of \S\ref{mi41}--\S\ref{mi42} to produce stack
function spaces $\uSF(\F,\Up,\La)$ modifying $\uSF(\F)$ (or
$\uoSF,\oSF(\F,\Up,\La)$, or $\uoSF,\oSF(\F,\Up,\La^\ci)$, or
$\uoSF,\oSF(\F,\Th,\Om)$, or $\uhSF,\hSF(\F,\chi,\Q)$: there are
several different versions), with the same operations and
functoriality properties. Here is one way to motivate these spaces.
The pushforward of constructible functions
$\CF^\stk(\phi):\CF(\F)\ra\CF(\G)$ is defined by `integration' over
the fibres of $\phi$ using the Euler characteristic $\chi$ as
measure.

If $\Up$ is a $\La$-valued motivic invariant as above, we could
instead take $\La$-valued constructible functions
$\CF(\F)_{\sst\La}$, and define pushforwards
$\CF^\stk_{\sst\La}(\phi):\CF(\F)_{\sst\La}\ra\CF(\G)_{\sst\La}$ by
`integration' using $\Up$ as measure. But then $\CF^\stk_{\sst\La}
(\psi\ci\phi)=\CF^\stk_{\sst\La}(\psi)\ci\CF^\stk_{\sst\La}(\phi)$
may no longer hold, as this depends on properties of $\chi$ on
non-Zariski-locally-trivial fibrations which are false for other
$\Up$ such as virtual Poincar\'e polynomials. This is a pity, as
there would be interesting applications such as the Ringel--Hall
algebras in \cite{Joyc3} if functoriality held.

Our spaces $\uSF(\F,\Up,\La),\ldots$ are designed to overcome this
problem. They are a substitute for $\CF(\F)_{\sst\La}$, and would
reduce to $\CF(\F)_{\sst\La}$ if every 1-morphism $\phi:X\ra\F$ for
$X$ a $\K$-variety could be broken into finitely many Zariski
locally trivial fibrations $\phi_i:X_i\ra\F_i\subseteq\F$, but in
general this is impossible. They have important applications in the
author's series \cite{Joyc2,Joyc3,Joyc4,Joyc5}, where we use them to
associate algebras and Lie algebras to a $\K$-linear abelian
category $\cA$, including quantized universal enveloping algebras,
and to define invariants in $\La$ which `count' $\tau$-semistable
objects in~$\cA$.

In a recent paper \cite{Toen} written independently, Toen defines a
{\it Grothendieck ring of Artin $n$-stacks} which is closely related
to ideas below. In particular, \cite[Th.~1.1]{Toen} is similar to
our Theorem \ref{mi4thm1}, with the same hypotheses Assumption
\ref{mi4ass}. Toen's ring $K(\mathcal{CH}^{\rm sp}(k))\ot_\Z\Q$ is
also more-or-less the same thing as $\uSF(\Spec\K,\Up_\uni,
\La_\uni)$, combining Example \ref{mi4ex3} and Definition
\ref{mi4def1} below.

In \cite[\S 2.4]{Joyc5} we will generalize parts of
\S\ref{mi3}--\S\ref{mi4} below. We define spaces of {\it essential
stack functions} $\mathop{\rm ESF}(\F)$ with $\SF(\F)\subseteq
\mathop{\rm ESF}(\F)\subseteq\LSF(\F)$ and a notion of {\it strong
convergence} of infinite sums in $\mathop{\rm ESF}(\F)$, and then we
extend the motivic invariants $\Up'$ of \S\ref{mi4} to $\mathop{\rm
ESF}(\F)$ in such a way that $\Up'$ takes strongly convergent sums
in $\mathop{\rm ESF}(\F)$ to convergent sums in $\La$, and commutes
with taking limits.
\medskip

\noindent{\it Acknowledgements.} I thank Kai Behrend, Tom
Bridgeland, Franziska Heinloth, J\"org Sch\"urmann, Bertrand Toen
and Burt Totaro for useful conversations, and a referee for helpful
comments. I held an EPSRC Advanced Research Fellowship whilst
writing this paper.

\section{Background material}
\label{mi2}

We introduce $\K$-groups and Artin stacks in
\S\ref{mi21}--\S\ref{mi22}, and then review the author's paper
\cite{Joyc1} on constructible functions on stacks in~\S\ref{mi23}.

\subsection{Algebraic $\K$-groups}
\label{mi21}

Let $\K$ be an algebraically closed field. A good reference on
algebraic $\K$-groups is Borel \cite{Bore}. Following Borel, we
define a $\K$-{\it variety} to be a $\K$-scheme which is reduced,
separated, and of finite type. We do not require our $\K$-varieties
to be {\it irreducible}, as many authors do. This allows algebraic
$\K$-groups with more than one connected component as
$\K$-varieties. An {\it algebraic $\K$-group} is then a $\K$-variety
$G$ with identity $1\in G$ (that is, $1:\Spec\K\ra G$),
multiplication $\mu:G\t G\ra G$ and inverse $i:G\ra G$ (as morphisms
of $\K$-varieties) satisfying the usual group axioms. We call $G$
{\it affine} if it is an affine $\K$-variety.

We will need the following notation and facts about algebraic
$\K$-groups and tori. Throughout $G$ is an affine algebraic
$\K$-group.
\begin{itemize}
\setlength{\itemsep}{0pt}
\setlength{\parsep}{0pt}
\item Write $\bG_m$ for $\K\sm\{0\}$ as a $\K$-group under
multiplication. Write $\bA^m$ for affine space $\K^m$, regarded as a
$\K$-variety. If $A$ is a finite-dimensional $\K$-algebra, write
$A^\t$ for the $\K$-group of invertible elements of $A$ under
multiplication.
\item By a {\it torus} we mean an algebraic $\K$-group isomorphic
to $\bG_m^k$ for some $k\ge 0$. A {\it subtorus} of $G$ means a
$\K$-subgroup of $G$ which is a torus.
\item A {\it maximal torus} in $G$ is a subtorus $T^G$ contained in
no larger subtorus $T$ in $G$. All maximal tori in $G$ are conjugate
by Borel \cite[Cor.~IV.11.3]{Bore}. The {\it rank\/} $\rk\,G$ is the
dimension of any maximal torus. A maximal torus in $\GL(m,\K)$ is
the subgroup $\bG_m^m$ of diagonal matrices.
\item Let $T$ be a torus and $H$ a closed $\K$-subgroup of $T$.
Then $H$ is isomorphic to $\bG_m^k\t K$ for some $k\ge 0$ and finite
abelian group~$K$.
\item If $S$ is a subset of $T^G$, define the {\it centralizer} of $S$
in $G$ to be $C_G(S)=\{\ga\in G:\ga s=s\ga$ $\forall s\in S\}$, and
the {\it normalizer} of $S$ in $G$ to be $N_G(S)=\{\ga\in
G:\ga^{-1}S\ga=S\}$. They are closed $\K$-subgroups of $G$
containing $T^G$, and $C_G(S)$ is normal in~$N_G(S)$.
\item The quotient group $W(G,T^G)=N_G(T^G)/C_G(T^G)$ is called the
{\it Weyl group} of $G$. As in \cite[IV.11.19]{Bore} it is a {\it
finite group}, which acts on~$T^G$.
\item Define the {\it centre} of $G$ to be $C(G)=\{\ga\in
G:\ga\de=\de\ga$ $\forall\de\in G\}$. It is a closed $\K$-subgroup
of~$G$.
\item There is a notion \cite[I.4.5]{Bore} of {\it semisimple
elements} $\ga\in G$, which are diagonalizable in any representation
of $G$. (It is essential that $G$ is {\it affine} here.) Morphisms
of affine algebraic $\K$-groups take semisimple elements to
semisimple elements, \cite[Th.~I.4.4(4)]{Bore}. If $G$ is {\it
connected\/} then $\ga\in G$ is semisimple if and only if it lies in
a maximal torus of $G$, \cite[Th.~IV.11.10]{Bore}.
\end{itemize}

We will also need the notion of {\it special\/} algebraic
$\K$-group, which is studied in the articles by Serre and
Grothendieck in the Chevalley seminar~\cite[\S\S 1, 5]{Chev}.

\begin{dfn} An algebraic $\K$-group $G$ is {\it special\/} if every
principal $G$-bundle over a $\K$-variety is locally trivial in the
Zariski topology.
\label{mi2def1}
\end{dfn}

The following facts may be found in \cite[\S\S 1.4, 1.5 \&
5.5]{Chev}, or easily deduced.
\begin{itemize}
\setlength{\itemsep}{0pt}
\setlength{\parsep}{0pt}
\item $\bG_m$, $\bG_m^n$ and $\GL(m,\K)$ are special. If $A$ is a
finite-dimensional $\K$-algebra then $A^\t$ is special. Products of
special $\K$-groups are special.
\item A $\K$-group $G$ is special if and only if it admits an
embedding $G\subseteq\GL(m,\K)$ with the $G$-principal bundle
$\GL(m,\K)\ra\GL(m,\K)/G$ Zariski locally trivial. If this holds for
some embedding $G\subseteq\GL(m,\K)$ it holds for any
embedding~$G\subseteq\GL(n,\K)$.
\item Special $\K$-groups are always {\it affine} and
{\it connected}. A {\it semisimple} $\K$-group is special if and
only if it is isomorphic to a product of $\K$-groups of the form
$\SL(m,\K)$ and $\Sp(2n,\K)$. {\it Connected, soluble} $\K$-groups
are special. If $H$ is normal in $G$ with $H,G/H$ special then $G$
is special.
\end{itemize}

\subsection{Introduction to Artin $\K$-stacks}
\label{mi22}

Fix an algebraically closed field $\K$ throughout. There are four
main classes of `spaces' over $\K$ used in algebraic geometry, in
increasing order of generality:
\begin{equation*}
\text{$\K$-varieties}\subset \text{$\K$-schemes}\subset
\text{algebraic $\K$-spaces}\subset \text{algebraic $\K$-stacks}.
\end{equation*}

{\it Algebraic stacks} (also known as Artin stacks) were introduced
by Artin, generalizing {\it Deligne--Mumford stacks}. Our principal
reference is Laumon and Moret-Bailly \cite{LaMo}, and a good
introduction is provided by G\'omez \cite{Gome}. Following
\cite{LaMo,Gome} we include in the definition of an algebraic stack
$\F$ that the diagonal morphism $\Delta_\F$ is representable,
quasi-compact and separated, but probably the separatedness
assumption can be omitted. We make the convention that all algebraic
$\K$-stacks in this paper are {\it locally of finite type}, and
$\K$-substacks are {\it locally closed}.

Algebraic $\K$-stacks form a 2-{\it category}. That is, we have {\it
objects} which are $\K$-stacks $\F,\G$, and also two kinds of
morphisms, 1-{\it morphisms} $\phi,\psi:\F\ra\G$ between
$\K$-stacks, and 2-{\it morphisms} $A:\phi\ra\psi$ between
1-morphisms. An analogy to keep in mind is a 2-category of
categories, where objects are categories, 1-morphisms are functors
between the categories, and 2-morphisms are isomorphisms (natural
transformations) between functors.

We define the set of $\K$-{\it points} of a stack.

\begin{dfn} Let $\F$ be a $\K$-stack. Write $\F(\K)$ for groupoid of
1-morphisms $x:\K\ra\F$, and $\oF(\K)$ for the set of isomorphism
classes in $\F(\K)$, so that elements of $\oF(\K)$ are 2-isomorphism
classes $[x]$ of 1-morphisms $x:\Spec\K\ra\F$. Elements of $\oF(\K)$
are called $\K$-{\it points}, or {\it geometric points}, of $\F$. If
$\phi:\F\ra\G$ is a 1-morphism then composition with $\phi$ induces
a map of sets~$\phi_*:\oF(\K)\ra\oG(\K)$.

Let $\F$ be an algebraic $\K$-stack and $x:\Spec\K\ra\F$ a
1-morphism. Then the group of 2-morphisms $x\ra x$ has the structure
of a {\it group $\K$-scheme}, which is not necessarily reduced.
Define $\Aut_\K(x)$ to be the associated reduced group $\K$-scheme.
Then $\Aut_\K(x)$ is an {\it algebraic $\K$-group}, which we call
the {\it stabilizer group} of $x$. We say that $\F$ {\it has affine
geometric stabilizers} if $\Aut_\K(x)$ is an affine algebraic
$\K$-group for all 1-morphisms $x:\Spec\K\ra\F$. As an algebraic
$\K$-group up to isomorphism, $\Aut_\K(x)$ depends only on the
isomorphism class $[x]\in\oF(\K)$ of $x\in\F(\K)$. If $\phi:\F\ra\G$
is a 1-morphism, composition induces a morphism of algebraic
$\K$-groups $\phi_*:\Aut_\K([x])\ra\Aut_\K \bigr(\phi_*([x])\bigr)$,
for~$[x]\in\oF(\K)$.
\label{mi2def2}
\end{dfn}

One important difference in working with 2-categories rather than
ordinary categories is that in diagram-chasing one only requires
1-morphisms to be 2-{\it isomorphic} rather than {\it equal}. The
simplest kind of {\it commutative diagram} is:
\begin{equation*}
\xymatrix@R=6pt{
& \G \ar@{=>}[d]^{\,F} \ar[dr]^\psi \\
\F \ar[ur]^\phi \ar[rr]_\chi && \H, }
\end{equation*}
by which we mean that $\F,\G,\H$ are $\K$-stacks, $\phi,\psi,\chi$
are 1-morphisms, and $F:\psi\ci\phi\ra\chi$ is a 2-isomorphism.
Usually we omit $F$, and mean that~$\psi\ci\phi\cong\chi$.

\begin{dfn} Let $\phi:\F\ra\H$, $\psi:\G\ra\H$ be 1-morphisms
of $\K$-stacks. Then one can define the {\it fibre product stack\/}
$\F\t_{\phi,\H,\psi}\G$, or $\F\t_\H\G$ for short, with 1-morphisms
$\pi_\F,\pi_\G$ fitting into a commutative diagram:
\e
\begin{gathered}
\xymatrix@R=-4pt{
& \F \ar[dr]^\phi \ar@{=>}[dd] \\
\F\t_\H\G
\ar[dr]_{\pi_\G} \ar[ur]^{\pi_\F} && \H.\\
& \G \ar[ur]_\psi \\
}
\end{gathered}
\label{mi2eq1}
\e
A commutative diagram
\begin{equation*}
\xymatrix@R=-4pt{
& \F \ar[dr]^\phi \ar@{=>}[dd] \\
\E
\ar[dr]_\eta \ar[ur]^\th && \H\\
& \G \ar[ur]_\psi \\
}
\end{equation*}
is a {\it Cartesian square} if it is isomorphic to \eq{mi2eq1}, so
there is a 1-isomorphism $\E\cong\F\t_\H\G$. Cartesian squares may
also be characterized by a universal property.
\label{mi2def3}
\end{dfn}

\subsection{Constructible functions on stacks}
\label{mi23}

Finally we discuss {\it constructible functions} on $\K$-stacks,
following \cite{Joyc1}. For this section we need $\K$ to have {\it
characteristic zero}.

\begin{dfn} Let $\F$ be an algebraic $\K$-stack. We call
$C\subseteq\oF(\K)$ {\it constructible} if $C=\bigcup_{i\in I}
\oF_i(\K)$, where $\{\F_i:i\in I\}$ is a finite collection of finite
type algebraic $\K$-substacks $\F_i$ of $\F$. We call
$S\subseteq\oF(\K)$ {\it locally constructible} if $S\cap C$ is
constructible for all constructible~$C\subseteq\oF(\K)$.

A function $f:\oF(\K)\ra\Q$ is called {\it constructible} if
$f(\oF(\K))$ is finite and $f^{-1}(c)$ is a constructible set in
$\oF(\K)$ for each $c\in f(\oF(\K))\sm\{0\}$. A function
$f:\oF(\K)\ra\Q$ is called {\it locally constructible} if
$f\cdot\de_C$ is constructible for all constructible
$C\subseteq\oF(\K)$, where $\de_C$ is the characteristic function of
$C$. Write $\CF(\F)$ and $\LCF(\F)$ for the $\Q$-vector spaces of
$\Q$-valued constructible and locally constructible functions on
$\F$. They are closed under multiplication.
\label{mi2def4}
\end{dfn}

We explain {\it pushforwards} and {\it pullbacks} of constructible
functions along a 1-morphism $\phi:\F\ra\G$, following~\cite[Def.s
4.8, 5.1 \& 5.5]{Joyc1}.

\begin{dfn} Let $\F$ be an algebraic $\K$-stack with affine
geometric stabilizers and $C\subseteq\oF(\K)$ be constructible. Then
\cite[Def.~4.8]{Joyc1} defines the {\it na\"\i ve Euler
characteristic} $\chi^\na(C)$ of $C$. It is called {\it na\"\i ve}
as it takes no account of stabilizer groups. For $f\in\CF(\F)$,
define $\chi^\na(\F,f)$ in $\Q$ by
\begin{equation*}
\chi^\na(\F,f)=\ts\sum_{c\in f(\oF(\K))\sm\{0\}}c\,\chi^\na
\bigl(f^{-1}(c)\bigr).
\end{equation*}

Let $\phi:\F\ra\G$ be a 1-morphism between algebraic $\K$-stacks
with affine geometric stabilizers. For $f\in\CF(\F)$, define
$\CF^\na(\phi)f:\oG(\K)\ra\Q$ by
\begin{equation*}
\CF^\na(\phi)f(y)=\chi^\na\bigl(\F,f\cdot \de_{\phi_*^{-1}(y)}\bigr)
\quad\text{for $y\in\oG(\K)$,}
\end{equation*}
where $\de_{\smash{\phi_*^{-1}(y)}}$ is the characteristic function
of $\phi_*^{-1}(\{y\})\subseteq\oG(\K)$ on $\oG(\K)$. Then
$\CF^\na(\phi):\CF(\F)\ra\CF(\G)$ is a $\Q$-linear map called the
{\it na\"\i ve pushforward}.

Now suppose $\phi$ is {\it representable}. Then for any
$x\in\oF(\K)$ we have an injective morphism
$\phi_*:\Aut_\K(x)\ra\Aut_\K\bigl( \phi_*(x)\bigr)$ of affine
algebraic $\K$-groups. The image $\phi_*\bigl(\Aut_\K(x)\bigr)$ is
an affine algebraic $\K$-group closed in
$\Aut_\K\bigl(\phi_*(x)\bigr)$, so the quotient
$\Aut_\K\bigl(\phi_*(x)\bigr)/\phi_*\bigl(\Aut_\K(x)\bigr)$ exists
as a quasiprojective $\K$-variety. Define a function
$m_\phi:\oF(\K)\ra\Z$ by $m_\phi(x)=\chi\bigl(\Aut_\K(\phi_*(x))
/\phi_*(\Aut_\K(x))\bigr)$ for $x\in\oF(\K)$. For $f\in\CF(\F)$,
define $\CF^\stk(\phi)f:\oG(\K)\ra\Q$ by $\CF^\stk(\phi)f=\CF^\na
(\phi)(m_\phi\cdot f)$. Then $\CF^\stk(\phi):\CF(\F)\ra\CF(\G)$ is a
$\Q$-linear map called the {\it stack pushforward}.

Let $\phi$ be of {\it finite type}, not necessarily representable.
If $C\subseteq\oG(\K)$ is constructible then so is
$\phi_*^{-1}(C)\subseteq\oF(\K)$. Thus if $f\in\CF(\G)$ then
$f\ci\phi_*\in\CF(\F)$. Define the {\it pullback\/}
$\phi^*:\CF(\G)\ra\CF(\F)$ by $\phi^*(f)=f\ci\phi_*$. It is
$\Q$-linear.
\label{mi2def5}
\end{dfn}

Here \cite[Th.s 4.9, 5.4, 5.6 \& Def.~5.5]{Joyc1} are some
properties of these.

\begin{thm} Let\/ $\E,\F,\G,\H$ be algebraic $\K$-stacks with
affine geometric stabilizers, and\/ $\be:\F\ra\G$, $\ga:\G\ra\H$ be
$1$-morphisms. Then
\ea
\CF^\na(\ga\ci\be)&=\CF^\na(\ga)\ci\CF^\na(\be):\CF(\F)\ra\CF(\H),
\label{mi2eq2}\\
\CF^\stk(\ga\ci\be)&=\CF^\stk(\ga)\ci\CF^\stk(\be):\CF(\F)\ra\CF(\H),
\label{mi2eq3}\\
(\ga\ci\be)^*&=\be^*\ci\ga^*:\CF(\H)\ra\CF(\F),
\label{mi2eq4}
\ea
supposing $\be,\ga$ representable in \eq{mi2eq3}, and of finite type
in \eq{mi2eq4}. If
\e
\begin{gathered}
\xymatrix{
\E \ar[r]_\eta \ar[d]^\th & \G \ar[d]_\psi \\
\F \ar[r]^\phi & \H }
\end{gathered}
\quad
\begin{gathered}
\text{is a Cartesian square with}\\
\text{$\eta,\phi$ representable and}\\
\text{$\th,\psi$ of finite type, then}\\
\text{the following commutes:}
\end{gathered}
\quad
\begin{gathered}
\xymatrix@C=35pt{
\CF(\E) \ar[r]_{\CF^\stk(\eta)} & \CF(\G) \\
\CF(\F) \ar[r]^{\CF^\stk(\phi)} \ar[u]_{\th^*} & \CF(\H).
\ar[u]^{\psi^*} }
\end{gathered}
\label{mi2eq5}
\e
\label{mi2thm1}
\end{thm}

As discussed in \cite[\S 3.3]{Joyc1} for the $\K$-scheme case,
equation \eq{mi2eq3} is {\it false} for algebraically closed fields
$\K$ of characteristic $p>0$. In \cite[\S 5.3]{Joyc1} we extend all
these results to {\it locally constructible functions}. The main
differences are in which 1-morphisms must be of finite type.

\section{Stack functions, the basic version}
\label{mi3}

We now introduce {\it stack functions}, a universal generalization
of constructible functions with similar properties under
multiplication, pushforwards and pullbacks. Here we study the basic
versions $\uSF(\F),\SF(\F)$, and in \S\ref{mi4}--\S\ref{mi6} we
generalize them to more complicated spaces $\uSF(\F,\Up,\La),
\ldots$. Throughout $\K$ will be an algebraically closed field of
{\it arbitrary characteristic}, except when we specify $\cha\K=0$
for results comparing stack and constructible functions. The
assumption that all $\K$-stacks are locally of finite type can be
relaxed too. For some related constructions for $\K$-varieties
rather than $\K$-stacks, see Bittner~\cite[\S 5--\S 6]{Bitt}.

\begin{dfn} Let $\F$ be an algebraic $\K$-stack with affine geometric
stabilizers. Consider pairs $(\R,\rho)$, where $\R$ is a finite type
algebraic $\K$-stack with affine geometric stabilizers and
$\rho:\R\ra\F$ is a 1-morphism. We call two pairs $(\R,\rho)$,
$(\R',\rho')$ {\it equivalent\/} if there exists a 1-isomorphism
$\io:\R\ra \R'$ such that $\rho'\ci\io$ and $\rho$ are 2-isomorphic
1-morphisms $\R\ra\F$. Write $[(\R,\rho)]$ for the equivalence class
of $(\R,\rho)$. If $(\R,\rho)$ is such a pair and $\fS$ is a closed
$\K$-substack of $\R$ then $(\fS,\rho\vert_\fS)$, $(\R\sm\fS,\rho
\vert_{\R\sm\fS})$ are pairs of the same kind. Define
\begin{itemize}
\setlength{\itemsep}{0pt}
\setlength{\parsep}{0pt}
\item[(a)] $\uSF(\F)$ to be the $\Q$-vector space generated by
equivalence classes $[(\R,\rho)]$ as above, with for each closed
$\K$-substack $\fS$ of $\R$ a relation
\e
[(\R,\rho)]=[(\fS,\rho\vert_\fS)]+[(\R\sm\fS,\rho\vert_{\R\sm\fS})].
\label{mi3eq1}
\e
\item[(b)] $\SF(\F)$ to be the $\Q$-vector space generated
by $[(\R,\rho)]$ with $\rho$ {\it representable}, with the same
relations~\eq{mi3eq1}.
\end{itemize}
Define a {\it multiplication} `$\,\cdot\,$' on $\uSF(\F)$ analogous
to multiplication of functions by
\e
[(\R,\rho)]\cdot[(\fS,\si)]=[(\R\t_{\rho,\F,\si}\fS,\rho\ci\pi_\R)].
\label{mi3eq2}
\e
This is compatible with the relations \eq{mi3eq1}, and so extends to
a $\Q$-bilinear product $\uSF(\F)\t\uSF(\F)\ra\uSF(\F)$. If
$\rho,\si$ are representable then so is $\rho\ci\pi_\R$, so
$\SF(\F)$ is closed under `$\,\cdot\,$'. As $\rho\ci\pi_\R$ is
2-isomorphic to $\si\ci\pi_\fS$, `$\,\cdot\,$' is {\it commutative},
and one can show it is {\it associative} using properties of fibre
products.
\label{mi3def1}
\end{dfn}

The assumption that $\R,\F$ have {\it affine geometric stabilizers}
here will be used in this section only in the results below
comparing $\uSF(\F),\SF(\F)$ and $\CF(\F)$ --- in particular,
without it the linear maps $\pi_\F^\na:\uSF(\F)\ra\CF(\F)$ and
$\pi_\F^\stk:\SF(\F)\ra\CF(\F)$ in Definition \ref{mi3def2} would
not be well-defined. But in \S\ref{mi4}--\S\ref{mi6} we use the
assumption in a much more essential way.

We refer to elements of $\uSF(\F),\SF(\F)$ as {\it stack functions}.
There is an obvious inclusion $\SF(\F)\subset\uSF(\F)$. We could
instead work over $\Z$ rather than $\Q$, and define $\uSF(\F)_\Z$ to
be the abelian group generated by equivalence classes $[(\R,\rho)]$
of pairs $(\R,\rho)$ with relations \eq{mi3eq1}, so that
$\uSF(\F)=\uSF(\F)_\Z\ot_\Z\Q$, and so on. Or we could work over any
ring or abelian group. But for simplicity we consider only $\Q$. We
define maps between $\CF(\F)$ and~$\SF(\F),\uSF(\F)$.

\begin{dfn} Let $\F$ be an algebraic $\K$-stack with affine geometric
stabilizers and $C\subseteq\oF(\K)$ a constructible subset. Then we
may write $C=\coprod_{i=1}^n\oR_i(\K)$, for $\R_1,\ldots,\R_n$
finite type $\K$-substacks of $\F$. Let $\rho_i:\R_i\ra\F$ be the
inclusion 1-morphism. Then $\rho_i$ is representable, so
$[(\R_i,\rho_i)]$ lies in $\SF(\F)\subseteq\uSF(\F)$ by Definition
\ref{mi3def1}. Define
\begin{equation*}
\bde_C=\ts\sum_{i=1}^n[(\R_i,\rho_i)]\in\SF(\F)\subseteq\uSF(\F).
\end{equation*}
We think of this stack function as the analogue of the
characteristic function $\de_C\in\CF(\F)$ of $C$. Using \eq{mi3eq1}
and the argument of \cite[Def.~3.7]{Joyc1} we find that $\bde_C$ is
independent of the choice of decomposition
$C=\coprod_{i=1}^n\oR_i(\K)$, and so is well-defined.

Define a $\Q$-linear map $\io_\F:\CF(\F)\ra\SF(\F)\subseteq\uSF(\F)$
by
\begin{equation*}
\io_\F(f)=\ts\sum_{0\ne c\in f(\oF(\K))}c\cdot\bde_{f^{-1}(c)}.
\end{equation*}
This is well-defined as $f(\oF(\K))$ is finite and $f^{-1}(c)$
constructible for all $0\ne c\in f(\oF(\K))$. Since $f=\sum_{0\ne
c\in f(\oF(\K))}c\cdot\de_{f^{-1}(c)}$, $\io_\F$ is the unique
$\Q$-linear map which takes $\de_C$ to $\bde_C$ for all
constructible $C\subseteq\oF(\K)$. When $\K$ has {\it characteristic
zero}, define $\Q$-linear maps $\pi_\F^\na: \uSF(\F)\ra\CF(\F)$ and
$\pi_\F^\stk:\SF(\F)\ra\CF(\F)$ by
\e
\begin{split}
\pi_\F^\na\bigl(\ts\sum_{i=1}^nc_i[(\R_i,\rho_i)]\bigr)&=
\ts\sum_{i=1}^nc_i\CF^\na(\rho_i)1_{\R_i}\\
\text{and}\qquad
\pi_\F^\stk\bigl(\ts\sum_{i=1}^nc_i[(\R_i,\rho_i)]\bigr)&=
\ts\sum_{i=1}^nc_i\CF^\stk(\rho_i)1_{\R_i},
\end{split}
\label{mi3eq3}
\e
where $1_{\R_i}$ is the function 1 in $\CF(\R_i)$, which is
constructible as $\R_i$ is of finite type. Here in the second line
$\rho_i$ is representable by definition of $\SF(\F)$, so
$\CF^\stk(\rho_i)1_{\R_i}$ makes sense. To see \eq{mi3eq3} is
well-defined, note that if $\R,\rho$ are as in Definition
\ref{mi3def1} and $\fS$ is a closed $\K$-substack of $\R$ then
\begin{equation*}
\CF^\na(\rho)1_\R\!=\!\CF^\na(\rho)(\de_{\oS(\K)}\!+
\!\de_{\overline{(\R\sm\fS)}(\K)})\!=\!\CF^\na(\rho\vert_\fS)1_\fS
\!+\!\CF^\na(\rho\vert_{\R\sm\fS})1_{\R\sm\fS}.
\end{equation*}
So $\pi_\F^\na$ is still well-defined after quotienting by relations
\eq{mi3eq1}, and for representable $\rho$ the same argument works
for~$\pi_\F^\stk$.
\label{mi3def2}
\end{dfn}

\begin{prop} For $\K$ of characteristic zero, $\pi_\F^\na\ci\io_\F$
and\/ $\pi_\F^\stk\ci\io_\F$ are both the identity on $\CF(\F)$.
Hence $\io_\F$ is injective and\/ $\pi_\F^\na,\pi_\F^\stk$ are
surjective. Also $\io_\F,\pi_\F^\stk$ commute with multiplication
in~$\CF(\F),\SF(\F)$.
\label{mi3prop1}
\end{prop}

\begin{proof} If $\R$ is a finite type $\K$-substack in $\F$
with inclusion 1-morphism $\rho:\R\ra\F$ then
$\io_\F(\de_{\oR(\K)})=[(\R,\rho)]$, and
$\pi_\F^\na\bigl([(\R,\rho)]\bigr)=
\pi_\F^\stk\bigl([(\R,\rho)]\bigr)=\de_{\oR(\K)}$. Thus
$\pi_\F^\na\ci\io_\F,\pi_\F^\stk\ci\io_\F$ take $\de_{\oR(\K)}$ to
itself. As such $\de_{\oR(\K)}$ generate $\CF(\F)$, we see that
$\pi_\F^\na\ci\io_\F,\pi_\F^\stk\ci\io_\F$ are the identity, so
$\io_\F$ is injective and $\pi_\F^\na,\pi_\F^\stk$ surjective.

If $\R$ is a $\K$-substack of $\F$ with inclusion $\rho:\R\ra\F$,
there is a 1-isomorphism $\R\cong\R\t_\F\R$ which by \eq{mi3eq2}
implies that $[(\R,\rho)]\cdot[(\R,\rho)]=[(\R,\rho)]$. Hence
\e
\io_\F(\de_{\oR(\K)}\de_{\oR(\K)})=\io_\F(\de_{\oR(\K)})=
\io_\F(\de_{\oR(\K)})\cdot\io_\F(\de_{\oR(\K)}).
\label{mi3eq4}
\e
When $\R,\fS$ are disjoint $\K$-substacks of $\F$ it is easy to see
that
\e
\io_\F(\de_{\oR(\K)}\de_{\oS(\K)})=0=
\io_\F(\de_{\oR(\K)})\cdot\io_\F(\de_{\oS(\K)}).
\label{mi3eq5}
\e
Given any $f,g\in\CF(\F)$ there exist a finite collection of
disjoint $\K$-substacks $\R_i$ of $\F$ such that $f,g$ are
$\Q$-linear combinations of the $\de_{\oR_i(\K)}$. Therefore
$\io_\F(fg)=\io_\F(f)\cdot\io_\F(g)$ follows from
\eq{mi3eq4}--\eq{mi3eq5} and bilinearity.

For $[(\R,\rho)],[(\fS,\si)]\in\SF(\F)$, apply \eq{mi2eq5} with
$\R\t_\F\fS,\fS,\R,\F$ in place of $\E,\F,\G,\H$ respectively to the
function $1_\fS\in\CF(\fS)$. This gives
\begin{equation*}
\CF^\stk(\pi_\R)1_{\R\t_\F\fS}\!=\!\CF^\stk(\pi_\R)\!\ci\!
\pi_\fS^*(1_\fS)\!=\!\rho^*\!\ci\!\CF^\stk(\si)1_\fS\!=\!
1_\R\cdot\rho^*\!\ci\!\CF^\stk(\si)1_\fS.
\end{equation*}
Applying $\CF^\stk(\rho)$ to this and using \eq{mi2eq3}, \eq{mi3eq2}
and \eq{mi3eq3} gives
\begin{gather*}
\pi_\F^\stk\bigl([(\R,\rho)]\cdot[(\fS,\si)]\bigr)=
\CF^\stk(\rho\ci\pi_\R)1_{\R\t_\F\fS}=\\
\CF^\stk(\rho)1_\R\cdot\CF^\stk(\si)1_\fS=
\pi_\F^\stk\bigl([(\R,\rho)])\pi_\F^\stk\bigl([(\fS,\si)]\bigr),
\end{gather*}
since multiplication by $\rho^*\ci\CF^\stk(\si)1_\fS$ and
$\CF^\stk(\si)1_\fS$ commute with $\CF^\stk(\rho)$. Thus
$\pi_\F^\stk(f\cdot g)=\pi_\F^\stk(f)\pi_\F^\stk(g)$ for
$f,g\in\SF(\F)$ follows by bilinearity.
\end{proof}

In general, $\io_\F$ is far from being surjective, and $\SF(\F),
\uSF(\F)$ are much larger than $\CF(\F)$. For example,
$\CF(\Spec\K)\cong\Q$, but one can show $\SF(\Spec\K)\cong
K_0(\Var_\K)\ot_\Z\Q$ and $\uSF(\Spec\K)\cong K_0({\rm
Sta}_\K)\ot_\Z\Q$, where $K_0(\Var_\K)$, $K_0({\rm Sta}_\K)$ are the
Grothendieck rings of the (2-)categories of $\K$-varieties and
algebraic $\K$-stacks respectively. The ring $K_0(\Var_\K)$ for
$\cha\K=0$ is studied by Bittner \cite{Bitt} and is clearly very
large, and $K_0({\rm Sta}_\K)$ is even larger. Also, $\pi^\na_\F$
does not usually commute with multiplication. Next we define {\it
pushforwards}, {\it pullbacks} and {\it tensor products} on stack
functions.

\begin{dfn} Let $\phi:\F\!\ra\!\G$ be a 1-morphism of algebraic
$\K$-stacks with affine geometric stabilizers. Define the {\it
pushforward\/} $\phi_*:\uSF(\F)\!\ra\!\uSF(\G)$~by
\begin{equation*}
\phi_*:\ts\sum_{i=1}^nc_i[(\R_i,\rho_i)]\longmapsto
\ts\sum_{i=1}^nc_i[(\R_i,\phi\ci\rho_i)].
\end{equation*}
This intertwines the relations \eq{mi3eq1} in $\uSF(\F),\uSF(\G)$,
and so is well-defined. If $\phi$ is {\it representable} then the
restriction maps $\phi_*:\SF(\F)\ra\SF(\G)$, since the
$\phi\ci\rho_i$ are representable as $\phi,\rho_i$ are.

Now let $\phi$ be of {\it finite type}. If $\R_i$ is a finite type
algebraic $\K$-stack and $\rho_i:\R_i\ra\G$ a 1-morphism then we may
form the Cartesian square:
\e
\begin{gathered}
\xymatrix@C=70pt@R=10pt{ \R_i\t_{\rho_i,\G,\phi}\F \ar[r]_{\pi_\F}
\ar[d]^{\,\pi_{\R_i}} & \F
\ar[d]_{\phi\,} \\
\R_i \ar[r]^{\rho_i} & \G. }
\end{gathered}
\label{mi3eq6}
\e
Since $\R_i$ and $\phi$ are of finite type, so are $\pi_{\R_i}$ and
$\R_i\t_{\rho_i,\G,\phi}\F$ as \eq{mi3eq6} is Cartesian. Define the
{\it pullback\/} $\phi^*:\uSF(\G)\!\ra\!\uSF(\F)$~by
\e
\phi^*:\ts\sum_{i=1}^nc_i[(\R_i,\rho_i)]\longmapsto
\ts\sum_{i=1}^nc_i[(\R_i\t_{\rho_i,\G,\phi}\F,\pi_\F)].
\label{mi3eq7}
\e
This is well-defined as $\R_i\t_{\rho_i,\G,\phi}\F$ is unique up to
1-isomorphism, and $\phi^*$ intertwines the relations \eq{mi3eq1} in
$\uSF(\G),\uSF(\F)$. The restriction maps
$\phi^*:\SF(\G)\ra\SF(\F)$, since the $\pi_\F$ are representable as
the $\rho_i$ are, and \eq{mi3eq6} is Cartesian. The {\it tensor
product\/} $\ot:\uSF(\F)\!\t\!\uSF(\G)\!\ra\!\uSF(\F\!\t\!\G)$ and
$\ot\!:\!\SF(\F)\!\t\!\SF(\G)\!\ra\!\SF(\F\!\t\!\G)$~is
\e
\bigl(\mathop{\ts\sum}\limits_{i\in I}c_i[(\R_i,\rho_i)]\bigr)\ot
\bigl(\mathop{\ts\sum}\limits_{j\in J}d_j[(\fS_j,\si_j)]\bigr)=
\mathop{\ts\sum}\limits_{i\in I,\; j\in
J}c_id_j[(\R_i\t\fS_j,\rho_i\t\si_j)],
\label{mi3eq8}
\e
for finite $I,J$. This is compatible with the relations, and so
well-defined. It is the analogue of the obvious map
$\ot:\CF(\F)\t\CF(\G)\ra\CF(\F\t\G)$ on constructible functions
given by~$(f\ot g)(x,y)=f(x)g(y)$.
\label{mi3def3}
\end{dfn}

We can now justify the name `stack function'. Each $[x]\in\oF(\K)$
is an isomorphism class of (finite type) 1-morphisms $x:\Spec\K
\ra\F$. These induce pullbacks $x^*:\uSF(\F)\ra\uSF(\Spec\K)$ and
$x^*:\SF(\F)\ra\SF(\Spec\K)$ depending only on $[x]$. Thus, to each
$f\in\uSF(\F)$ or $\SF(\F)$ we associate a function $\oF(\K)\ra\uSF
(\Spec\K)$ or $\SF(\Spec\K)$ by $[x]\mapsto x^*(f)$.

By definition $\uSF(\Spec\K)$ and $\SF(\Spec\K)$ are the $\Q$-vector
spaces generated by 1-isomorphism classes $[\R]$ of finite type
algebraic $\K$-stacks $\R$ with affine geometric stabilizers, and
finite type algebraic $\K$-spaces $\R$ respectively, with a relation
$[\R]=[\fS]+[\R\sm\fS]$ whenever $\fS$ is a closed $\K$-substack of
$\R$. Thus, stack functions on $\F$ are like `functions on $\oF(\K)$
with values in stacks'.

Here is the analogue of Theorem~\ref{mi2thm1}.

\begin{thm} Let\/ $\E,\F,\G,\H$ be algebraic $\K$-stacks with
affine geometric stabilizers and\/ $\be:\F\ra\G$, $\ga:\G\ra\H$ be
$1$-morphisms. Then
\e
\begin{aligned}
(\ga\!\ci\!\be)_*\!&=\!\ga_*\!\ci\!\be_*:\uSF(\F)\!\ra\!\uSF(\H),&
(\ga\!\ci\!\be)_*\!&=\!\ga_*\!\ci\!\be_*:\SF(\F)\!\ra\!\SF(\H),\\
(\ga\!\ci\!\be)^*\!&=\!\be^*\!\ci\!\ga^*:\uSF(\H)\!\ra\!\uSF(\F),&
(\ga\!\ci\!\be)^*\!&\!=\!\be^*\!\ci\!\ga^*:\SF(\H)\!\ra\!\SF(\F),
\end{aligned}
\label{mi3eq9}
\e
for $\be,\ga$ representable in the second equation, and of finite
type in the third and fourth. If\/ $f,g\in\uSF(\G)$ and\/ $\be$ is
finite type then $\be^*(f\cdot g)=\be^*(f)\cdot\be^*(g)$. If
\e
\begin{gathered}
\xymatrix@R=15pt{
\E \ar[r]_\eta \ar[d]^{\,\th} & \G \ar[d]_{\psi\,} \\
\F \ar[r]^\phi & \H }
\end{gathered}
\quad
\begin{gathered}
\text{is a Cartesian square with}\\
\text{$\th,\psi$ of finite type, then}\\
\text{the following commutes:}
\end{gathered}
\quad
\begin{gathered}
\xymatrix@C=35pt@R=10pt{
\uSF(\E) \ar[r]_{\eta_*} & \uSF(\G) \\
\uSF(\F) \ar[r]^{\phi_*} \ar[u]_{\,\th^*} & \uSF(\H).
\ar[u]^{\psi^*\,} }
\end{gathered}
\label{mi3eq10}
\e
The same applies for $\SF(\E),\ldots,\SF(\H)$ if\/ $\eta,\phi$ are
representable.
\label{mi3thm1}
\end{thm}

\begin{proof} The first and second equations of \eq{mi3eq9} follow from
\begin{equation*}
(\ga\ci\be)_*\bigl([(\R,\rho)]\bigr)=[(\R,\ga\ci\be\ci\rho)]=
\ga_*\bigl([(\R,\be\ci\rho)]\bigr)=\ga_*\ci\be_*\bigl([(\R,\rho)]\bigr),
\end{equation*}
as $(\ga\ci\be)\ci\rho=\ga\ci(\be\ci\rho)$. For the third and fourth
equations, we need to prove that for $\rho:\R\ra\H$ as in Definition
\ref{mi3def1} we have $(\ga\ci\be)^*\bigl([(\R,\rho)]\bigr)
=\be^*\ci\ga^*\bigl([(\R,\rho)]\bigr)$. This follows from the
existence of a 1-isomorphism
\e
\io:(\R\t_{\rho,\H,\ga}\G)\t_{\pi_\G,\G,\be}\F
\longra\R\t_{\rho,\H,\ga\ci\be}\F
\label{mi3eq11}
\e
with $\pi_\F\ci\io$ 2-isomorphic to $\pi_\F$, as 1-morphisms
$(\R\t_{\rho,\H,\ga}\G)\t_{\pi_\G,\G,\be}\F\ra\F$. We can construct
$\io$ easily using the explicit definition \cite[2.2.2]{LaMo} of
fibre products of $\K$-stacks. This proves \eq{mi3eq9}. One can show
$\be^*(f\cdot g)=\be^*(f)\cdot\be^*(g)$ for $f=[(\R,\rho)]$,
$g=[(\fS,\si)]$ using properties of fibre products.

For both cases of \eq{mi3eq10}, let $\rho:\R\ra\F$ be as in
Definition \ref{mi3def1}. Then
$\psi^*\ci\phi_*\bigl([(\R,\rho)]\bigr)=\eta_*\ci\th^*
\bigl([(\R,\rho)]\bigr)$ if
$[(\R\t_{\phi\ci\rho,\H,\psi}\G,\pi_\G)]=
[(\R\t_{\rho,\F,\th}\E,\eta\ci\pi_\E)]$. From Definition
\ref{mi2def3} and equivalence in Definition \ref{mi3def1} we see
that we may replace $\E$ here by $\F\t_{\phi,\H,\psi}\G$ and
$\th,\eta$ by $\pi_\F,\pi_\G$, so this is equivalent to
\begin{equation*}
\bigl[\bigl(\R\t_{\phi\ci\rho,\H,\psi}\G,\pi_\G\bigr)\bigr]=
\bigl[\bigl( \R\t_{\rho,\F,\pi_\F}(\F\t_{\phi,\H,\psi}\G),
\pi_\G\ci\pi_{\F\t_{\phi,\H,\psi}\G}\bigr)\bigr].
\end{equation*}
This follows from the existence of a 1-isomorphism
\begin{equation*}
\io':\R\t_{\rho,\F,\pi_\F}(\F\t_{\phi,\H,\psi}\G)\longra
\R\t_{\phi\ci\rho,\H,\psi}\G
\end{equation*}
with $\pi_\G\ci\io'$ 2-isomorphic to $\pi_\G\ci\pi_{\F\t_{\phi,
\H,\psi}\G}$, which can constructed using \cite[2.2.2]{LaMo} as for
$\io$ in \eq{mi3eq11} above. This completes the proof.
\end{proof}

Here are some compatibilities between $\ot$ and the other
operations. The proofs are all elementary.

\begin{prop} Let\/ $\phi:\E\!\ra\!\F$ and\/ $\psi:\G\ra\H$ be
$1$-morphisms of algebraic $\K$-stacks with affine geometric
stabilizers and $e,f,g,h$ lie in $\uSF(\E),\ldots,\uSF(\H)$. Then
$\phi_*(e)\ot\psi_*(g)=(\phi\ot\psi)_*(e\ot g)$, and\/
$\phi^*(f)\ot\psi^*(h)=(\phi\ot\psi)^*(f\ot h)$ when $\phi,\psi$ are
of finite type. Also $e\ot(f\ot g)=(e\ot f)\ot g\in\uSF
(\E\t\F\t\G)$. If\/ $\K$ has characteristic zero then
$\pi_{\E\t\F}^\na(e\ot f)=\pi_\E^\na(e)\ot\pi^\na_\F(f)$ in
$\CF(\E\t\F)$, and\/ $\pi_{\E\t\F}^\stk(e\ot
f)=\pi_\E^\stk(e)\ot\pi^\stk_\F(f)$ when $e,f$ lie in
$\SF(\E),\SF(\F)$.
\label{mi3prop2}
\end{prop}

The next two results consider the relationships between pushforwards
and pullbacks of stack and constructible functions, via the
maps~$\io_\F,\pi_\F^\na,\pi_\F^\stk$.

\begin{prop} Let\/ $\phi:\F\!\ra\!\G$ be a finite type $1$-morphism
of algebraic $\K$-stacks with affine geometric stabilizers.
Then~$\phi^*\!\ci\!\io_\G\!=\!\io_\F\!\ci\!\phi^*\!:\!
\CF(\G)\!\ra\!\SF(\F)$.
\label{mi3prop3}
\end{prop}

\begin{proof} Let $\R$ be a finite type $\K$-substack of $\G$ with
inclusion $\rho:\R\ra\G$. Then $\R\t_{\rho,\G,\phi}\F$ is a finite
type $\K$-substack of $\F$ with inclusion $\pi_\F$, and
$\phi^*(\oR(\K))=\overline{(\R\t_{\rho,\G,\phi}\F)}(\K)\subseteq
\oF(\K)$. Hence
\begin{align*}
\phi^*\ci\io_\G(\de_{\oR(\K)})&=\phi^*\bigl(\bde_{\oR(\K)})\bigr)=
\phi^*\bigl([(\R,\rho)]\bigr)=[(\R\t_{\rho,\G,\phi}\F,\pi_\F)]\\
&=\bde_{\overline{(\R\t_{\rho,\G,\phi}\F)}(\K)}=\io_\F(\de_{\phi^*
(\oR(\K))})=\io_\F\ci\phi^*(\de_{\oR(\K)}).
\end{align*}
As such $\de_{\oR(\K)}$ generate $\CF(\G)$, the proposition follows
by linearity.
\end{proof}

\begin{thm} Let\/ $\K$ have characteristic zero, $\F,\G$ be
algebraic $\K$-stacks with affine geometric stabilizers, and\/
$\phi:\F\ra\G$ a $1$-morphism. Then
\begin{itemize}
\setlength{\itemsep}{0pt}
\setlength{\parsep}{0pt}
\item[{\rm(a)}] $\pi^\na_\G\ci\phi_*=\CF^\na(\phi)\ci\pi_\F^\na:
\uSF(\F)\ra\CF(\G);$
\item[{\rm(b)}] $\pi^\stk_\G\ci\phi_*=\CF^\stk(\phi)\ci\pi_\F^\stk:
\SF(\F)\ra\CF(\G)$ if\/ $\phi$ is representable; and
\item[{\rm(c)}] $\pi^\stk_\F\ci\phi^*=\phi^*\ci\pi_\G^\stk:
\SF(\G)\ra\CF(\F)$ if\/ $\phi$ is of finite type.
\end{itemize}
\label{mi3thm2}
\end{thm}

\begin{proof} Let $[(\R,\rho)]\in\uSF(\F)$, so that $\rho:\R\ra\F$.
Then
\begin{align*}
\CF^\na(\phi)\ci\pi_\F^\na\bigl([(\R,\rho)]\bigr)&=
\CF^\na(\phi)\ci\CF^\na(\rho)1_\R=\CF^\na(\phi\ci\rho)1_\R\\
&=\pi_\G^\na\bigl([(\R,\phi\ci\rho)]\bigr)=
\pi_\G^\na\ci\phi^*\bigl([(\R,\rho)]\bigr),
\end{align*}
using Definitions \ref{mi3def2} and \ref{mi3def3} and equation
\eq{mi2eq2}. Part (a) follows by linearity. The proof of (b) is the
same, using $\CF^\stk,\pi_\F^\stk,\pi_\G^\stk$ and \eq{mi2eq3}. Let
$[(\R,\rho)]\in\SF(\G)$. Then by linearity, part (c) follows from
\begin{gather*}
\pi^\stk_\F\ci\phi^*\bigl([(\R,\rho)]\bigr)
=\pi^\stk_\F\bigl([(\R\t_{\rho,\G,\phi}\F,\pi_\F)]\bigr)
=\CF^\stk(\pi_\F)1_{\R\t_{\rho,\G,\phi}\F}\\
=\CF^\stk(\pi_\F)\ci(\pi_\R)^*1_\R =\phi^*\ci\CF^\stk(\rho)1_\R
=\phi^*\ci\pi_\G^\stk\bigl([(\R,\rho)]\bigr),
\end{gather*}
using Definitions \ref{mi3def2} and \ref{mi3def3} and equation
\eq{mi2eq5} applied to the Cartesian square
\begin{equation*}
\xymatrix@C=60pt@R=7pt{ \R\t_{\rho,\G,\phi}\F \ar[r]_{\pi_\F}
\ar[d]^{\,\pi_\R} & \F \ar[d]_{\phi\,} \\ \R \ar[r]^\rho & \G, }
\end{equation*}
with $\rho,\pi_\F$ representable and $\phi,\pi_\R$ of finite type.
\end{proof}

The other possible commutation relations are in general false. That
is, we expect $\phi_*\ci\io_\F\ne\io_\G\ci\CF^\na(\phi)$,
$\phi_*\ci\io_\F\ne\io_\G\ci\CF^\stk(\phi)$ and $\pi^\na_\F\ci
\phi^*\ne\phi^*\ci\pi_\G^\na$. This is why we use only the
$\pi_\F^\stk$ and not the $\pi_\F^\na$ in the applications of
\cite{Joyc3,Joyc4,Joyc5}, as the $\pi_\F^\stk$ commute with both
pushforwards and pullbacks, but the $\pi_\F^\na$ do not.

Suppose $\F$ is a $\K$-variety, $\K$-scheme or algebraic $\K$-space,
and $[(\R,\rho)]\in\SF(\F)$. Then $\rho:\R\ra\F$ is representable,
so $\R$ is a finite type algebraic $\K$-space. Thus $\R$ can be
written as the disjoint union of finitely many quasiprojective
$\K$-subvarieties $X_i$, and $[(\R,\rho)]=\sum_i[(X_i,\rho_i)]$.
Therefore $\SF(\F)$ is generated over $\Q$ by $[(X,\rho)]$ for $X$ a
quasiprojective $\K$-variety, with relations $[(X,\rho)]=[(Y,\rho)]+
[(X\sm Y,\rho)]$ for closed subvarieties~$Y\subseteq X$.

This implies that for $\F$ a $\K$-variety, $\SF(\F)$ equals
$K_0(\Var_\F)\ot_\Z\Q$, where $K_0(\Var_\F)$ is the {\it
Grothendieck group of\/ $\F$-varieties} studied by Bittner \cite[\S
5]{Bitt}. So this section generalizes the ideas of Bittner to Artin
stacks. The operations `$\,\cdot\,$',$\phi_*,\phi^*,\ot$ of
Definitions \ref{mi3def1} and \ref{mi3def3} agree with
$\ot,\phi_!,\phi^*,\boxtimes$ in~\cite[\S 6]{Bitt}.

This raises two interesting questions. Firstly, Bittner also defines
an involution ${\mathcal D}_\F$ and operations $\phi_*,\phi^!$ on a
modified space $K_0(\Var_\F)[[\bA^1]^{-1}]$. Do these have analogues
for Artin stacks? Secondly, modifications of $K_0(\Var_\F)$ are the
natural value groups for {\it motivic integrals}, which is the main
reason for studying them. Can the theory of motivic integration be
extended to Artin stacks, using modifications of our
spaces~$\uSF,\SF(\F)$?

Finally, we define {\it local stack functions}, the analogue of
locally constructible functions. Roughly speaking, we want to repeat
Definition \ref{mi3def1} using pairs $(\R,\rho)$ for which $\R$ is
not necessarily of finite type, but $\rho$ is. However, this must be
modified in two ways. Firstly, we allow sums $\sum_{i\in
I}c_i[(\R_i,\phi_i)]$ over {\it infinite} indexing sets $I$, because
locally constructible functions can take infinitely many values.
Secondly, the relations \eq{mi3eq1} are no longer sufficient,
because for $\R$ not of finite type we should be able to cut $\R$
into infinitely many disjoint pieces, but \eq{mi3eq1} allows only
for finite decompositions.

\begin{dfn} Let $\F$ be an algebraic $\K$-stack with affine
geometric stabilizers. Consider pairs $(\R,\rho)$, where $\R$ is an
algebraic $\K$-stack with affine geometric stabilizers and
$\rho:\R\ra\F$ is a finite type 1-morphism, with {\it equivalence}
of pairs as in Definition \ref{mi3def1}. Let $\underline{V\!}\,_\F$
be the $\Q$-vector space of formal $\Q$-linear combinations
$\sum_{i\in I}c_i[(\R_i,\rho_i)]$, where $I$ is a possibly infinite
indexing set, $c_i\in\Q$ and $[(\R_i,\rho_i)]$ is an equivalence
class as above, such that for all finite type $\K$-substacks $\G$ in
$\F$ with inclusion 1-morphism $\phi:\G\ra\F$, there are only
finitely many $i\in I$ with $c_i\ne 0$ and
$\R_i\t_{\rho_i,\F,\phi}\G$ nonempty.

Let $\underline{W\!}\,_\F$ be the vector subspace of $\sum_{i\in
I}c_i[(\R_i,\rho_i)]$ in $\underline{V\!}\,_\F$ such that for all
finite type $\K$-substacks $\G$ in $\F$ with inclusion 1-morphism
$\phi:\G\ra\F$, we have $\sum_{i\in I}
c_i[(\R_i\t_{\rho_i,\F,\phi}\G,\pi_\G)]=0$ in $\uSF(\G)$. There are
only finitely many nonzero terms in this sum by definition of
$\underline{V\!}\,_\F$, so this makes sense. Define $\uLSF(\F)$ to
be the quotient $\underline{V\!}\,_\F/\underline{W\!}\,_\F$. Define
$V_\F,W_\F,\LSF(\F)$ in exactly the same way, but with all
1-morphisms $\rho_i$ {\it representable}, and interpreting the
relation $\sum_{i\in I}c_i[(\R_i\t_{\rho_i,\F,\phi}\G, \pi_\G)]=0$
in~$\SF(\G)$.

We define commutative, associative multiplications `$\,\cdot\,$' on
$\uLSF(\F),\LSF(\F)$ by extending \eq{mi3eq2} bilinearly to sums
$\sum_{i\in I}c_i[(\R_i,\rho_i)]$, $\sum_{j\in
J}d_j[(\fS_j,\si_j)]$.
\label{mi3def4}
\end{dfn}

If $\F$ is of {\it finite type} and $\rho:\R\ra\F$ a 1-morphism then
$\R$ is of finite type if and only if $\rho$ is, and taking $\G=\F$
shows sums in $\underline{V\!}\,_\F$ have only finitely many nonzero
terms. It follows easily that $\uLSF(\F)=\uSF(\F)$ and
$\LSF(\F)=\SF(\F)$ in this case, just as $\LCF(\F)=\CF(\F)$. All the
definitions and results above for $\uSF(\F),\SF(\F)$ have
straightforward generalizations to $\uLSF(\F),\LSF(\F)$, analogous
to \cite[\S 5.3]{Joyc1}. We just state these, leaving the proofs as
an exercise. Note the differences in which 1-morphisms are required
to be of finite type.

\begin{dfn} Let $\F$ be an algebraic $\K$-stack with affine geometric
stabilizers and $S\subseteq\oF(\K)$ a locally constructible subset.
Then we may write $S=\coprod_{i\in I}\oR_i(\K)$, for $\K$-substacks
$\R_i$ of $\F$ with only finitely many intersecting any
constructible set $C\subseteq\oF(\K)$. Let $\rho_i:\R_i\ra\F$ be the
inclusion 1-morphism, which is representable and of finite type.
Define a local stack function
\begin{equation*}
\bde_S=\ts\sum_{i\in I}[(\R_i,\rho_i)]\in\LSF(\F)\subseteq\uLSF(\F).
\end{equation*}
This is independent of the choice of $I,\R_i$. Define $\io_\F:
\LCF(\F)\!\ra\!\LSF(\F)\!\subseteq\!\uLSF(\F)$ by $\io_\F(f)
\!=\!\ts\sum_{c\in f(\oF(\K))}c\cdot\bde_{f^{-1}(c)}$. This
potentially infinite sum makes sense as only finitely many terms are
nonzero over any constructible subset. For $\K$ of characteristic
zero, define $\Q$-linear maps $\pi_\F^\na:\uLSF(\F)\!\ra\!\LCF(\F)$
and $\pi_\F^\stk:\LSF(\F)\!\ra\!\LCF(\F)$~by
\begin{align*}
\pi_\F^\na\bigl(\ts\sum_{i\in I}c_i[(\R_i,\rho_i)]\bigr)&=
\ts\sum_{i\in I}c_i\LCF^\na(\rho_i)1_{\R_i}\\
\text{and}\qquad \pi_\F^\stk\bigl(\ts\sum_{i\in
I}c_i[(\R_i,\rho_i)]\bigr)&= \ts\sum_{i\in
I}c_i\LCF^\stk(\rho_i)1_{\R_i}.
\end{align*}
Here $\LCF^\na(\rho_i),\LCF^\stk(\rho_i)$ make sense as $\rho_i$ is
of finite type. On any constructible subset there are only finitely
many nonzero terms on the right hand sides of these equations, so
they are well-defined and lie in $\LCF(\F)$. The analogue of
Proposition \ref{mi3prop1} holds
for~$\io_\F,\pi_\F^\na,\pi_\F^\stk$.
\label{mi3def5}
\end{dfn}

\begin{dfn} Let $\phi:\F\!\ra\!\G$ be a {\it finite type} 1-morphism
of algebraic $\K$-stacks with affine geometric stabilizers. Define
$\phi_*:\uLSF(\F)\!\ra\!\uLSF(\G)$~by
\begin{equation*}
\phi_*:\ts\sum_{i\in I}c_i[(\R_i,\rho_i)]\longmapsto \ts\sum_{i\in
I}c_i[(\R_i,\phi\ci\rho_i)].
\end{equation*}
If $\phi$ is also {\it representable} define $\phi_*:\LSF(\F)
\ra\LSF(\G)$ the same way. For any $\phi:\F\ra\G$, define
$\phi^*:\uLSF(\G)\!\ra\!\uLSF(\F)$ and $\phi^*:\LSF(\G)
\!\ra\!\LSF(\F)$~by
\begin{equation*}
\phi^*:\ts\sum_{i\in I}c_i[(\R_i,\rho_i)]\longmapsto \ts\sum_{i\in
I}c_i[(\R_i\t_{\rho_i,\G,\phi}\F,\pi_\F)].
\end{equation*}
As in Proposition \ref{mi3prop3} we have $\phi^*\ci\io_\G=
\io_\F\ci\phi^*:\LCF(\G)\ra\LSF(\F)$. Define $\ot:\uLSF(\F)\!\t\!
\uLSF(\G)\!\ra\!\uLSF(\F\!\t\!\G)$ and $\ot\!:\!\LSF(\F)\!\t\!
\LSF(\G)\!\ra\!\LSF(\F\!\t\!\G)$ by \eq{mi3eq8}, allowing $I,J$
infinite.
\label{mi3def6}
\end{dfn}

\begin{thm} Let\/ $\E,\F,\G,\H$ be algebraic $\K$-stacks with
affine geometric stabilizers and\/ $\be:\F\ra\G$, $\ga:\G\ra\H$ be
$1$-morphisms. Then
\begin{align*}
(\ga\!\ci\!\be)_*\!&=\!\ga_*\!\ci\!\be_*:\uLSF(\F)\!\ra\!\uLSF(\H),&
(\ga\!\ci\!\be)_*\!&=\!\ga_*\!\ci\!\be_*:\LSF(\F)\!\ra\!\LSF(\H),\\
(\ga\!\ci\!\be)^*\!&=\!\be^*\!\ci\!\ga^*:\uLSF(\H)\!\ra\!\uLSF(\F),&
(\ga\!\ci\!\be)^*\!&\!=\!\be^*\!\ci\!\ga^*:\LSF(\H)\!\ra\!\LSF(\F),
\end{align*}
supposing $\be,\ga$ are of finite type in the first and second
equations, and representable in the second. If\/ $f,g\in\uLSF(\G)$
then $\be^*(f\cdot g)=\be^*(f)\cdot\be^*(g)$. If
\begin{equation*}
\begin{gathered}
\xymatrix@R=15pt{
\E \ar[r]_\eta \ar[d]^{\,\th} & \G \ar[d]_{\psi\,} \\
\F \ar[r]^\phi & \H }
\end{gathered}
\quad
\begin{gathered}
\text{is a Cartesian square with}\\
\text{$\eta,\phi$ of finite type, then}\\
\text{the following commutes:}
\end{gathered}
\quad
\begin{gathered}
\xymatrix@C=35pt@R=10pt{
\uLSF(\E) \ar[r]_{\eta_*} & \uLSF(\G) \\
\uLSF(\F) \ar[r]^{\phi_*} \ar[u]_{\,\th^*} & \uLSF(\H).
\ar[u]^{\psi^*\,} }
\end{gathered}
\end{equation*}
The same applies for $\LSF(\E),\ldots,\LSF(\H)$ if also $\eta,\phi$
are representable.
\label{mi3thm3}
\end{thm}

\begin{thm} Let\/ $\K$ have characteristic zero, $\F,\G$ be algebraic
$\K$-stacks with affine geometric stabilizers, and\/
$\phi:\F\!\ra\!\G$ a $1$-morphism. Then
\begin{itemize}
\setlength{\itemsep}{0pt}
\setlength{\parsep}{0pt}
\item[{\rm(a)}] $\pi^\na_\G\ci\phi_*\!=\!\LCF^\na(\phi)\!\ci\!\pi_\F^\na:
\uLSF(\F)\!\ra\!\LCF(\G)$ if\/ $\phi$ is of finite type;
\item[{\rm(b)}] $\pi^\stk_\G\ci\phi_*=\LCF^\stk(\phi)\ci\pi_\F^\stk:
\LSF(\F)\ra\LCF(\G)$ if\/ $\phi$ is representable and of finite
type; and
\item[{\rm(c)}] $\pi^\stk_\F\ci\phi^*=\phi^*\ci\pi_\G^\stk:
\LSF(\G)\ra\LCF(\F)$.
\end{itemize}
\label{mi3thm4}
\end{thm}

\section{Motivic invariants of stacks}
\label{mi4}

Let $\K$ be an algebraically closed field, and suppose $\Up$ is some
invariant of quasiprojective $\K$-varieties $X$ up to isomorphism,
taking values in a commutative ring or algebra $\La$. We call $\Up$
{\it motivic} if whenever $Y\subseteq X$ is a closed subvariety we
have $\Up([X])=\Up([X\sm Y])+\Up([Y])$, and whenever $X,Y$ are
varieties we have $\Up([X\t Y])=\Up([X]) \Up([Y])$. The name
`motivic' is a reference to the subject of {\it motives} and {\it
motivic integration}, where such constructions are common.
Well-known examples of motivic invariants are the {\it Euler
characteristic}, and {\it virtual Poincar\'e polynomials}.

This section extends such invariants $\Up$ to Artin stacks, in the
special case when $\ell-1$ and some other elements are invertible in
$\La$, where $\ell=\Up([\bA^1])$. Roughly speaking, we need this
because for a quotient stack $[X/G]$ we want to set
$\Up\bigl(\bigl[[X/G]\bigr]\bigr)=\Up([X])/\Up([G])$, but $\Up([G])$
is divisible by $\ell-1$ for any algebraic $\K$-group $G$ with
$\rk\,G>0$, so $(\ell-1)^{-1}$ must exist if $\Up([G])^{-1}$ does.

For virtual Poincar\'e polynomials we can make $\ell-1$ invertible
by defining $\La$ appropriately. But if $\Up$ is the Euler
characteristic $\chi$ then $\ell=\chi(\bA^1)=1$, so $\ell-1$ cannot
be invertible, and the approach of this section fails. Section
\ref{mi6} defines refined versions of the constructions of this
section, which do work when $\ell-1$ is not invertible, and so for
Euler characteristics.

Section \ref{mi41} explains the properties of $\Up$ we need and
gives examples, and \S\ref{mi42} explains how to extend $\Up$
naturally to $\Up'([\R])$ for finite type algebraic $\K$-stacks $\R$
with affine geometric stabilizers. This $\Up'$ is motivic and
satisfies $\Up\bigl(\bigl[[X/G]\bigr]\bigr)=\Up([X])/\Up([G])$ when
$G$ is a {\it special\/} algebraic $\K$-group. Section \ref{mi43}
combines these ideas with stack functions to define modified spaces
$\uSF(\F,\Up,\La)$ which will be powerful tools
in~\cite{Joyc3,Joyc4,Joyc5}.

\subsection{Initial assumptions and examples}
\label{mi41}

Here is the data we shall need for our constructions.

\begin{ass} Suppose $\La$ is a commutative $\Q$-algebra with
identity 1, and
\begin{equation*} \Up:\{\text{isomorphism classes
$[X]$ of quasiprojective $\K$-varieties $X$}\}\longra\La
\end{equation*}
a map, satisfying the following conditions:
\begin{itemize}
\setlength{\itemsep}{0pt} \setlength{\parsep}{0pt}
\item[(i)] If $Y\subseteq X$ is a closed subvariety then
$\Up([X])=\Up([X\sm Y])+\Up([Y])$;
\item[(ii)] If $X,Y$ are quasiprojective $\K$-varieties then
$\Up([X\!\t\!Y])\!=\!\Up([X])\Up([Y])$;
\item[(iii)] Write $\ell=\Up([\bA^1])$ in $\La$, where $\bA^1$ is
the affine line $\K$ regarded as a $\K$-variety. Then $\ell$ and
$\ell^k-1$ for $k=1,2,\ldots$ are invertible in~$\La$.
\end{itemize}
\label{mi4ass}
\end{ass}

We chose the notation `$\ell$' as in motivic integration $[\bA^1]$
is called the {\it Tate motive} and written $\mathbb L$. We will
often use following easy consequence of (i),(ii).

\begin{lem} Suppose Assumption \ref{mi4ass} holds, and\/ $\phi:X\ra
Y$ is a Zariski locally trivial fibration of quasiprojective
$\K$-varieties with fibre $F$, that is, $F$ is a quasiprojective
$\K$-variety and\/ $Y$ can be covered by Zariski open sets $U$ such
that\/ $\phi^{-1}(U)\cong F\t U$. Then~$\Up([X])=\Up([F])\Up([Y])$.
\label{mi4lem1}
\end{lem}

\begin{proof} Let $n\ge 0$ and $k\ge 1$ be given. Suppose by
induction that the lemma holds when either $\dim Y<n$ or $\dim Y=n$
and $Y$ has fewer than $k$ irreducible components. Let $\phi,X,Y$ be
as above, and suppose $\dim Y=n$ and $Y$ has $k$ irreducible
connected components. Then $Y\ne\emptyset$, so we can choose a
nonempty open set $U\subseteq Y$ with $\phi^{-1}(U)\cong F\t U$. Set
$Y'=Y\sm U$ and $X'=\phi^{-1}(Y')$. Then $X',Y'$ are quasiprojective
$\K$-varieties, and $\phi\vert_{X'}:X'\ra Y'$ is a Zariski locally
trivial fibration with fibre $F$, and either $\dim Y'<\dim Y$, or
$\dim Y'=\dim Y$ and $Y'$ has fewer irreducible components than $Y$.
So by the inductive hypothesis we have
$\Up([X'])=\Up([F])\Up([Y'])$. But then
\begin{align*}
\Up\bigl([X]\bigr)&=\Up\bigl([X']\bigr)+\Up\bigl([\phi^{-1}(U)]\bigr)
=\Up\bigl([F]\bigr)\Up\bigl([Y']\bigr)+\Up\bigl([F\t U]\bigr)\\
&=\Up\bigl([F]\bigr)\bigl(\Up\bigl([Y'])+\Up\bigl([U]\bigr)\bigr)=
\Up\bigl([F]\bigr)\Up\bigl([Y]\bigr),
\end{align*}
using Assumption \ref{mi4ass}(i),(ii). The lemma follows by
induction on~$n,k$.
\end{proof}

Here are some examples of suitable $\La,\Up$. The first, for
$\K=\C$, uses the {\it virtual Hodge polynomials} introduced by
Danilov and Khovanskii \cite[\S 1]{DaKh}, and discussed by
Cheah~\cite[\S 0.1]{Chea}.

\begin{ex} Let $\K=\C$. Define $\La_\Ho=\Q(x,y)$, the $\Q$-algebra
of rational functions in $x,y$ with coefficients in $\Q$. Elements
of $\La_\Ho$ are of the form $P(x,y)/Q(x,y)$, for $P,Q$ rational
polynomials in $x,y$ with~$Q\not\equiv 0$.

Let $X$ be a quasiprojective $\C$-variety of dimension $m$, and
$H^k_{\rm c}(X,\C)$ the {\it compactly-supported cohomology} of $X$.
Deligne defined a {\it mixed Hodge structure} on $H^k_{\rm
c}(X,\C)$. Let $h^{p,q}\bigl(H^k_{\rm c}(X,\C)\bigr)$ be the
corresponding {\it Hodge--Deligne numbers}. Following \cite[\S
1.5]{DaKh}, \cite[\S 0.1]{Chea} define the {\it virtual Hodge
polynomial\/} $e(X;x,y)$ of $X$ to be $e(X;x,y)\!=\!
\sum_{p,q=0}^m\sum_{k=0}^{2m}(-1)^kh^{p,q}\bigl(H^k_{\rm c}
(X,\C)\bigr)x^py^q$. Set $\Up_\Ho([X])=e(X;x,y)$. This lies in
$\Z[x,y]$, and so in $\La_\Ho$. Assumption \ref{mi4ass}(i),(ii) for
$\Up_\Ho$ follow from \cite[Prop.s 1.6 \& 1.8]{DaKh}, and
\cite[Ex.~1.10]{DaKh} gives $\ell=\Up_\Ho([\bA^1])=xy$, implying
(iii). Thus Assumption \ref{mi4ass} holds.
\label{mi4ex1}
\end{ex}

If $X$ is a smooth projective $\C$-variety then $h^{p,q}
\bigl(H^k_{\rm c}(X,\C)\bigr)=h^{p,q}(X)$ if $p+q=k$ and 0
otherwise, so $e(X;x,y)=\sum_{p,q=0}^m(-1)^{p+q}h^{p,q}(X)x^py^q$
just encodes the usual Hodge numbers of $X$. The point about virtual
Hodge polynomials is that they extend ordinary Hodge polynomials to
the non-smooth, non-projective case with the additive and
multiplicative properties we need.

As Hodge numbers refine Betti numbers, so the virtual Hodge
polynomial $e(X;x,y)$ refines the {\it virtual Poincar\'e
polynomial\/} $P(X;z)=e(X;-z,-z)$, as in Cheah \cite[\S 0.1]{Chea}.
However, virtual Poincar\'e polynomials work for all algebraically
closed $\K$, not just $\K=\C$. I have not been able to find a good
reference for the general $\K$ case, though some of the ideas can be
found in Deligne \cite{Deli}. I am grateful to Burt Totaro for
explaining it to me.

\begin{ex} Define $\La_\Po=\Q(z)$, the algebra of rational
functions in $z$ with coefficients in $\Q$. Let $\K=\C$ and $X$ be a
quasiprojective $\C$-variety. Deligne defined a {\it weight
filtration} on $H^k_{\rm c}(X,\C)$. Write $W_j\bigl(H^k_{\rm
c}(X,\C)\bigr)$ for the $j^{\rm th}$ quotient space of this
filtration. Define $P(X;z)=\sum_{j,k=0}^{2m}(-1)^{k-j} \dim
W_j\bigl(H^k_{\rm c}(X,\C)\bigr)z^j$ to be the {\it virtual
Poincar\'e polynomial\/} of $X$. Then $P(X;z)=e(X;-z,-z)$ and
$P(X;-1)=\chi(X)$, the {\it Euler characteristic} of $X$. Set
$\Up_\Po([X])=P(X;z)$. As in Example \ref{mi4ex1}, Assumption
\ref{mi4ass} holds for $\La_\Po,\Up_\Po$, with~$\ell=z^2$.

Here is how to extend this to general algebraically closed $\K$. If
$\K$ has characteristic zero and $X$ is a quasiprojective
$\K$-variety then $X$ is actually defined over a subfield $\K_0$ of
$\K$ which is finitely generated over $\Q$. That is,
$X=X_0\t_{\Spec\K_0}\Spec\K$, for $X_0$ a quasiprojective
$\K_0$-variety, and regarding $\Spec\K$ as a $\K_0$-scheme. We can
embed $\K_0$ as a subfield of $\C$, and form a quasiprojective
$\C$-variety $X_\C=X_0\t_{\Spec\K_0}\Spec\C$. Define
$P(X;z)=P(X_\C;z)$, reducing to the $\K=\C$ case, and
$\Up_\Po([X])=P(X;z)$. This is independent of choices, and
Assumption \ref{mi4ass} holds with~$\ell=z^2$.

If $\K$ has characteristic $p>0$ we use some different ideas,
sketched in Deligne \cite{Deli}. Write $\bF_p$ for the finite field
with $p$ elements, and $\bar\bF_p$ for its algebraic closure. Let
$l$ be a prime different from $p$. First we explain how to define
the virtual Poincar\'e polynomial of a quasiprojective
$\bF_p$-variety $X$. Then $X_{\bar\bF_p}
=X\t_{\Spec\bF_p}\Spec\bar\bF_p$ is a quasiprojective
$\bar\bF_p$-variety, so we can form the compactly-supported $l$-adic
cohomology $H^k_{\rm c}(X_{\bar\bF_p},\Q_l)$, a finite-dimensional
vector space over $\Q_l$.

The {\it geometric Frobenius\/} ${\rm Fr}$ acts on $X_{\bar\bF_p}$,
and so ${\rm Fr}^*$ acts on $H^k_{\rm c}(X_{\bar\bF_p},\Q_l)$. In
his proof of the Weil Conjecture, Deligne showed that the
eigenvalues of ${\rm Fr}^*$ are `Weil numbers of weight $j$' for
$j\ge 0$. Thus we may define a {\it weight filtration} on $H^k_{\rm
c}(X_{\bar\bF_p},\Q_l)$, whose $j^{\rm th}$ quotient space
$W^j\bigl(H^k_{\rm c}(X_{\bar\bF_p},\Q_l)\bigr)$ is the eigenspaces
of ${\rm Fr}^*$ with eigenvalues of weight $j$. Then we set
\begin{equation*}
\ts P(X;z)=\sum_{j,k=0}^{2m}(-1)^{k-j} \dim_{\Q_l}W_j\bigl(H^k_{\rm
c}(X_{\bar\bF_p},\Q_l)\bigr)z^j.
\end{equation*}

Now let $\K$ have characteristic $p>0$, and $X$ be a quasiprojective
$\K$-variety. Then $X$ is defined over a subfield $\K_0$ of $\K$
finitely generated over $\bF_p$, so $X=X_0\t_{\Spec\K_0}\Spec\K$,
for $X_0$ a quasiprojective $\K_0$-variety. Regard $X_0$ as a
$\bF_p$-variety with a dominant morphism $X_0\ra\Spec\K_0$, that is,
as a family of quasiprojective $\bF_p$-varieties. We specialize this
to get a quasiprojective $\bF_p$-variety $X_0^{\rm sp}$ --- that is,
$X_0^{\rm sp}$ is the fibre of
$X_0\t_{\Spec\bF_p}\Spec\bar\bF_p\ra\Spec\K_0\t_{\Spec\bF_p}
\bar\bF_p$ over a general point. Then we set $P(X;z)= P(X_0^{\rm
sp};z)$, reducing to the finite field case, and
$\Up_\Po([X])=P(X;z)$. Again, Assumption \ref{mi4ass} holds
with~$\ell=z^2$.
\label{mi4ex2}
\end{ex}

Here is the universal example, through which all other examples
factor.

\begin{ex} Let $\K$ be an algebraically closed field.
Define $\La_\uni$ to be the $\Q$-algebra generated by isomorphism
classes $[X]$ of quasiprojective $\K$-varieties $X$ and by
$\ell^{-1}$ and $(\ell^k-1)^{-1}$ for $k=1,2,\ldots$ for
$\ell=[\bA^1]$, with relations $[X]=[X\sm Y]+[Y]$ for $Y$ a closed
subvariety of $X$, and $[X\t Y]=[X]\cdot[Y]$ for $X,Y$
quasiprojective $\K$-varieties, and identity $[\Spec\K]$. Define
$\Up_\uni([X])=[X]$. Then Assumption \ref{mi4ass} holds trivially.
\label{mi4ex3}
\end{ex}

The drawback is that $\La_\uni$ is difficult to describe ---
Examples \ref{mi4ex1} and \ref{mi4ex2} are much more explicit. It is
a modification of the {\it Grothendieck group $K_0({\rm Var}_\K)$
of\/ $\K$-varieties}, as in Bittner \cite{Bitt}: $\La_\uni$ is
$(K_0({\rm Var}_\K)\ot_\Z\Q)[\ell^{-1},(\ell^k-1)^{-1}$,
$k=1,2,\ldots]$. Rings and algebras of this kind are often used in
{\it motivic integration}.

Notice that we have not included {\it Euler characteristics} in our
list of examples, though the Euler characteristic $\chi$ is the most
well-known and useful motivic invariant. This is because
$\Up([X])=\chi(X)$ does not satisfy Assumption \ref{mi4ass}, since
$\ell=\Up([\bA^1])=\chi(\bA^1)=1$, so $\ell^k-1=0$ is not invertible
in $\La=\Q$ for any $k=1,2,\ldots$, and Assumption \ref{mi4ass}(iii)
fails. Section \ref{mi6} will modify our approach for the case
$\ell=1$, to include Euler characteristics.

Assumption \ref{mi4ass} is equivalent to the hypotheses of Toen
\cite[Cor.~3.18]{Toen}. In \cite[\S 3.5]{Toen} he gives Example
\ref{mi4ex1} above and two new examples of suitable data $\Up,\La$,
{\it motivic Euler characteristics} and $l$-{\it adic Euler
characteristics}.

\subsection{Extending $\Up$ to a homomorphism $\uSF(\Spec\K)\ra\La$}
\label{mi42}

We now extend $\Up$ in \S\ref{mi41} from quasiprojective
$\K$-varieties to finite type $\K$-stacks with affine geometric
stabilizers. We express this as an {\it algebra homomorphism}
$\Up':\uSF(\Spec\K)\ra\La$. The next lemma is the reason for
Assumption~\ref{mi4ass}(iii).

\begin{lem} Let Assumption \ref{mi4ass} hold. Then for all\/
$m=1,2,\ldots$ we have
\e
\Up\bigl([\GL(m,\K)]\bigr)=\ell^{m(m-1)/2}\ts\prod_{k=1}^m(\ell^k-1),
\quad\text{which is invertible in $\La$.}
\label{mi4eq1}
\e
\label{mi4lem2}
\end{lem}

\begin{proof} Consider the projection morphism $\GL(m,\K)\ra
\bA^m\sm\{0\}$ taking a matrix to its first column. This is a
Zariski locally trivial fibration, with fibre
$\bA^{m-1}\t\GL(m-1,\K)$. So Lemma \ref{mi4lem1} gives
\begin{align*}
\Up\bigl([\GL(m,\K)]\bigr)&=\Up\bigl([\bA^m\sm\{0\}]\bigr)\cdot
\Up\bigl(\bA^{m-1}\bigr)\cdot\Up\bigl(\GL(m-1,\K)\bigr)\\
&=(\ell^m-1)\ell^{m-1}\cdot\Up\bigl(\GL(m-1,\K)\bigr).
\end{align*}
We deduce \eq{mi4eq1} by induction on $m$, and invertibility by
Assumption~\ref{mi4ass}(iii).
\end{proof}

\begin{lem} Let Assumption \ref{mi4ass} hold, and\/ $G$ be a special
algebraic $\K$-group. Then $\Up([G])$ is invertible in~$\La$.
\label{mi4lem3}
\end{lem}

\begin{proof} Embed $G\subseteq\GL(m,\K)$ with $\GL(m,\K)\ra
\GL(m,\K)/G$ a Zariski locally trivial fibration with fibre $G$.
Then $\Up([\GL(m,\K)])=\Up([G])\Up([\GL(m,\K)/G])$ by Lemma
\ref{mi4lem1}. But $\Up([\GL(m,\K)])$ is invertible by Lemma
\ref{mi4lem2}, so $\Up([G])$ is invertible in~$\La$.
\end{proof}

If a $\K$-stack $\R$ is 1-isomorphic to $[X/G]$ for $X$ a
quasiprojective $\K$-variety and $G$ a special algebraic $\K$-group,
we intend to define $\Up'([\R])=\Up([X])/\Up([G])$ in $\La$. This is
independent of choices.

\begin{prop} Let Assumption \ref{mi4ass} hold, and\/ $\R$
be a finite type algebraic $\K$-stack. Suppose $\R$ is
$1$-isomorphic to a quotient stack\/ $[X/G]$, where $X$ is a
quasiprojective $\K$-variety, and\/ $G$ is a special algebraic
$\K$-group acting on $X$. Then $\Up([X])/\Up([G])$ depends only on
$\R$, not on the choice of\/~$X,G$.
\label{mi4prop1}
\end{prop}

\begin{proof} Suppose $\R$ is 1-isomorphic to $[X/G]$ and
$[Y/H]$, for $X,Y$ quasiprojective $\K$-varieties and $G,H$ special
algebraic $\K$-groups acting on $X,Y$. The 1-isomorphisms
$[X/G]\cong\R\cong[Y/H]$ give 1-morphisms $\phi:X\ra\R$,
$\psi:Y\ra\R$ which are atlases, invariant under the $G,H$ actions.
Form the fibre product $Z=X\t_{\phi,\R,\psi}Y$. This is a finite
type algebraic $\K$-space with an action of $G\t H$, such that the
$G$- and $H$-actions are free with $Z/H\cong X$ and~$Z/G\cong Y$.

The projections $Z\ra X$ and $Z\ra Y$ are principal $H$- and
$G$-bundles respectively. Therefore $Z$ is a quasiprojective
$\K$-variety, not just an algebraic $\K$-space, as $H,X$ are
quasiprojective. Also $Z\ra X$ and $Z\ra Y$ are Zariski locally
trivial fibrations by Definition \ref{mi2def1}, since $G,H$ are
special. So by Lemma \ref{mi4lem1} we have
$\Up([Z])=\Up([H])\Up([X])=\Up([G])\Up([Y])$. Dividing by
$\Up([G])\Up([H])$, which is invertible by Lemma \ref{mi4lem3},
proves that~$\Up([X])/\Up([G])=\Up([Y])/\Up([H])$.
\end{proof}

We show by example that the condition that $G$ is special (or
something like it) is necessary in Proposition~\ref{mi4prop1}.

\begin{ex} Let $\K=\C$ and $\La_\Ho,\Up_\Ho$ be as in Example
\ref{mi4ex1}. Take $X$ to be the quasiprojective variety
$\C\sm\{0\}$, with affine $\C$-groups $\{1\}$ and $\{\pm 1\}$ acting
freely on $X$ by $\ep:x\mapsto\ep x$. Take $\R$ to be $X$, regarded
as a $\C$-stack. Then $\R$ is 1-isomorphic to $[X/\{1\}]$, and also
to $[X/\{\pm 1\}]$ by $\{\pm x\}\mapsto x^2$. We have
$\Up_\Ho([X])=xy-1$, $\Up_\Ho([\{\pm 1\}])=2$ and
$\Up_\Ho([\{1\}])=1$. Thus
\begin{equation*}
\Up_\Ho([X])/\Up_\Ho([\{1\}])=xy-1\ne\ha(xy-1)=
\Up_\Ho([X])/\Up_\Ho([\{\pm 1\}]),
\end{equation*}
so $\Up_\Ho([X])/\Up_\Ho([G])$ depends on the choice of $X,G$ with
$\R\!\cong\![X/G]$. This does not contradict Proposition
\ref{mi4prop1}, since $\{\pm 1\}$ is not special. Note also that
$x\mapsto x^2$ is a principal $\{\pm 1\}$-bundle which is not a
Zariski locally trivial fibration.
\label{mi4ex4}
\end{ex}

By definition $\uSF(\Spec\K)$ is spanned by $[(\R,\rho)]$ for
$\rho:\R\ra\Spec\K$ a 1-morphism. Since for any $\R$ there is a
projection $\rho:\R\ra\Spec\K$ unique up to 2-isomorphism, we omit
$\rho$, and write $[\R]$ instead of $[(\R,\rho)]$. Recall from
\S\ref{mi3} that $\uSF(\Spec\K)$ is a {\it commutative
$\Q$-algebra}. We shall construct an {\it algebra morphism} $\Up':
\uSF(\Spec\K)\ra\La$.

\begin{thm} Let Assumption \ref{mi4ass} hold. Then there exists
a unique morphism of\/ $\Q$-algebras $\Up':\uSF(\Spec\K)\ra\La$ such
that if\/ $G$ is a special algebraic $\K$-group acting on a
quasiprojective $\K$-variety $X$ then $\Up'\bigl(\bigl[[X/G]
\bigr]\bigr)=\Up([X])/\Up([G])$, where $\Up([G])$ is invertible by
Lemma~\ref{mi4lem3}.
\label{mi4thm1}
\end{thm}

\begin{proof} By linearity it is enough to define $\Up'([\R])$
for $[\R]\in\uSF(\Spec\K)$. Then $\R$ is a finite type algebraic
$\K$-stack with affine geometric stabilizers. Thus by Kresch
\cite[Prop.~3.5.9]{Kres} $\R$ can be {\it stratified by global
quotient stacks}. This means that the {\it associated reduced
stack\/} $\R^{\rm red}$ is the disjoint union of finitely many
locally closed substacks $\U_i$ for $i\in I$ with each $\U_i$
1-isomorphic to a global quotient $[X_i/G_i]$, with $X_i$ a
quasiprojective $\K$-variety and $G_i$ an affine $\K$-group acting
on $X_i$. (Kresch takes the $X_i$ to be $\K$-schemes, but using
varieties is equivalent.) As in \cite[Lem.~3.5.1]{Kres} we can take
$G_i=\GL(m_i,\K)$, so in particular we can suppose $G_i$ is special

Since $\R^{\rm red}$ is a closed $\K$-substack of $\R$ with
$\R\sm\R^{\rm red}$ empty we have
\begin{equation*}
[\R]=[\R^{\rm red}]=\ts\sum_{i\in I}[\U_i]=\ts\sum_{i\in I}
\bigl[[X_i/G_i]\bigr]
\end{equation*}
in $\uSF(\Spec\K)$. Thus, if $\Up'$ exists at all we must have
\e
\Up'\bigl([\R]\bigr)=\ts\sum_{i\in I}\Up\bigl([X_i]\bigr)/
\Up\bigl([G_i]\bigr).
\label{mi4eq2}
\e
This proves uniqueness of $\Up'$, if it exists. To show it does,
suppose $\R^{\rm red}$ is also the disjoint union of finitely many
locally closed substacks $\V_j$ for $j\in J$ with $\V_j$
1-isomorphic to $[Y_j/H_j]$, with $H_j$ a special algebraic
$\K$-group acting on a quasiprojective $\K$-variety~$Y_j$.

Since $\U_i$ is the disjoint union of locally closed $\K$-substacks
$\U_i\cap\V_j$ for $j\in J$, and $\U_i\cong[X_i/G_i]$, we can write
$X_i$ as the disjoint union of locally closed $G_i$-invariant
quasiprojective $\K$-subvarieties $X_{ij}$ for $j\in J$, with
$\U_i\cap\V_j\cong [X_{ij}/G_i]$. Similarly, we write $Y_j$ as the
disjoint union of locally closed, $H_j$-invariant $\K$-subvarieties
$Y_{ij}$ for $i\in I$, with $\U_i\cap\V_j\cong[Y_{ij}/H_j]$. Thus
$[X_{ij}/G_i]\cong[Y_{ij}/H_j]$, so $\Up([X_{ij}])/\Up([G_i])=
\Up([Y_{ij}])/\Up([H_j])$ by Proposition \ref{mi4prop1}. Therefore
\begin{align*}
\ts\sum_{i\in I}\Up\bigl([X_i]\bigr)/\Up\bigl([G_i]\bigr)&=
\ts\sum_{i\in I}\sum_{j\in J}\Up\bigl([X_{ij}]\bigr)/
\Up\bigl([G_i]\bigr)=\\
\ts\sum_{j\in J}\sum_{i\in I}\Up([Y_{ij}])/\Up([H_j])&=
\ts\sum_{j\in J}\Up([Y_j])/\Up([H_j]),
\end{align*}
by Assumption \ref{mi4ass}(i). Thus the right hand side of
\eq{mi4eq2} is independent of choices, and we can take \eq{mi4eq2}
as the {\it definition} of~$\Up'([\R])$.

Using Assumption \ref{mi4ass}(i) we find that if $\fS$ is a closed
$\K$-substack in $\R$ then $\Up'([\R])=\Up'([\fS])+
\Up'([\R\sm\fS])$. So $\Up'$ is compatible with the relations
\eq{mi3eq1} defining $\uSF(\Spec\K)$, and extends uniquely to a
$\Q$-linear map $\Up':\uSF(\Spec\K)\ra\La$. By Assumption
\ref{mi4ass}(ii) we see this is a $\Q$-algebra morphism. If $X,G$
are as in the theorem then taking $\R=[X/G]$, $I=\{1\}$, $X_1=X$,
$G_1=G$ in the definition \eq{mi4eq2} gives $\Up'
\bigl(\bigl[[X/G]\bigr]\bigr)=\Up([X])/\Up([G])$, as we want.
\end{proof}

Theorem \ref{mi4thm1} is similar to Toen \cite[Th.~1.1]{Toen}.
Combining it with Examples \ref{mi4ex1} and \ref{mi4ex2}, for $\R$ a
finite type algebraic $\K$-stack with affine geometric stabilizers,
we can define the {\it virtual Hodge function}
$e(\R;x,y)=\Up'_\Ho([\R])$ when $\K=\C$, and the {\it virtual
Poincar\'e function} $P(\R;z)=\Up'_\Po([\R])$ for general $\K$.
These are $\Q$-rational functions in $x,y$ and $z$ respectively,
agree with the usual virtual Hodge and Poincar\'e polynomials when
$\R$ is a quasiprojective $\K$-variety, and have attractive additive
and multiplicative properties.

Let $G$ be an affine algebraic $\K$-group, and take $\R=
[\Spec\K/G]$, which is a single point $r$ with $\Aut_\K(r)=G$. Then
the theorem gives $\Up'([\R])=\Up([\GL(m,\K)/G])/\Up([\GL(m,\K)])$
for any embedding $G\subseteq\GL(m,\K)$. If $G$ is special this
reduces to $\Up'([\R])=\Up([G])^{-1}$, but in general this is false
--- when $\K=\C$ and $G=\{\pm 1\}$, Example \ref{mi4ex4} gives
$\Up_\Ho'([\R])=1$ but~$\Up_\Ho([G])^{-1}=\ha$.

This has surprising implications. In problems involving `counting'
points on Deligne--Mumford stacks, say, one would expect a point $r$
with (finite) stabilizer group $G$ to `count' with weight
$1/\md{G}$. But our discussion shows that $\Up'$ `counts' points
with stabilizer $G$ with weight $\Up([\GL(m,\K)/G])/\Up
([\GL(m,\K)])$, which is not $1/\md{G}$ in general. So these ideas,
especially the r\^ole of special algebraic $\K$-groups, may be
telling us about the `right' way to approach counting problems on
stacks, such as the invariants studied in~\cite{Joyc5}.

\subsection{The spaces $\uSF(\F,\Up,\La)$ and their operations}
\label{mi43}

We now integrate the ideas of \S\ref{mi41}--\S\ref{mi42} with the
stack function material of \S\ref{mi3}. Here is an extension of
Definition~\ref{mi3def1}.

\begin{dfn} Let Assumption \ref{mi4ass} hold, and $\F$ be an
algebraic $\K$-stack with affine geometric stabilizers. Consider
pairs $(\R,\rho)$, where $\R$ is a finite type algebraic $\K$-stack
with affine geometric stabilizers and $\rho:\R\ra\F$ is a
1-morphism, with {\it equivalence} of pairs as in Definition
\ref{mi3def1}. Define $\uSF(\F,\Up,\La)$ to be the $\La$-module
generated by equivalence classes $[(\R,\rho)]$ as above, with the
following {\it relations}:
\begin{itemize}
\setlength{\itemsep}{0pt}
\setlength{\parsep}{0pt}
\item[(i)] Given $[(\R,\rho)]$ as above and $\fS$ a closed $\K$-substack
of $\R$ we have $[(\R,\rho)]=[(\fS,\rho\vert_\fS)]+[(\R\sm\fS,
\rho\vert_{\R\sm\fS})]$, as in~\eq{mi3eq1}.
\item[(ii)] Let $\R$ be a finite type algebraic $\K$-stack with
affine geometric stabilizers, $U$ a quasiprojective $\K$-variety,
$\pi_\R:\R\t U\ra\R$ the natural projection, and $\rho:\R\ra\F$ a
1-morphism. Then~$[(\R\t U,\rho\ci\pi_\R)]=\Up([U])[(\R,\rho)]$.
\item[(iii)] Given $[(\R,\rho)]$ as above and a 1-isomorphism
$\R\cong[X/G]$ for $X$ a quasiprojective $\K$-variety and $G$ a
special algebraic $\K$-group acting on $X$, we have
$[(\R,\rho)]=\Up([G])^{-1}[(X,\rho\ci\pi)]$, where
$\pi:X\ra\R\cong[X/G]$ is the natural projection 1-morphism. Here
$\Up([G])^{-1}$ exists in $\La$ by Lemma~\ref{mi4lem3}.
\end{itemize}

Similarly, we could define $\SF(\F,\Up,\La)$ to be the $\La$-module
generated by $[(\R,\rho)]$ with $\rho$ {\it representable}, and
relations (i)--(iii) as above. But using (i),(iii) above we find
$\uSF(\F,\Up,\La)$ is spanned over $\La$ by $[(X,\rho\ci\pi)]$ for
$\K$-varieties $X$. Then $\rho\ci\pi$ is automatically
representable, so $[(X,\rho\ci\pi)]\in\SF(\F,\Up, \La)$, and
$\SF(\F,\Up,\La)=\uSF(\F,\Up,\La)$. Thus we shall not bother with
the~$\SF(\F,\Up,\La)$.

Define a $\Q$-linear projection $\Pi^{\Up,\La}_\F:\uSF(\F)\ra
\uSF(\F,\Up,\La)$ by
\e
\Pi^{\Up,\La}_\F:\ts\sum_{i\in I}c_i[(\R_i,\rho_i)]\longmapsto
\ts\sum_{i\in I}c_i[(\R_i,\rho_i)],
\label{mi4eq3}
\e
using the embedding $\Q\subseteq\La$ to regard $c_i\in\Q$ as an
element of $\La$. Then $\Pi^{\Up,\La}_\F$ is well-defined, as the
relation \eq{mi3eq1} in $\uSF(\F)$ maps to relation (i) above.
\label{mi4def1}
\end{dfn}

The important point here is the {\it relations} (i)--(iii) above.
These are not arbitrary, but lead to interesting spaces, as our
results below will show. In defining a space by generators and
relations, one should consider two issues. The first is that any
{\it operations} on the spaces we define by their action on
generators must be compatible with all the relations, or they will
not be well-defined. We deal with this in Theorem \ref{mi4thm2}
below.

The second is that if we impose too many relations, or inconsistent
relations, then the space may be much smaller than we expect, even
zero. We will show in Proposition \ref{mi4prop3} below that
$\uSF(\F,\Up,\La)$ is at least as large as $\CF(\F)\ot_\Q\La$. So
the spaces $\uSF(\F,\Up,\La)$ are quite large (though much smaller
than $\uSF(\F)\ot_\Q\La$), and (i)--(iii) have some kind of
consistency about them. As in \S\ref{mi3}, $\uSF(\F,\Up,\La)$ has
{\it multiplication, pushforwards, pullbacks} and {\it tensor
products}.

\begin{dfn} Let Assumption \ref{mi4ass} hold, $\F,\G$ be
algebraic $\K$-stacks with affine geometric stabilizers, and
$\phi:\F\ra\G$ a 1-morphism. Define a $\La$-bilinear {\it
multiplication} `$\,\cdot\,$' on $\uSF(\F,\Up,\La)$ by \eq{mi3eq2}.
This is {\it commutative} and {\it associative} as in Definition
\ref{mi3def1}. Define the {\it pushforward\/} $\phi_*:
\uSF(\F,\Up,\La)\ra\uSF(\G,\Up,\La)$ by \eq{mi3eq6}, taking the
$c_i\in\La$ rather than $c_i\in\Q$. For $\phi$ of {\it finite type},
define the {\it pullback\/}
$\phi^*:\uSF(\G,\Up,\La)\!\ra\!\uSF(\F,\Up,\La)$ by \eq{mi3eq7}.
Define the {\it tensor product\/}
$\ot:\uSF(\F,\Up,\La)\t\uSF(\G,\Up,\La)\ra\uSF(\F\t\G,\Up,\La)$
by~\eq{mi3eq8}.
\label{mi4def2}
\end{dfn}

Notice that we do not define $\pi^\stk_\F$ for the
$\uSF(\F,\Up,\La)$, as in \S\ref{mi3}. This is because $\pi^\stk_\F$
for $\SF(\F)$ is defined using the Euler characteristic $\chi$, and
to define their analogues for $\uSF(\F,\Up,\La)$ we would need an
algebra morphism $\Phi:\La\ra\Q$ with
$\chi(X)\equiv\Phi\ci\Up([X])$. But $\Up([\bG_m])$ is invertible in
$\La$ and $\chi([\bG_m])\!=\!0$, so no such $\Phi$ exists. The
analogue of $\io_\F$ for $\uSF(\F,\Up,\La)$
is~$\Pi^{\Up,\La}_\F\!\ci\!\io_\F$.

\begin{thm} These operations `$\,\cdot\,$', $\phi_*$, $\phi^*$
and\/ $\ot$ are compatible with the relations {\rm(i)--(iii)} in
Definition \ref{mi4def1}, and so are well-defined.
\label{mi4thm2}
\end{thm}

\begin{proof} Regard `$\,\cdot\,$'$,\phi_*,\phi^*,\ot$ as defined
on generators $[(\R,\rho)],[(\fS,\si)]$ by \eq{mi3eq2}, \eq{mi3eq6},
\eq{mi3eq7}, \eq{mi3eq8} giving well-defined elements of
$\uSF(*,\Up,\La)$. We have to show that applying
`$\,\cdot\,$'$,\phi_*,\phi^*$ or $\ot$ to each relation (i)--(iii)
above gives a finite $\La$-linear combination of relations
(i)--(iii), that is, relations map to relations. All four are
compatible with (i), as for the $\uSF(\F)$ case in \S\ref{mi3}. For
$\phi_*$ and $\ot$ compatibility with (ii)--(iii) is easy. So we
must show `$\,\cdot\,$', $\phi^*$ are compatible with (ii)--(iii).

For `$\,\cdot\,$' and $\phi^*$, compatibility with (ii) follows as
the factor $U$ passes through the appropriate fibre products. So,
for instance, we have
\begin{align*}
[(\R\t U,\rho\ci\pi_\R)]\cdot[(\fS,\si)]&=[((\R\t
U)\t_{\rho\ci\pi_\R,\F,\si}\fS,\rho\ci\pi_{\R\t U})]\\
&=[(\R\t_{\rho,\F,\si}\fS)\t U,\rho\ci\pi_\R\ci\pi_{\R\t_\F\fS})].
\end{align*}
Therefore right multiplication `$\,\cdot[(\fS,\si)]$' maps (ii) to
(ii), and left multiplication does too by commutativity, so
`$\,\cdot\,$' is compatible with (ii). A similar argument works for
$\phi^*$ and~(ii).

Let $[(\R,\rho)],[(\fS,\si)]\in\uSF(\F,\Up,\La)$, with
$\R\cong[X/G]$ for $X$ a quasiprojective $\K$-variety acted on by a
special algebraic $\K$-group $G$. Using Kresch
\cite[Prop.~3.5.9]{Kres} as in Theorem \ref{mi4thm1}, we can find
finite sets $I,J$ and $\K$-substacks $\R_i,\fS_j,\F_{ij}$ in
$\R,\fS,\F$ for all $i\in I$ and $j\in J$, such that
$\R=\coprod_{i\in I}\R_i$, $\fS=\coprod_{j\in J}\fS_j$, and
$\rho,\si$ map $\rho:\R_i\ra\F_{ij}$, $\si:\fS_j\ra\F_{ij}$, and
$\R_i\cong[X_i/G]$ for $X_i$ a $G$-invariant subvariety of $X$ with
$X=\coprod_{i\in I}X_i$, and $\fS_j\cong[Y_j/H_j]$,
$\F_{ij}\cong[Z_{ij}/K_{ij}]$ for quasiprojective $\K$-varieties
$Y_j,Z_{ij}$ acted on by special algebraic $\K$-groups~$H_j,K_{ij}$.

Refining the decompositions if necessary, we can suppose the
1-morphisms $[X_i/G]\ra[Z_{ij}/K_{ij}]$ and $[Y_j/H_j]\ra
[Z_{ij}/K_{ij}]$ corresponding to $\rho:\R_i\ra\F_{ij}$ and
$\si:\fS_j\ra\F_{ij}$ are induced by $\K$-variety morphisms
$\al_{ij}:X_i\ra Z_{ij}$, $\be_{ij}:Y_j\ra Z_{ij}$ equivariant with
respect to $\K$-group morphisms $\ga_{ij}:G\ra K_{ij}$ and
$\de_{ij}:H_j\ra K_{ij}$. By \eq{mi3eq2} and (i) we see that in
$\uSF(\F,\Up,\La)$ we have
\e
[(\R,\rho)]\cdot[(\fS,\si)]=\ts\sum_{i\in I,\; j\in J}
\bigl[\bigl([X_i/G]\t_{[Z_{ij}/K_{ij}]}[Y_j/H_j],
\rho\ci\pi_{\R_i})].
\label{mi4eq4}
\e

The definitions of fibre products and quotients yield a
1-isomorphism
\e
\begin{split}
&[X_i/G]\t_{[Z_{ij}/K_{ij}]}[Y_j/H_j]\cong \\
&\bigl[\bigl((X_i\t Y_j)\t_{\al_{ij}\t\be_{ij},Z_{ij}\t
Z_{ij},\pi_{ij}}(Z_{ij}\t K_{ij})\bigr)/G\t H_{ij}\bigr],
\label{mi4eq5}
\end{split}
\e
using the fibre product of $\K$-varieties $X_i\t Y_j$ and $Z_{ij}\t
K_{ij}$ over $Z_{ij}\t Z_{ij}$, where $\pi_{ij}:Z_{ij}\t K_{ij}\ra
Z_{ij}\t Z_{ij}$ is the $\K$-variety morphism
$\pi_{ij}:(z,k)\mapsto(z,k\cdot z)$. Here $G\t H_{ij}$ acts on
$X_i\t Y_j$ by $(g,h):(x,y)\mapsto(g\cdot x,h\cdot y)$, on $Z_{ij}\t
K_{ij}$ by $(g,h):(z,k)\mapsto\bigl(\ga_{ij}(g)\cdot z,
\de_{ij}(h)k\ga_{ij}(g^{-1})\bigr)$, and on $Z_{ij}\t Z_{ij}$ by
$(g,h):(z,z')\mapsto(\ga_{ij}(g)\cdot z,\de_{ij}(h)\cdot z')$. These
actions commute with $\al_{ij}\t\be_{ij}$, $\pi_{ij}$, and so push
down to the fibre product~$(X_i\t Y_j)\t_{Z_{ij}\t Z_{ij}}(Z_{ij}\t
K_{ij})$.

Using \eq{mi4eq4} and the compatibility of `$\,\cdot\,$' with (i),
it is enough to show that
\begin{align*}
&\bigl[\bigl(\bigl[\bigl((X_i\t Y_j)\t_{\al_{ij}\t\be_{ij},Z_{ij}\t
Z_{ij},\pi_{ij}}(Z_{ij}\t
K_{ij})\bigr)/G\t H_{ij}\bigr],\rho\ci\pi_{\R_i}\bigr)\bigr]=\\
&\Up([G])^{-1}\bigl[\bigl(\bigl[\bigl((X_i\!\t\! Y_j)
\t_{\al_{ij}\!\t\!\be_{ij},Z_{ij}\!\t\!Z_{ij},\pi_{ij}}
(Z_{ij}\!\t\!K_{ij})\bigr)/\{1\}\!\t\!H_{ij}\bigr],
\rho\!\ci\!\pi_{\R_i}\!\ci\!\pi_{ij}\bigr)\bigr]
\end{align*}
in $\uSF(\F,\Up,\La)$. This holds as by (iii) both sides are equal
to
\begin{equation*}
\Up([G\t H_{ij}])^{-1}\bigl[\bigl((X_i\t
Y_j)\t_{\al_{ij}\t\be_{ij},Z_{ij}\t Z_{ij},\pi_{ij}}(Z_{ij}\t
K_{ij}),\rho\ci\pi_{\R_i}\ci\pi_{ij}'\bigr)\bigr],
\end{equation*}
showing right multiplication `$\,\cdot[(\fS,\si)]$' is compatible
with (iii). Left multiplication is too, so `$\,\cdot\,$' is
compatible with (i)--(iii), and is well-defined.

To show $\phi^*$ is compatible with (iii), let $\phi:\F\ra\G$ be of
finite type and $[(\R,\rho)]\in \uSF(\G,\Up,\La)$ with
$\R\cong[X/G]$ as usual. Since $\R$ is of finite type its image is
constructible in $\G$, so we can find a finite collection of
disjoint finite type $\K$-substacks $\G_i$ in $\G$ such that
$\coprod_{i\in I}\G_i$ contains the image of $\rho$. Setting
$\R_i=\rho^*(\G_i)$ then gives $\R=\coprod_{i\in I}\R_i$, with
$\rho:\R_i\ra\G_i$. Also $\R_i\cong[X_i/G]$, for $\K$-subvarieties
$X_i$ of $X$ with $X=\coprod_{i\in I}X_i$. Refining the
decomposition if necessary and using \cite[Prop.~3.5.9]{Kres} as
above, we can assume that $\G_i\cong[Z_i/K_i]$ for $Z_i$ a
quasiprojective $\K$-variety acted on by a special algebraic
$\K$-group $K_i$, and that the 1-morphisms $[X_i/G]\ra[Z_i/K_i]$
corresponding to $\rho:\R_i\ra\G_i$ are induced by $\K$-variety and
$\K$-group morphisms $\al_i:X_i\ra Z_i$ and~$\ga_i:G\ra K_i$.

Since $\G_i$ and $\phi$ are of finite type, $\phi^*(\G_i)$ is of
finite type in $\F_i$, so by \cite[Prop.~3.5.9]{Kres} again we can
write $\phi^*(\G_i)=\coprod_{j\in J_i}\F_{ij}$, for $J_i$ finite and
$\K$-substacks $\F_{ij}$ in $\F$ with $\F_{ij}\cong[Y_{ij}/H_{ij}]$,
where $Y_{ij}$ is a quasiprojective $\K$-variety acted on by a
special algebraic $\K$-group $H_{ij}$, and the 1-morphisms
$[Y_{ij}/H_{ij}]\ra[Z_i/K_i]$ corresponding to $\phi$ are induced by
$\K$-variety and $\K$-group morphisms $\be_{ij}:Y_{ij}\ra Z_i$ and
$\de_{ij}:H_{ij}\ra K_i$. Equation \eq{mi3eq7} implies an equation
similar to~\eq{mi4eq4}:
\begin{equation*}
\phi^*\bigl([(\R,\rho)]\bigr)=\ts\sum_{i\in I,\; j\in J_i}\bigl[
\bigl([X_i/G]\t_{[Z_i/K_i]}[Y_{ij}/H_{ij}],\pi_\F\bigr)\bigr].
\end{equation*}
The compatibility of $\phi^*$ with (iii) now follows using the same
argument as for `$\,\cdot\,$', changing sums and subscripts as
necessary.
\end{proof}

Many properties of the spaces $\uSF(\F)$ and their operations now
immediately follow for the $\uSF(\F,\Up,\La)$, since the operations
are defined by the same formulae on generators. In particular, we
deduce:

\begin{cor} The projections $\Pi^{\Up,\La}_*$ commute with the
operations `$\,\cdot\,$', $\phi_*$, $\phi^*$, $\ot$ on
$\uSF(*),\uSF(*,\Up,\La)$, so that\/
$\phi_*\ci\Pi^{\Up,\La}_\F=\Pi^{\Up,\La}_\G\ci\phi_*$ for
$\phi:\F\ra\G$, and so on. The analogue of Theorem \ref{mi3thm1}
holds for the spaces~$\uSF(*,\Up,\La)$.
\label{mi4cor}
\end{cor}

Next we identify~$\uSF(\Spec\K,\Up,\La)$.

\begin{prop} The map $i_\La:\La\ra\uSF(\Spec\K,\Up,\La)$ taking
$i_\La:c\mapsto c[\Spec\K]$ is an isomorphism of algebras.
\label{mi4prop2}
\end{prop}

\begin{proof} As in the proof of Theorem \ref{mi4thm1},
$\uSF(\Spec\K,\Up,\La)$ is generated over $\La$ by elements $[X/G]$
for $X$ a quasiprojective $\K$-variety acted on by a special
algebraic $\K$-group $G$. But using Definition
\ref{mi4def1}(ii),(iii) and $X\cong\Spec\K\t X$ we see that
\begin{equation*}
[X/G]=\Up([G])^{-1}[X]=\Up([G])^{-1}[\Spec\K\t X]=
\Up([G])^{-1}\Up([X])[\Spec\K],
\end{equation*}
so $\uSF(\Spec\K,\Up,\La)$ is generated over $\La$ by $[\Spec\K]$,
and $i_\La$ is {\it surjective}.

Define $\pi_\La:\uSF(\Spec\K,\Up,\La)\ra\La$ by $\pi_\La:
\sum_{i\in I}c_i[\R_i]\mapsto\sum_{i\in I}c_i\Up'([\R_i])$, for
$I$ a finite set, $c_i\in\La$ and $\Up'$ as in Theorem
\ref{mi4thm1}. Using Theorem \ref{mi4thm1} it is easy to check
$\pi_\La$ is compatible with Definition \ref{mi4def1}(i)--(iii) for
$\uSF(\Spec\K,\Up,\La)$, and so is well-defined. But
$\pi_\La([\Spec\K])=1$, so $\pi_\La\ci i_\La$ is the identity on
$\La$ by $\La$-linearity. Thus $i_\La$ is injective, and so it is an
isomorphism.
\end{proof}

Using this we show the spaces $\uSF(\F,\Up,\La)$ are at least as big
as~$\CF(\F)\ot_\Q\La$.

\begin{prop} The following map is $\La$-linear and injective:
\e
(\Pi^{\Up,\La}_\F\ci\io_\F)\ot_\Q\id_\La:
\CF(\F)\ot_\Q\La\ra\uSF(\F,\Up,\La).
\label{mi4eq6}
\e
\label{mi4prop3}
\end{prop}

\begin{proof} $\La$-linearity is obvious. Let $f\in\CF(\F)\ot_\Q\La$
and $x:\Spec\K\ra\F$ be a 1-morphism. It is easy to show from the
definitions that $i_\La^{-1}\ci x^*\bigl(((\Pi^{\Up,\La}_\F
\ci\io_\F) \ot_\Q\id_\La)(f)\bigr)=f([x])\in\La$. Thus
$((\Pi^{\Up,\La}_\F \ci\io_\F)\ot_\Q\id_\La)(f)=0$ only if
$f([x])=0$ for all $x:\Spec\K\ra\F$, that is, only if $f=0$, so
\eq{mi4eq6} is injective.
\end{proof}

This prompts the following {\it intuitive explanation} of the spaces
$\uSF(\F,\Up,\La)$, which was the author's motivation for inventing
them. In \S\ref{mi23} we considered constructible functions
$\CF(\F)$, with pushforwards $\CF^\stk(\phi)$ defined by
`integration' using the Euler characteristic $\chi$. We can think of
$\uSF(\F,\Up,\La)$ as being like constructible functions
$\CF(\F)\ot_\Q\La$ with values in $\La$, with pushforwards $\phi_*$
defined by `integration' using $\Up$ instead of~$\chi$.

In fact, pushforwards on $\CF(\F)\ot_\Q\La$ using $\Up$ do {\it
not\/} usually satisfy the analogue of \eq{mi2eq3}, because for a
non-Zariski-locally-trivial fibration $\pi:X\ra Y$ with fibre $F$ we
have $\chi(X)=\chi(F)\chi(Y)$ but $\Up([X])\ne\Up([F])\Up([Y])$ in
general, as in Example \ref{mi4ex4}. So to get a theory with the
properties we want (Theorem \ref{mi3thm1}), we must allow
$\uSF(\F,\Up,\La)$ to be larger than $\CF(\F)\ot_\Q\La$ to keep
track of $\rho:\R\ra\F$ which are non-Zariski-locally-trivial
fibrations over substacks of $\F$. All fibrations over $\Spec\K$ are
Zariski locally trivial, so $\uSF(\Spec\K,\Up,\La)$ reduces to
$\La=\CF(\Spec\K)\ot_\Q\La$, as in Proposition~\ref{mi4prop2}.

The spaces $\uSF(\F,\Up,\La)$ will be important tools in the series
\cite{Joyc2,Joyc3,Joyc4,Joyc5}. Given a $\K$-linear abelian category
$\cA$ we shall define the moduli $\K$-stack $\mathfrak{Obj}_\cA$ of
objects in $\cA$. Then $\uSF(\mathfrak{Obj}_\cA,\Up,\La)$ is
well-defined, and in \cite{Joyc3} using the Ringel--Hall algebra
idea we define an associative multiplication $*$ on it, different
from `$\,\cdot\,$', making it into a noncommutative $\La$-algebra.
Examples of this yield quantized universal enveloping algebras of
Kac--Moody algebras.

An advantage of working with spaces $\uSF(*,\Up,\La)$ rather than
$\uSF(*)$ is that because of the relations Definition
\ref{mi4def1}(i)--(iii), special properties of $\cA$ such as ${\rm
Ext}^i(X,Y)=0$ for all $X,Y\in\cA$ and $i>1$ are translated in
\cite{Joyc3} to extra identities in $\uSF(\mathfrak{Obj}_\cA,
\Up,\La)$, telling us something special about this algebra. In
\cite{Joyc5} we use Proposition \ref{mi4prop2} to project elements
of $\uSF(\mathfrak{Obj}_\cA,\Up,\La)$ to $\La$, and so define
interesting invariants in $\La$ which `count' $\tau$-(semi)stable
objects in~$\cA$.

\section{Virtual rank and projections $\Pi^\vi_n$ on $\uSF(\F)$}
\label{mi5}

Section \ref{mi4} assumed $\ell-1$ is invertible in $\La$, and we
want to relax this assumption. The basic reason for it is that
$[[\Spec\K/G]]=\Up([G])^{-1}[\Spec\K]$ in $\uSF(\Spec\K,\Up,\La)$,
and if $G$ has maximal torus $T^G\cong\bG_m^n$ then
$\Up([T^G])=(\ell-1)^n$ divides $\Up([G])$. In this section we shall
define new spaces $\uoSF,\oSF(\F,\Up,\La)$ with finer relations,
which keep track of maximal tori. These will satisfy
\begin{equation*}
\bigl[[\Spec\K/G]\bigr]=\Up\bigl([G/T^G]\bigr){}^{-1}
\bigl[[\Spec\K/T^G]\bigr]+\text{`lower order terms',}
\end{equation*}
and because $\ell-1$ does not divide $\Up\bigl([G/T^G]\bigr)$ it
will no longer be necessary for $\ell-1$ to be invertible, as we
will see in~\S\ref{mi6}.

To do this we need the difficult idea of {\it virtual rank}. The
{\it rank\/} $\rk\,G$ of an affine algebraic $\K$-group $G$ is the
dimension of any maximal torus $T^G$. We begin in \S\ref{mi51} by
defining the {\it real rank\/} projections $\Pi^\re_n:\uSF(\F)
\ra\uSF(\F)$ which project $[(\R,\rho)]$ to $[(\R_n,\rho)]$, where
$\R_n$ is the $\K$-substack of points $r\in\oR(\K)$ with stabilizer
groups $\Aut_\K(r)$ of rank $n$. This is primarily for motivation.

Section \ref{mi52} then defines analogous {\it virtual rank\/}
projections $\Pi^\vi_n:\uSF(\F)\ra\uSF(\F)$. These coincide with the
$\Pi^\re_n$ on $[(\R,\rho)]$ when $\R$ has abelian stabilizer
groups, but points $r$ with $\Aut_\K(r)$ nonabelian of rank $k$
split into components with `virtual rank' $n\le k$. Using these
ideas, \S\ref{mi53} defines spaces $\uoSF,\oSF(\F,\Up,\La)$ similar
to those of \S\ref{mi43} on which operations `$\,\cdot\,$',
$\phi_*$, $\phi^*$ and $\Pi^\vi_n$ are well-defined.

\subsection{Real rank and projections $\Pi^\re_n$}
\label{mi51}

We define a family of commuting projections
$\Pi^\re_n:\uSF(\F)\ra\uSF(\F)$ for $n=0,1,\ldots$ which project to
the part of $\uSF(\F)$ spanned by $[(\R,\rho)]$ such that the
stabilizer group $\Aut_\K(r)$ has rank $n$ for all $r\in\oR(\K)$.
The superscript `re' is short for `real', meaning that the
$\Pi^\re_n$ decompose $\uSF(\F)$ by the real (actual) rank of
stabilizer groups.

\begin{dfn} If $\R$ is an algebraic $\K$-stack and $r\in\oR(\K)$ then
$\Aut_\K(r)$ is an algebraic $\K$-group, so the rank
$\rk(\Aut_\K(r))$ is well-defined. There is a natural topology on
$\oR(\K)$, in which the open sets are $\overline{\U}(\K)$ for open
$\K$-substacks $\U\subseteq\R$. In this topology the function
$r\mapsto\rk(\Aut_\K(r))$ is {\it upper semicontinuous}. Thus, there
exist locally closed $\K$-substacks $\R_n$ in $\R$ for
$n=0,1,\ldots$, such that $\oR(\K)=\coprod_{n\ge 0}\R_n(\K)$, and
$r\in\oR(\K)$ has $\rk(\Aut_\K(r))=n$ if and only if
$r\in\oR_n(\K)$. If $\R$ is of finite type then $\R_n=\emptyset$
for~$n\gg 0$.

Now let $\F$ be an algebraic $\K$-stack with affine geometric
stabilizers, and $\uSF(\F)$ be as in \S\ref{mi3}. Define $\Q$-linear
maps $\Pi^\re_n:\uSF(\F)\ra\uSF(\F)$ for $n=0,1,\ldots$ on the
generators $[(\R,\rho)]$ of $\uSF(\F)$ by
$\Pi^\re_n:[(\R,\rho)]\mapsto [(\R_n,\rho\vert_{\R_n})]$, for $\R_n$
defined as above. If $\fS$ is a closed substack of $\R$ it is easy
to see that $\fS_n$ is a closed substack of $\R_n$ and
$(\R\sm\fS)_n=\R_n\sm\fS_n$. Thus, $\Pi^\re_n$ is compatible with
the relations \eq{mi3eq1} in $\uSF(\F)$, and is well-defined. If
$\rho:\R\ra\F$ is representable then so is $\rho\vert_{\R_n}$, so
the restriction to $\SF(\F)$ maps~$\Pi^\re_n:\SF(\F)\ra\SF(\F)$.
\label{mi5def1}
\end{dfn}

Here are some easy properties of the $\Pi^\re_n$. The proofs are
left as an exercise.

\begin{prop} In the situation above, we have:
\begin{itemize}
\setlength{\itemsep}{0pt} \setlength{\parsep}{0pt}
\item[{\rm(i)}] $(\Pi^\re_n)^2=\Pi^\re_n$, so that\/ $\Pi^\re_n$ is
a projection, and\/ $\Pi^\re_m\ci\Pi^\re_n=0$ for~$m\ne n$.
\item[{\rm(ii)}] For all\/ $f\in\uSF(\F)$ we have $f=\sum_{n\ge 0}
\Pi^\re_n(f)$, where the sum makes sense as $\Pi^\re_n(f)=0$
for~$n\gg 0$.
\item[{\rm(iii)}] If\/ $\phi:\F\ra\G$ is a $1$-morphism of algebraic
$\K$-stacks with affine geometric stabilizers
then~$\Pi^\re_n\ci\phi_*=\phi_*\ci\Pi^\re_n:\uSF(\F)\ra\uSF(\G)$.
\item[{\rm(iv)}] If\/ $f\in\uSF(\F)$, $g\in\uSF(\G)$
then~$\Pi^\re_n(f\ot
g)=\sum_{m=0}^n\Pi^\re_m(f)\ot\Pi^\re_{n-m}(g)$.
\end{itemize}
\label{mi5prop1}
\end{prop}

\subsection{Operators $\Pi^\mu$ and projections $\Pi^\vi_n$}
\label{mi52}

Next we study a family of commuting operators $\Pi^\mu$ on
$\uSF(\F)$ defined by a {\it weight function} $\mu$, which include
as special cases projections $\Pi^\vi_n$ for $n\ge 0$ similar to the
$\Pi^\re_n$ of \S\ref{mi51}. But the $\Pi^\mu,\Pi^\vi_n$ are much
more subtle and difficult than the $\Pi^\re_n$, as applied to
$[(\R,\rho)]$ they modify $\R$ in a very nontrivial way, rather than
just restricting to substacks $\R_n$. Roughly speaking, $\Pi^\mu$
replaces a point in $\R$ with stabilizer group $G$ by a linear
combination of points with stabilizer groups $C_G(T)$, for certain
subgroups $T$ of the maximal torus $T^G$ of~$G$.

From Definition \ref{mi5def2} until Lemma \ref{mi5lem3} we take $X$
to be a quasiprojective $\K$-variety acted on by an affine algebraic
$\K$-group $G$, with maximal torus~$T^G$.

\begin{dfn} If $S\subseteq T^G$ define $X^S$ to be the
$\K$-subvariety of $X$ fixed by all elements of $S$. Then $X^S$ is
closed, but not necessarily irreducible, and $\overline{X^S}(\K)=
\{x\in \overline X(\K):t\cdot x=x$ for all $t\in S\}$. For such
$X,S$ define $P$ to be the $\K$-subgroup of $T^G$ fixing the
subvariety $X^S$. Then $P$ is a closed $\K$-subgroup of $T^G$,
containing $S$, and $\overline P(\K)=\{t\in\overline{T^G}(\K):
t\cdot x=x$ for all $x\in\overline{X^S}(\K)\}$. As $S\subseteq P$ we
have $X^P\subseteq X^S$. But also $X^S\subseteq X^P$ by definition
of $P$, so $X^P=X^S$. Thus, $X^P$ and $P$ determine each other.
Define $\cP(X,T^G)$ to be the set of closed $\K$-subgroups $P$ of
$T^G$ such that $P$ is the $\K$-subgroup of $T^G$ fixing~$X^P$.
\label{mi5def2}
\end{dfn}

\begin{lem} {\rm(i)} $\cP(X,T^G)$ is finite.
\begin{itemize}
\setlength{\itemsep}{0pt}
\setlength{\parsep}{0pt}
\item[{\rm(ii)}] $\cP(X,T^G)$ is closed under intersections, with
maximal element\/ $T^G$ and minimal element\/ $P_{\rm min}$ the
subgroup of\/ $T^G$ acting trivially on~$X$.
\item[{\rm(iii)}] If\/ $S\subseteq T^G$ then $X^S=X^P$, where $P$
is the unique smallest element of\/ $\cP(X,T^G)$ containing~$S$.
\end{itemize}
\label{mi5lem1}
\end{lem}

\begin{proof} The map $x\mapsto\Stab_{T^G}(x)$ is a constructible
map from $X$ to $\K$-subgroups of $T^G$, and so realizes finitely
many values $H_1,\ldots,H_n$ say. These stratify $X$ into locally
closed subvarieties $X_1,\ldots,X_n$ with $x\in X_i$ if and only if
$\Stab_{T^G}(x)=H_i$. For any $S\subset T^G$, $X^S$ is the union of
those $X_i$ for which $S\subset H_i$, and the corresponding $P$
constructed above is the intersection of the corresponding $H_i$, or
$T^G$ if there are no $H_i$. Therefore $\cP(X,T^G)$ is exactly the
set of intersections of nonempty subsets of
$\{T^G,H_1,\ldots,H_n\}$. This proves (i) and the first two parts of
(ii). For the last part of (ii), the minimal element of $\cP(X,T^G)$
is $T^G\cap H_1\cap\cdots\cap H_n$, which is $P_{\rm min}$. Part
(iii) follows easily from the discussion in
Definition~\ref{mi5def2}.
\end{proof}

\begin{dfn} If $S\subset T^G$ then $Q=T^G\cap C(C_G(S))$ is a closed
$\K$-subgroup of $T^G$ containing $S$. As $S\subseteq Q$ we have
$C_G(Q)\subseteq C_G(S)$. But $Q$ commutes with $C_G(S)$, so
$C_G(S)\subseteq C_G(Q)$. Thus $C_G(S)=C_G(Q)$. So $Q=T^G\cap
C(C_G(Q))$, and $Q$ and $C_G(Q)$ determine each other, given
$G,T^G$. Define $\cQ(G,T^G)$ to be the set of closed $\K$-subgroups
$Q$ of $T^G$ such that~$Q=T^G\cap C(C_G(Q))$.
\label{mi5def3}
\end{dfn}

\begin{lem} {\rm(i)} $\cQ(G,T^G)$ is finite.
\begin{itemize}
\setlength{\itemsep}{0pt}
\setlength{\parsep}{0pt}
\item[{\rm(ii)}] $\cQ(G,T^G)$ is closed under intersections, with
maximal element\/ $T^G$ and minimal element\/~$Q_{\rm min}=T^G\cap
C(G)$.
\item[{\rm(iii)}] If\/ $S\subseteq T^G$ then $C_G(S)=C_G(Q)$,
where $Q$ is the unique smallest element of\/ $\cQ(G,T^G)$
containing~$S$.
\end{itemize}
\label{mi5lem2}
\end{lem}

\begin{proof} The proof is similar to Lemma \ref{mi5lem1}. The map
$t\mapsto C_G(\{t\})$ is constructible from $T^G$ to closed
$\K$-subgroups of $G$, and realizes finitely many values
$H_1,\ldots,H_n$. Set $Q_i=T^G\cap C(H_i)$. Then $\cQ(G,T^G)$ is the
set of intersections of nonempty subsets of
$\{T^G,Q_1,\ldots,Q_n\}$. We leave the details to the reader.
\end{proof}

We calculate $\cQ(G,T^G)$ for the case~$G=\GL(m,\K)$.

\begin{ex} Set $G=\GL(m,\K)$ with maximal torus $T^G=\bG_m^m$, the
subgroup of diagonal matrices. Fix $t\in T^G$, which may be written
$\diag(t_1,\ldots,t_m)$ for $t_i\in\bG_m$. Let $t_1,\ldots,t_m$
realize $n$ distinct values $u_1,\ldots,u_n$. Then there is a unique
surjective map $\phi:\{1,\ldots,m\}\ra\{1,\ldots,n\}$ with
$t_i=u_{\phi(i)}$. It is easy to show that $C_G(\{t\})$ is the
subgroup of matrices $(A_{ij})_{i,j=1}^m$ in $\GL(m,\K)$ with
$A_{ij}=0$ if $\phi(i)\ne\phi(j)$. Hence
\e
C_G(\{t\})\cong\ts\prod_{k=1}^n
\GL\bigl(\md{\phi^{-1}(\{k\})},\K\bigr).
\label{mi5eq1}
\e

The centre $C(C_G(\{t\}))$ of $C_G(\{t\})$, which agrees with
$T^G\!\cap\!C(C_G(\{t\}))$, is
\e
\bigl\{\diag(q_1,\ldots,q_m):\text{$q_i\in\bG_m$, $q_i=q_j$ if
$\phi(i)=\phi(j)$, all $i,j$}\bigr\}\cong\bG_m^n.
\label{mi5eq2}
\e
Since $\cQ(G,T^G)$ is the set of $T^G\cap C(C_G(\{t\}))$ for all
$t\in T^G$, we see $\cQ(G,T^G)$ is the set of tori \eq{mi5eq2} for
all $1\le n\le m$ and surjective~$\phi:\{1,\ldots,m\}\!\ra\!
\{1,\ldots,n\}$.
\label{mi5ex}
\end{ex}

\begin{dfn} For $P'\subseteq P$ in $\cP(X,T^G)$ and $Q'\subseteq Q$
in $\cQ(G,T^G)$, set
\ea
m^X_{T^G}(P',P)=\sum_{\substack{\text{$A\subseteq\{\hat
P\in\cP(X,T^G): \hat P\subseteq P\}:P\in A$, $\bigcap_{\hat P\in
A}\hat P=P'$}}}(-1)^{\md{A}-1},
\label{mi5eq3}
\\
n^G_{T^G}(Q',Q)=\sum_{\substack{\text{$B\subseteq\{\hat
Q\in\cQ(G,T^G): \hat Q\subseteq Q\}:Q\in B$, $\bigcap_{\hat Q\in
B}\hat Q=Q'$}}}(-1)^{\md{B}-1}.
\label{mi5eq4}
\ea

Define $\cR(X,G,T^G)=\{P\cap Q:P\in\cP(X,T^G)$, $Q\in\cQ(G,T^G)$.
For $P\in\cP(X,T^G)$, $Q\in\cQ(G,T^G)$ and $R\in\cR(X,G,T^G)$ with
$R\subseteq P\cap Q$, define
\e
M^X_G(P,Q,R)\!=\!\left\vert\frac{N_G(T^G)}{C_G(Q)\!\cap\!N_G(T^G)}
\right\vert^{-1}\!\cdot\!\!\!\!\!
\sum_{\begin{subarray}{l} \text{$P'\in\cP(X,T^G)$, $Q'\in\cQ(G,T^G)$:}\\
\text{$P'\subseteq P$, $Q'\subseteq Q$, $R=P'\cap Q'$}
\end{subarray}\!\!\!\!\!\!\!\!\!\!\!\!\!\!\!\!\!\!\!\!\!\!\!\! }
\!\!\!\!\!\!\!\!\!\!\!\!\!\!\!\! m^X_{T^G}(P',P)n^G_{T^G}(Q',Q).
\label{mi5eq5}
\e
Here $C_G(T^G)\subseteq C_G(Q)$ as $Q\subseteq T^G$, so $C_G(Q)\cap
N_G(T^G)$ is a subgroup of $N_G(T^G)$ containing $C_G(T^G)$. But
$C_G(T^G)$ is of finite index in $N_G(T^G)$ as
$W(G,T^G)=N_G(T^G)/C_G(T^G)$ is finite. Hence $N_G(T^G)/(C_G(Q)\cap
N_G(T^G))$ is finite, and $M^X_G(P,Q,R)$ is well-defined.
\label{mi5def4}
\end{dfn}

\begin{lem} If\/ $M^X_G(P,Q,R)\ne 0$ then $P$ is the smallest
element of\/ $\cP(X,T^G)$ containing $P\cap Q$, and\/ $Q$ the
smallest element of\/ $\cQ(G,T^G)$ containing $P\cap Q$. Therefore
$X^P=X^{P\cap Q}$ and\/~$C_G(Q)=C_G(P\cap Q)$.
\label{mi5lem3}
\end{lem}

\begin{proof} Using \eq{mi5eq3}--\eq{mi5eq4} we may rewrite
the sum in \eq{mi5eq5} as
\e
\sum_{\substack{\text{$A\subseteq\{\hat P\in\cP(X,T^G):
\hat P\subseteq P\}$, $B\subseteq\{\hat Q\in\cQ(G,T^G):
\hat Q\subseteq Q\}:$}\\
\text{$P\in A$, $Q\in B$, $\bigcap_{\hat P\in A}\cap\bigcap_{\hat
Q\in B}\hat Q=R$}}}(-1)^{\md{A}+\md{B}}.
\label{mi5eq6}
\e
Suppose there exists $P'\in\cP(X,T^G)$ with $P\cap Q\subseteq
P'\subset P$ and $P'\ne P$. Then in \eq{mi5eq6} the intersection
$\bigcap_{\hat P\in A}\cap\bigcap_{\hat Q\in B}\hat Q$ is unchanged
by whether $P'\in A$, as it lies in $P\cap Q$. Thus for each pair
$A,B$ in \eq{mi5eq6} with $P'\notin A$ there corresponds another
pair $A\cup\{P'\},B$, and the total contribution of both is
$(-1)^{\md{A}+\md{B}}+(-1)^{\md{A}+1+\md{B}}=0$. So \eq{mi5eq6} and
$M^X_G(P,Q,R)$ are zero. Conversely, if $M^X_G(P,Q,R)\ne 0$ there
exists no such $P'$, so $P$ is the smallest element of $\cP(X,T^G)$
containing $P\cap Q$. The argument for $Q$ is the same. The final
part follows from Lemmas \ref{mi5lem1}(iii) and~\ref{mi5lem2}(iii).
\end{proof}

Now we define some linear maps~$\Pi^\mu:\uSF(\F)\ra\uSF(\F)$.

\begin{dfn} A {\it weight function} $\mu$ is a map
\begin{equation*}
\mu:\bigl\{\text{$\K$-groups $\bG_m^k\!\t\!K$, $k\!\ge\!0$, $K$
finite abelian, up to isomorphism}\bigr\}\!\ra\!\Q.
\end{equation*}
For any algebraic $\K$-stack $\F$ with affine geometric stabilizers,
we will define a linear map $\Pi^\mu:\uSF(\F)\ra\uSF(\F)$. Now
$\uSF(\F)$ is generated by elements $[(\R,\rho)]$ with $\R$
1-isomorphic to a global quotient $[X/G]$, for $X$ a quasiprojective
$\K$-variety and $G$ a special algebraic $\K$-group, with maximal
torus $T^G$. For such $\R,\rho,X,G,T^G$ define
\e
\begin{gathered}
\Pi^\mu\bigl([(\R,\rho)]\bigr)\!=
\!\!\!\!\!\!\!\!\!\!\!\!\!\!\!\!\!\!\!\!\!\!\!
\sum_{\substack{\text{$P\in\cP(X,T^G)$, $Q\in\cQ(G,T^G)$ and}\\
\text{$R\!\in\!\cR(X,G,T^G):R\!\subseteq\!P\cap Q$,
$M^X_G(P,Q,R)\!\ne\!0$}}}
\begin{aligned}[t]
&M^X_G(P,Q,R)\mu(R)\,\cdot\\
&\bigl[\bigl([X^P/C_G(Q)],\rho\!\ci\!\io^{P\cap Q}\bigr)\bigr].
\end{aligned}
\end{gathered}
\label{mi5eq7}
\e
Here $X^P=X^{P\cap Q}$ and $C_G(Q)=C_G(P\cap Q)$ by Lemma
\ref{mi5lem3}, so $X^P$ is $C_G(Q)$-invariant, and the stack
$[X^P/C_G(Q)]$ is well-defined. The inclusions $X^P\subseteq X$,
$C_G(Q)\subseteq G$ induce a 1-morphism $\io^{P\cap
Q}:[X^P/C_G(Q)]\ra[X/G]$. As the $\io^{P\cap Q}$ are representable,
if $[(\R,\rho)]\in\SF(\F)$ then~$\Pi^\mu([(\R,\rho)])\in\SF(\F)$.
\label{mi5def5}
\end{dfn}

An informal but helpful way to rewrite \eq{mi5eq7} is:
\e
\begin{split}
&\Pi^\mu\bigl([(\R,\rho)]\bigr)=\\
&\int_{t\in T^G}\frac{\md{\{w\in W(G,T^G):w\cdot
t=t\}}}{\md{W(G,T^G)}}\,\bigl[\bigl([X^{\{t\}}/
C_G(\{t\})],\rho\ci\io^{\{t\}}\bigr)\bigr]\d\mu.
\end{split}
\label{mi5eq8}
\e
Here $\d\mu$ is a {\it measure} on a class of subsets of $T^G$
described below. Lemmas \ref{mi5lem1}(iii), \ref{mi5lem2}(iii) give
$[([X^{\{t\}}/C_G(\{t\})],\rho\ci\io^{\{t\}})]
\!=\![([X^P/C_G(Q)],\rho\ci\io^{P\cap Q})]$, with $P,Q$ the unique
smallest elements of $\cP(X,T^G),\cQ(G,T^G)$ containing $t$. Also
\e
\begin{aligned}
&\frac{\md{\{w\in W(G,T^G):w\cdot t=t\}}}{\md{W(G,T^G)}}
=\frac{\md{(C_G(\{t\})\cap N_G(T^G))/C_G(T^G)}}{\md{W(G,T^G)}}=\\
&\frac{\md{(C_G(Q)\cap N_G(T^G))/C_G(T^G)}}{\md{N_G(T^G)/
C_G(T^G)}}=\left\vert\frac{N_G(T^G)}{C_G(Q)\cap N_G(T^G)}
\right\vert^{-1}.
\end{aligned}
\label{mi5eq9}
\e
Thus the integrand in \eq{mi5eq8} at $t$ depends only on~$P,Q$.

Therefore the subsets of $T^G$ the measure $\d\mu$ must be defined
upon for \eq{mi5eq8} to make sense, are those generated from
$\cP(X,T^G),\cQ(G,T^G)$ by Boolean operations. This is determined
uniquely by setting $\d\mu(R)=\mu(R)$ for $R\in\cR(X,G,T^G)$. We
find that for $P\in\cP(X,T^G)$ and $Q\in\cQ(G,T^G)$ we have
\e
\begin{gathered}
\d\mu
\raisebox{-5pt}{\begin{Large}$\displaystyle\Bigl[$\end{Large}}
P\cap Q\sm\!\!\bigcup_{\begin{subarray}{l}
\text{$P'\in\cP(X,T^G)$, $Q'\in\cQ(G,T^G)$:}\\
\text{$P'\subseteq P$, $Q'\subseteq Q$,
$(P,Q)\ne(P',Q')$}\end{subarray}}P'\cap Q'
\raisebox{-5pt}{\begin{Large}$\displaystyle\Bigr]$\end{Large}}=
\\
\sum_{\text{$P'\in\cP(X,T^G)$, $Q'\in\cQ(G,T^G):P'\subseteq P$,
$Q'\subseteq Q$}\!\!\!\!\!\!\!\!\!\!\!\!\!\!\!\!\!\!\!\!\! } \!\!\!
m^X_{T^G}(P',P)n^G_{T^G}(Q',Q)\mu(P'\cap Q').
\end{gathered}
\label{mi5eq10}
\e
As $T^G$ is the disjoint union over $P,Q$ of the sets $[\cdots]$ on
the top line of \eq{mi5eq10}, comparing \eq{mi5eq5} and
\eq{mi5eq7}--\eq{mi5eq10} we see \eq{mi5eq7} and \eq{mi5eq8} are
equivalent.

As the integrand in \eq{mi5eq8} is invariant under the action of
$W(G,T^G)$, we can simplify \eq{mi5eq8} further by pushing the
integration down to~$T^G/W(G,T^G)$:
\e
\Pi^\mu\bigl([(\R,\rho)]\bigr)=\int_{tW(G,T^G)\in T^G/W(G,T^G)}
\bigl[\bigl([X^{\{t\}}/C_G(\{t\})],\rho\ci\io^{\{t\}}\bigr)\bigr]
\d\mu.
\label{mi5eq11}
\e
Now $T^G/W(G,T^G)$ is a natural object in algebraic group theory, as
it is isomorphic to $G^{\rm ss}/\Ad(G)$, where $G^{\rm ss}$ is the
open set of {\it semisimple elements} of $G$. In the quotient stack
$[G^{\rm ss}/\Ad(G)]$ the stabilizer group $\Aut_\K(tW(G,T^G))$ is
$C_G(\{t\})$. So \eq{mi5eq11} is an integral over $[G^{\rm
ss}/\Ad(G)]$ of a function of the stabilizer group. Probably there
is some extension of this construction to integrate over all of
$[G/\Ad(G)]$, replacing $T^G$ by a {\it Borel subgroup} perhaps, but
we do not consider it. We show $\Pi^\mu$ is independent of choices
in its definition.

\begin{thm} In the situation above, $\Pi^\mu\bigl([(\R,\rho)]
\bigr)$ is independent of the choices of\/ $X,G,T^G$ and\/
$1$-isomorphism $\R\cong[X/G]$, and\/ $\Pi^\mu$ extends to unique
linear maps $\Pi^\mu:\uSF(\F)\ra\uSF(\F)$
and\/~$\Pi^\mu:\SF(\F)\ra\SF(\F)$.
\label{mi5thm1}
\end{thm}

\begin{proof} Suppose $\R$ is 1-isomorphic to $[X/G]$ and
$[Y/H]$, for $X,Y$ quasiprojective $\K$-varieties and $G,H$ special
algebraic $\K$-groups acting on $X,Y$. Define $Z,I,Z_i,X_i,Y_i$ as
in the proof of Proposition \ref{mi4prop1}. Since \eq{mi5eq7} is
additive over $X=\coprod_{i\in I}X_i$, it is enough to prove it
gives the same answer for $[X_i/G]$ and $[Y_i/H]$, for all $i\in I$.
Thus for simplicity we replace $X_i,Y_i,Z_i$ by $X,Y,Z$. That is, we
have a quasiprojective $\K$-variety $Z$ with a $G\t H$-action, such
that the $G$- and $H$-actions are free and induce projections
$\pi_Y:Z\ra Y$ and $\pi_X:Z\ra X$ which are $G$- and $H$-principal
bundles respectively. Also $[Z/(G\t H)]\cong[X/G]\cong[Y/H]$. Fix
maximal tori $T^G,T^H$ in~$G,H$.

Let $t\in T^G$. We will relate $X^{\{t\}}$ to $Z^{\{(t,\ti t)\}}$
for $\ti t\in T^H$. Suppose $x\in X^{\{t\}}$, and $z\in\pi_X^{-1}
(\{x\})\subset Z$. Then the projection $\pi_G:G\t H\ra G$ induces an
isomorphism $\Stab_{G\t H}(z)\ra\Stab_G(x)$. So as $t\in\Stab_G(x)$
there is a unique $h\in H$ with $(t,h)\in\Stab_{G\t H}(z)$. Now
$G,H$ and $G\t H$ are connected, as $G,H$ are special. Thus elements
of $G,H,G\t H$ are semisimple if and only if they lie in a maximal
torus.

As $t\in T^G$ it is semisimple in $G$. So $t$ is semisimple in
$\Stab_G(x)\subseteq G$. Thus $(t,h)$ is semisimple in $\Stab_{G\t
H}(z)$ as $\Stab_G(x)\cong\Stab_{G\t H}(z)$, and so $(t,h)$ is
semisimple in $G\t H$. Therefore $(t,h)$ lies in a maximal torus of
$G\t H$, which we may take to be $T^G\t\ti T^H$. As all maximal tori
in $H$ are conjugate, $\ti h\ti T^H\ti h^{-1}=T^H$ for some $\ti
h\in H$, and so $\ti t=\ti hh\ti h^{-1}$ lies in $T^H$. Hence $\ti
z=\ti h\cdot z$ also lies in $\pi_X^{-1}(\{x\})\subset Z$, and is
fixed by~$(t,\ti t)\in T^G\t T^H$.

Since $\pi_X:Z\ra X$ is a principal $H$-bundle, $\pi_X^{-1}(\{x\})$
is a copy of $H$. It is now easy to see that
$z'\in\pi_X^{-1}(\{x\})$ is fixed by $(t,t')$ for some $t'\in T^H$
only if $t'=w\cdot\ti t$ for some $w\in W(H,T^H)$, and the set of
such $z'$ is a copy of $C_H(\{t'\})$. Thus $\pi_H$ induces a
morphism of $\K$-varieties
\e
\ts\coprod_{t'\in T^H:Z^{\{(t,t')\}}\ne\emptyset}
Z^{\{(t,t')\}}/C_H(\{t'\})\longra X^{\{t\}},
\label{mi5eq12}
\e
whose fibre over each $x\in X^{\{t\}}$ consists of one point in
$Z^{\{(t,t')\}}/C_H(\{t'\})$ for each $t'$ in exactly one orbit of
$W(H,T^H)$ in $T^H$. Dividing by $C_G(\{t\})$ and writing as
elements of $\uSF(\F)$ shows that
\e
\begin{split}
\bigl[\bigl([X^{\{t\}}/C_G(\{t\})],\rho\ci\io^{\{t\}}\bigr)\bigr]=
\!\!\!\sum_{t'\in T^H:Z^{\{(t,t')\}}\ne\emptyset}
\!\!\!\!\!\!\!\!\!\!
\frac{\md{\{w\!\in\!W(H,T^H):w\cdot t'\!=\!t'\}}}{\md{W(H,T^H)}}&\\
\bigl[\bigl([Z^{\{(t,t')\}}/C_G(\{t\})\t C_H(\{t'\}),
\rho\ci\io^{\{(t,t')\}}\bigr)\bigr]&,
\end{split}
\label{mi5eq13}
\e
where $\md{\{w\in W(H,T^H):w\cdot t'=t'\}}/\md{W(H,T^H)}$ is
$1/\md{W(H,T^H)\cdot t'}$, and compensates for the multiplicity
of~\eq{mi5eq12}.

Multiplying \eq{mi5eq13} by $\md{\{w\in W(G,T^G):w\cdot
t=t\}}/\md{W(G,T^G)}$, integrating over $t\in T^G$ as in \eq{mi5eq8}
and using $C_{G\t H}(\{(t,t')\})=C_G(\{t\})\t C_H(\{t'\})$ and
$W(G\t H,T^G\t T^H)=W(G,T^G)\t W(G,T^H)$ shows that \eq{mi5eq8}
gives the same answer using $X,G,T^G$ and $Z,G\t H,T^G\t T^H$. Here
in relating `integrals' over $T^G$ and $T^G\t T^H$ we use the fact
that the measure $\d\mu$ is defined using $\mu$, which depends only
on isomorphism classes of $\K$-groups, and not on anything special
to $T^G$ or~$T^G\t T^H$.

Exchanging $X,G,T^G$ and $Y,H,T^H$, we see \eq{mi5eq8} gives the
same answer using $Y,H,T^H$ and $Z,G\t H,T^G\t T^H$, so it gives the
same answer using $X,G,T^G$ and $Y,H,T^H$. As \eq{mi5eq8} is
equivalent to \eq{mi5eq7}, this proves the first part of the
theorem. The second part is then obvious, since such $[(\R,\rho)]$
span $\uSF(\F)$, and \eq{mi5eq7} is compatible with the relations
\eq{mi3eq1} defining~$\uSF(\F)$.
\end{proof}

Here are some properties of the operators~$\Pi^\mu$.

\begin{thm} {\rm(a)} $\Pi^1$ defined using $\mu\equiv 1$ is the
identity on~$\uSF(\F)$.
\begin{itemize}
\setlength{\itemsep}{0pt}
\setlength{\parsep}{0pt}
\item[{\rm(b)}] If\/ $\phi:\F\ra\G$ is a $1$-morphism of algebraic
$\K$-stacks with affine geometric stabilizers
then~$\Pi^\mu\ci\phi_*=\phi_*\ci\Pi^\mu:\uSF(\F)\ra\uSF(\G)$.
\item[{\rm(c)}] If\/ $\mu_1,\mu_2$ are weight functions as in
Definition \ref{mi5def5} then $\mu_1\mu_2$ is also a weight function
and\/~$\Pi^{\mu_2}\ci\Pi^{\mu_1}=\Pi^{\mu_1}\ci\Pi^{\mu_2}=
\Pi^{\mu_1\mu_2}$.
\end{itemize}
\label{mi5thm2}
\end{thm}

\begin{proof} Arguing as in Lemma \ref{mi5lem3} using
$P'=P_{\rm min}$ and $Q'=Q_{\rm min}$ we find
\begin{equation*}
\sum_{R\in\cR(X,G,T^G):R\subseteq P\cap Q}M^X_G(P,Q,R)=
\begin{cases} 1, & \text{$P=P_{\rm min}$, $Q=Q_{\rm min}$,}\\
0, & \text{otherwise.}\end{cases}
\end{equation*}
Substituting this into \eq{mi5eq7} with $\mu\equiv 1$ gives
\begin{equation*}
\Pi^1\bigl([(\R,\rho)]\bigr)=\bigl[\bigl([X^{P_{\rm min}}/
C_G(Q_{\rm min})],\rho\!\ci\!\io^{P_{\rm min}\cap Q_{\rm
min}}\bigr)\bigr]=[(\R,\rho)],
\end{equation*}
since $X^{P_{\min}}=X$ and $C_G(Q_{\rm min})=G$. This proves (a),
and (b) is immediate.

For (c), note that if $P'\in\cP(X,T^G)$ then $\cP(X^{P'},T^G)=
\{P\in\cP(X,T^G):P'\subseteq P\}$, and for such $P$ we have
$(X^{P'})^P=X^P$. Similarly, if $Q'\in\cQ(G,T^G)$ then $T^G$ is a
maximal torus in $C_G(Q')$, and $\cQ(C_G(Q'),T^G)=\{Q\in\cQ(G,T^G):
Q'\subseteq Q\}$. Therefore $\cR(X^{P'},C_G(Q'),T^G)=\{R\in
\cR(X,G,T^G):P'\cap Q'\subseteq R\}$. Using these and \eq{mi5eq7} in
the situation of Definition \ref{mi5def5} gives
\e
\begin{gathered}
\Pi^{\mu_2}\ci\Pi^{\mu_1}\bigl([(\R,\rho)]\bigr)=
\!\!\!\!\!\!\!\!\!\!\!\!\!\!\!\!\!\!\!\!\!\!\!\!\!\!\!\!\!\!\!
\sum_{\substack{\text{$P,P'\in\cP(X,T^G)$,
$Q,Q'\in\cQ(G,T^G)$}\\
\text{and $R,R'\in\cR(X,G,T^G):$}\\
\text{$P'\!\subseteq\!P$, $Q'\!\subseteq\!Q$,
$R'\!\subseteq\!P'\cap Q'\!\subseteq\!R\!\subseteq\!P\cap Q$,}\\
\text{$M^X_G(P',Q',R')\!\ne\!0$,
$M^{X^{P'}}_{C_G(Q')}(P,Q,R)\!\ne\!0$}}}
\begin{aligned}[t]
&M^X_G(P',Q',R')\mu_1(R')\cdot\\
&M^{X^{P'}}_{C_G(Q')}(P,Q,R)\mu_2(R)\cdot\\
&\bigl[\bigl([X^P/C_G(Q)],\rho\!\ci\!\io^{P\cap Q}\bigr)\bigr].
\end{aligned}
\end{gathered}
\label{mi5eq14}
\e

Now a combinatorial calculation with \eq{mi5eq3}--\eq{mi5eq5} shows
for fixed $P,Q,R,R'$ in \eq{mi5eq14} with $R'\subseteq R\subseteq
P\cap Q$ we have
\e
\begin{gathered}
\sum_{\substack{\text{$P'\in\cP(X,T^G)$, $Q'\in\cQ(G,T^G)$:}\\
\text{$P'\!\subseteq\!P$, $Q'\!\subseteq\!Q$, $R'\!\subseteq\!P'\cap
Q'\!\subseteq\!R$}}} M^X_G(P',Q',R')M^{X^{P'}}_{C_G(Q')}(P,Q,R)\\
=\begin{cases}
M^X_G(P,Q,R), & \text{$R'=R$,}\\
0, & \text{otherwise.}
\end{cases}
\end{gathered}
\label{mi5eq15}
\e
Combining \eq{mi5eq7}, \eq{mi5eq14} and \eq{mi5eq15} shows
$\Pi^{\mu_2}\ci\Pi^{\mu_1}([(\R,\rho)])=\Pi^{\mu_1\mu_2}
([(\R,\rho)])$, so that $\Pi^{\mu_2}\ci\Pi^{\mu_1}=
\Pi^{\mu_1\mu_2}$. Exchanging $\mu_1,\mu_2$ then
gives~$\Pi^{\mu_1}\ci\Pi^{\mu_2}=\Pi^{\mu_1\mu_2}$.
\end{proof}

In contrast to (b), the $\Pi^\mu$ do {\it not\/} in general commute
with pullbacks $\phi^*:\uSF(\G)\ra\uSF(\F)$ for finite type
1-morphisms $\phi:\F\ra\G$. We can now define operators $\Pi^\vi_n$,
similar to the operators $\Pi^\re_n$ of~\S\ref{mi51}.

\begin{dfn} For $n\ge 0$, define $\Pi^\vi_n$ to be the operator
$\Pi^{\mu_n}$ defined with weight $\mu_n$ given by $\mu_n([H])=1$ if
$\dim H=n$ and $\mu_n([H])=0$ otherwise, for all $\K$-groups
$H\cong\bG_m^k\t K$ with $K$ a finite abelian group.
\label{mi5def6}
\end{dfn}

The analogue of Proposition \ref{mi5prop1} holds for
the~$\Pi^\vi_n$.

\begin{prop} In the situation above, we have:
\begin{itemize}
\setlength{\itemsep}{0pt} \setlength{\parsep}{0pt}
\item[{\rm(i)}] $(\Pi^\vi_n)^2=\Pi^\vi_n$, so that\/ $\Pi^\vi_n$ is
a projection, and\/ $\Pi^\vi_m\ci\Pi^\vi_n=0$ for~$m\ne n$.
\item[{\rm(ii)}] For all\/ $f\in\uSF(\F)$ we have $f=\sum_{n\ge 0}
\Pi^\vi_n(f)$, where the sum makes sense as $\Pi^\vi_n(f)=0$
for~$n\gg 0$.
\item[{\rm(iii)}] If\/ $\phi:\F\ra\G$ is a $1$-morphism of algebraic
$\K$-stacks with affine geometric stabilizers
then~$\Pi^\vi_n\ci\phi_*=\phi_*\ci\Pi^\vi_n:\uSF(\F)\ra\uSF(\G)$.
\item[{\rm(iv)}] If\/ $f\in\uSF(\F)$, $g\in\uSF(\G)$
then~$\Pi^\vi_n(f\ot
g)=\sum_{m=0}^n\Pi^\vi_m(f)\ot\Pi^\vi_{n-m}(g)$.
\end{itemize}
\label{mi5prop2}
\end{prop}

\begin{proof} Part (i) is immediate from Theorem \ref{mi5thm2}(c) and
Definition \ref{mi5def6}. For (ii), in Definition \ref{mi5def5} we
have $\Pi^\vi_n([(\R,\rho)])=0$ for $n>\rk\,G$ as $\mu(R)=0$ for all
$R$ in \eq{mi5eq7}, so $\Pi^\vi_n(f)=0$ for $n\gg 0$. The first part
of (ii) then follows from Theorem \ref{mi5thm2}(a), as $\Pi^\mu$ is
additive in $\mu$ and $1=\sum_{n\ge 0}\mu_n$. Theorem
\ref{mi5thm2}(b) gives (iii), and (iv) is not difficult to prove
directly from Definition \ref{mi5def5} and the fact that~$\mu_n(R\t
R')=\sum_{m=0}^n\mu_m(R)\mu_{n-m}(R')$.
\end{proof}

To get a feel for what the operators $\Pi^\mu$ and $\Pi^\vi_n$ do,
consider the case $X=\Spec\K$ and $G=\bG_m^k$, so that
$\R=[\Spec\K/\bG_m^k]$ is a point with torus stabilizer $\bG_m^k$.
Then $\Pi^\mu\bigl([(\R,\rho)]\bigr)=
\mu\bigl([\bG_m^k]\bigr)[(\R,\rho)]$, so that $\Pi^\vi_n
\bigl([(\R,\rho)]\bigr)=[(\R,\rho)]$ if $k=n$ and 0 otherwise. More
generally, if $\R=[\Spec\K/G]$ for $G$ {\it abelian} then $\Pi^\vi_n
\bigl([(\R,\rho)]\bigr)=[(\R,\rho)]$ if $\rk\,G=n$ and 0 otherwise.
Thus $\Pi^\vi_n$ and $\Pi^\re_n$ coincide on points with abelian
stabilizers.

However, if $G$ is nonabelian then $\Pi^\vi_n[([\Spec\K/G],\rho)]$
may be nonzero when $\rk\,C(G)\le n\le\rk\,G$, and is zero outside
this range. We think of $[\Spec\K/G]$ as being like a linear
combination of points with {\it virtual rank\/} in the range $\rk\,
C(G)\le n\le\rk\,G$, and $\Pi^\vi_n$ as projecting to the part of
$[\Spec\K/G]$ with virtual rank~$n$.

We briefly sketch a conjectural alternative approach to the
operators $\Pi^\mu$, which may make them seem more natural. Let $G$
be an affine algebraic $\K$-group, and $\R$ a finite type algebraic
$\K$-stack. Then we can form a $\K$-stack
$\mathfrak{Hom}([\Spec\K/G],\R)$ by defining for each $\K$-scheme
$U$ the groupoid
\begin{equation*}
\mathfrak{Hom}([\Spec\K/G],\R)(U)=
\Hom\bigl(U\t[\Spec\K/G],\R\bigr),
\end{equation*}
and for each morphism of $\K$-schemes $\phi:U\ra V$ the functor
\begin{equation*}
\mathfrak{Hom}([\Spec\K/G],\R)(\phi):
\mathfrak{Hom}([\Spec\K/G],\R)(V)\longra
\mathfrak{Hom}([\Spec\K/G],\R)(U)
\end{equation*}
induced by composition of 1-morphisms and 2-morphisms with the
1-morphism $\phi\t\id_{[\Spec\K/G]}:U\t[\Spec\K/G]\ra
V\t[\Spec\K/G]$ and its identity 2-morphism. Taking $U=\Spec\K$, we
see that the $\K$-points of $\mathfrak{Hom}([\Spec\K/G],\R)$ are
1-morphisms $[\Spec\K/G]\ra\R$. There is a projection
$\Pi:\mathfrak{Hom}([\Spec\K/G], \R)\ra\R$ corresponding to
composition of 1-morphisms~$\Spec\K\ra[\Spec\K/G]\ra\R$.

If we used a general $\K$-stack $\F$ in place of $[\Spec\K/G]$ here
then $\mathfrak{Hom}(\F,\R)$ would be not even locally of finite
type -- essentially, infinite-dimensional. But
$\mathfrak{Hom}([\Spec\K/G],\R)$ is locally of finite type. It may
not be of finite type because the fibre of $\Pi$ over $r\in\oR(\K)$
is $[\Hom(G,\Aut_\K(r))/\Ad(\Aut_\K(r))]$, where the $\K$-group
morphisms $\Hom(G,\Aut_\K(r))$ may have infinitely many components.

Roughly speaking, one might construct the $\Pi^\mu$ as follows. For
an algebraic $\K$-group $T$ of the form $\bG_m^k\t K$ for $K$ finite
abelian, we restrict to a $\K$-substack $\mathfrak{Hom}([\Spec\K/T],
\R)^{\rm reg}$ of points in $\mathfrak{Hom}([\Spec\K/T],\R)$ with
some extra properties, and $\mathfrak{Hom}([\Spec\K/T],\R)^{\rm
reg}/\Aut(T)$ is of finite type. Then $\Pi^\mu([(\R,\rho)])$ is a
linear combination over $T$ of terms like $[(\mathfrak{Hom}(
[\Spec\K/T],\R)^{\rm reg}/\Aut(T),\rho\ci\Pi)]$. These ideas might
be worth further investigation, if anyone is interested.

Finally we discuss a generalization $\hat\Pi^\nu_{\smash{\F}}$ of
the operators $\Pi^\mu$ which will be useful in \cite[\S 5]{Joyc3}.
The idea is that $\hat\Pi^\nu_\F\bigl([(\R,\rho)]\bigr)$ depends not
just on subtori $T$ of stabilizer groups $\Aut_\K(r)$ in $\R$, but
also on the morphism $\rho_*:T\ra\Aut_\K(f)$ to the stabilizer group
of $f=\rho_*(r)$ in $\F$. Thus the weight function $\nu$ is a
function of all morphisms $\rho_*:T\ra\Aut_\K(f)$, which makes it
unwieldy to define.

\begin{dfn} Let $\F$ be an algebraic $\K$-stack with affine
geometric stabilizers. An $\F$-{\it weight function} is a map
\begin{align*}
\nu:\bigl\{(T,f,\phi):\,&\text{$T$ a $\K$-group isomorphic to
$\bG_m^k\!\t\!K$, $K$ finite abelian,}\\
&\text{$f\in\oF(\K)$, $\phi:T\ra\Aut_\K(f)$ a $\K$-group
morphism}\bigr\}\longra\Q,
\end{align*}
which satisfies $\nu(T,f,\phi)=\nu(T',f,\phi\ci\io)$ if $\io:T'\ra
T$ is a $\K$-group isomorphism, and is {\it locally constructible}
in $f,\phi$. That is, $\nu$ induces a locally constructible function
$\mathfrak{Hom}([\Spec\K/T],\F)\ra\Q$ for each fixed $T$, in the
notation above.

Let $[(\R,\rho)]\in\uSF(\F)$ with $\R\cong[X/G]$ for $X$ a
quasiprojective $\K$-variety and $G$ a special algebraic $\K$-group
with maximal torus $T^G$. For $P\in\cP(X,T^G)$ and
$R\in\cR(X,G,T^G)$ with $R\subseteq P$ and $c\in\Q$, define
\begin{equation*}
X^{P,R}_{\nu,c}=\bigl\{x\in \overline{X^P}(\K):
\nu\bigl(R,(\rho\ci\pi)_*(x),\rho_*\vert_R\bigr)=c\bigr\},
\end{equation*}
writing $\pi:X\ra[X/G]\cong\R$ for the projection, so that if $x\in
\overline{X^P}(\K)$ then $r=\pi_*(x)\in\oR(\K)$, and
$f=(\rho\ci\pi)_*(x)\in\oF(\K)$, and
$\rho_*:\Aut_\K(r)\ra\Aut_\K(f)$ is a $\K$-group morphism.
Identifying $\Aut_\K(r)$ with $\Stab_G(x)$ we have $R\subseteq
P\subseteq\Stab_G(x)=\Aut_\K(r)$, so $\rho_*\vert_R:R\ra\Aut_\K(f)$
is well-defined.

As $\nu$ is locally constructible $X^{P,R}_{\nu,c}$ is a {\it
constructible set\/} in $X^P$, and
$\overline{X^P}(\K)=\coprod_{c\in\Q} X^{P,R}_{\nu,c}$ with
$X^{P,R}_{\nu,c}\ne\emptyset$ for only finitely many $c\in\Q$. So
$X^{P,R}_{\nu,c}$ can be written as the disjoint union of finitely
many quasiprojective $\K$-varieties. But for simplicity we neglect
this, and pretend $X^{P,R}_{\nu,c}$ is a variety. Define
\e
\begin{gathered}
\hat\Pi^\nu_\F\bigl([(\R,\rho)]\bigr)\!=
\!\!\!\!\!\!\!\!\!\!\!\!\!\!\!\!\!\!\!\!\!\!\!\!\!\!\!\!\!\!
\sum_{\substack{\text{$P\!\in\!\cP(X,T^G)$, $Q\!\in\!\cQ(G,T^G)$,
$R\!\in\!\cR(X,G,T^G)$}\\
\text{and $c\!\in\!\Q:R\!\subseteq\!P\cap\!Q$,
$M^X_G(P,Q,R)\!\ne\!0$, $X^{P,R}_{\nu,c}\!\ne\!\emptyset$}}}
\begin{aligned}[t]
&M^X_G(P,Q,R)\,c\,\cdot\\
&\bigl[\bigl([X^{P,R}_{\nu,c}/C_G(Q)],\rho\!\ci\!\io^{P\cap
Q}\bigr)\bigr].
\end{aligned}
\end{gathered}
\label{mi5eq16}
\e
As for \eq{mi5eq7} we have $X^P=X^{P\cap Q}$ and $C_G(Q)=C_G(P\cap
Q)$, so $C_G(Q)$ commutes with $R\subseteq P\cap Q$, which implies
$X^{P,R}_{\nu,c}$ is $C_G(Q)$-invariant, and \eq{mi5eq16} is
well-defined. If $\nu(T,f,\phi)\equiv\mu(T)$ then $X^{P,R}_{\nu,c}$
is $X^P$ when $c=\mu(R)$ and $\emptyset$ otherwise, and \eq{mi5eq16}
reduces to \eq{mi5eq7} giving $\hat\Pi^\nu_\F=\Pi^\mu$. So
$\hat\Pi^\nu_\F$ does generalize~$\Pi^\mu$.
\label{mi5def7}
\end{dfn}

We can generalize Theorems \ref{mi5thm1} and \ref{mi5thm2}(a),(c) to
$\hat\Pi^\nu_{\smash{\F}}$ -- informally, we can generalize
\eq{mi5eq8} to regard $\hat\Pi^\nu_\F\bigl([(\R,\rho)]\bigr)$ as a
kind of double integral over $t\in T^G$ and $x\in X^{\{t\}}$ with
respect to a measure $\d\nu$ derived from $\nu$, and then the proof
of Theorem \ref{mi5thm1} needs few changes. Thus we have
well-defined linear maps $\hat\Pi^\nu_\F:\uSF(\F)\ra\uSF(\F)$
and~$\hat\Pi^\nu_\F:\SF(\F)\ra\SF(\F)$.

\subsection{The spaces $\uoSF,\oSF(\F,\Up,\La)$ and their operations}
\label{mi53}

The operators $\Pi^\mu,\Pi^\vi_n,\hat\Pi^\nu_\F$ of \S\ref{mi52}
cannot be defined on the $\uSF(\F,\Up,\La)$ of \S\ref{mi43},
basically since for $[(\R,\rho)]$ the spaces $\uSF(\F,\Up,\La)$
identify $\R=\Spec\K$ and $\R=T\t[\Spec\K/T]$ for a torus $T$, but
the $\Pi^\mu,\Pi^\vi_n,\hat\Pi^\nu_\F$ distinguish them. That is,
the relations in $\uSF(\F,\Up,\La)$ are too coarse, and identify
things separated by $\Pi^\mu,\Pi^\vi_n,\hat\Pi^\nu_\F$. We now
construct new spaces $\uoSF,\oSF (\F,\Up,\La)$ with finer relations,
on which $\Pi^\mu,\Pi^\vi_n,\hat\Pi^\nu_\F$ are well-defined.

\begin{dfn} An affine algebraic $\K$-group $G$ is called {\it very
special\/} if $C_G(Q)$ and $Q$ are special for all $Q\in\cQ(G,T^G)$,
for any maximal torus $T^G$ in $G$. (Since $Q$ is of the form
$\bG_m^k\t K$ for $K$ finite abelian, $Q$ is special if and only if
$\md{K}=1$, that is, if it is connected.) Then $G$ is special, as
$G=C_G(Q)$ for $Q=T^G\cap C(G)$. Since $\GL(k,\K)$ is special and
products of special groups are special, Example \ref{mi5ex} and
\eq{mi5eq1} imply that $\GL(m,\K)$ is very special.

When Assumption \ref{mi4ass} holds, $G$ is very special, $T^G$ is a
maximal torus in $G$ and $Q\in\cQ(G,T^G)$, define $E(G,T^G,Q)\in\La$
by
\e
E(G,T^G,Q)\!=\! \Up([Q])\!\!\!\!\!\!\!
\sum_{\substack{Q'\in\cQ(G,T^G):\\
Q\subseteq Q'}} \left\vert\frac{N_G(T^G)}{C_G(Q')\!\cap\!
N_G(T^G)}\right\vert^{-1} \!\!\cdot
\frac{n^G_{T^G}(Q',Q)}{\Up\bigl([C_G(Q')]\bigr)}\,.
\label{mi5eq17}
\e
Here $\Up([C_G(Q')])^{-1}$ exists in $\La$ by Lemma \ref{mi4lem3},
as $G$ is very special.
\label{mi5def8}
\end{dfn}

Here is our refinement of Definition~\ref{mi4def1}.

\begin{dfn} Let Assumption \ref{mi4ass} hold, and $\F$ be an
algebraic $\K$-stack with affine geometric stabilizers. Consider
pairs $(\R,\rho)$, where $\R$ is a finite type algebraic $\K$-stack
with affine geometric stabilizers and $\rho:\R\ra\F$ is a
1-morphism, with {\it equivalence} of pairs as in Definition
\ref{mi3def1}. Define $\uoSF(\F,\Up,\La)$ to be the $\La$-module
generated by equivalence classes $[(\R,\rho)]$ as above, with the
following {\it relations}:
\begin{itemize}
\setlength{\itemsep}{0pt}
\setlength{\parsep}{0pt}
\item[(i)] Given $[(\R,\rho)]$ as above and $\fS$ a closed $\K$-substack
of $\R$ we have $[(\R,\rho)]=[(\fS,\rho\vert_\fS)]+[(\R\sm\fS,
\rho\vert_{\R\sm\fS})]$, as in~\eq{mi3eq1}.
\item[(ii)] Let $\R$ be a finite type algebraic $\K$-stack with
affine geometric stabilizers, $U$ a quasiprojective $\K$-variety,
$\pi_\R:\R\t U\ra\R$ the natural projection, and $\rho:\R\ra\F$ a
1-morphism. Then~$[(\R\t U,\rho\ci\pi_\R)]=\Up([U])[(\R,\rho)]$.
\item[(iii)] Given $[(\R,\rho)]$ as above and a 1-isomorphism
$\R\cong[X/G]$ for $X$ a quasiprojective $\K$-variety and $G$ a very
special algebraic $\K$-group acting on $X$ with maximal torus $T^G$,
we have
\e
[(\R,\rho)]=\ts\sum_{Q\in\cQ(G,T^G)}E(G,T^G,Q)
\bigl[\bigl([X/Q],\rho\ci\io^Q\bigr)\bigr],
\label{mi5eq18}
\e
where $\io^Q:[X/Q]\ra\R\cong[X/G]$ is the natural projection
1-morphism.
\end{itemize}
Similarly, define $\oSF(\F,\Up,\La)$ to be the $\La$-module
generated by $[(\R,\rho)]$ with $\rho$ {\it representable}, and
relations (i)--(iii) as above. Since the $\io^Q$ are representable,
$\rho\ci\io^Q$ is representable in \eq{mi5eq18}, so these relations
make sense.

Define projections $\bar\Pi^{\Up,\La}_\F:\uSF(\F)
\ra\uoSF(\F,\Up,\La)$ and $\Pi^{\Up,\La}_\F:\uoSF(\F,\Up,\La)\ra
\uSF(\F,\Up,\La)$ by \eq{mi4eq3}. Here $\bar\Pi^{\Up,\La}_\F$ is
well-defined as relation \eq{mi3eq1} in $\uSF(\F)$ maps to (i)
above, and restricts to~$\bar\Pi^{\Up,\La}_\F:\SF(\F)\ra\oSF
(\F,\Up,\La)$.

To see $\Pi^{\Up,\La}_\F$ is well-defined we must show (i)--(iii)
above map to Definition \ref{mi4def1}(i)--(iii), which is obvious
for (i)--(ii) but nontrivial for (iii). By Definition
\ref{mi4def1}(ii),(iii) the l.h.s.\ of \eq{mi5eq18} maps under
$\Pi^{\Up,\La}_\F$ to $\Up([G])^{-1}[(X,\rho\ci\pi)]$, and the term
$[([X/Q],\rho\ci\io^Q)]$ on the r.h.s.\ maps to $\Up([Q])^{-1}
[(X,\rho\ci\pi)]$, since $Q$ is special. Therefore \eq{mi5eq18} maps
to relations in $\uSF(\F,\Up,\La)$ provided
\begin{equation*}
\Up([G])^{-1}=\ts\sum_{Q\in\cQ(G,T^G)}\Up([Q])^{-1}E(G,T^G,Q).
\end{equation*}
This follows from \eq{mi5eq4} and \eq{mi5eq17} as in the proof of
Theorem \ref{mi5thm2}(a), since $\sum_{Q\in\cQ(G,T^G):Q\subseteq
Q'}n^G_{T^G}(Q',Q)$ is 1 if $Q'=Q_{\rm min}$ and 0 otherwise, and
$C_G(Q_{\rm min})\ab=G$. Thus $\Pi^{\Up,\La}_\F$ is well-defined.
\label{mi5def9}
\end{dfn}

Here is the analogue of Definition \ref{mi4def2}, but also
including~$\Pi^\mu,\Pi^\vi_n,\hat\Pi^\nu_\F$.

\begin{dfn} Let Assumption \ref{mi4ass} hold, $\F,\G$ be
algebraic $\K$-stacks with affine geometric stabilizers, and
$\phi:\F\ra\G$ a 1-morphism. Define a $\La$-bilinear {\it
multiplication} `$\,\cdot\,$' on $\uoSF(\F,\Up,\La)$ by \eq{mi3eq2}.
This is commutative and associative as in Definition \ref{mi3def1},
and $\oSF(\F,\Up,\La)$ is closed under `$\,\cdot\,$'. Define the
{\it pushforward\/} $\phi_*:\uoSF(\F,\Up,\La)\ra \uoSF(\G,\Up,\La)$
by \eq{mi3eq6}, taking the $c_i\in\La$ rather than $c_i\in\Q$. If
$\phi$ is representable this restricts to
$\phi_*:\oSF(\F,\Up,\La)\ra\oSF(\G,\Up,\La)$. For $\phi$ of finite
type, define the {\it pullback\/} $\phi^*:\uoSF(\G,\Up,\La)\ra
\uoSF(\F,\Up,\La)$ by \eq{mi3eq7}. This restricts to $\phi^*:
\oSF(\G,\Up,\La)\ra\oSF(\F,\Up,\La)$. Define the {\it tensor
product\/} $\ot:\uoSF(\F,\Up,\La)\t\uoSF(\G,\Up,\La)
\ra\uoSF(\F\t\G,\Up,\La)$ by \eq{mi3eq8}. It restricts to
$\ot:\oSF(\F,\Up,\La)\t\oSF(\G,\Up,\La)\ra\oSF(\F\t\G,\Up,\La)$.
Define $\Pi^\mu,\Pi^\vi_n,\hat\Pi^\nu_\F:\ab\uoSF(\F,\Up,\La)
\ab\ra\ab\uoSF(\F,\Up,\La)$ and $\oSF(\F,\Up,\La)\ra
\oSF(\F,\Up,\La)$ by \eq{mi5eq7}, as in~\S\ref{mi52}.
\label{mi5def10}
\end{dfn}

Here is the analogue of Theorem~\ref{mi4thm2}.

\begin{thm} These operations `$\,\cdot\,$'$,\phi_*,\phi^*,\ot,\Pi^\mu,
\Pi^\vi_n$ and\/ $\hat\Pi^\nu_\F$ are compatible with the relations
{\rm(i)--(iii)} in Definition \ref{mi5def9}, and so are
well-defined.
\label{mi5thm3}
\end{thm}

\begin{proof} The proof of Theorem \ref{mi4thm2} shows $\phi_*$, $\ot$
are compatible with (i)--(iii) above and `$\,\cdot\,$', $\phi^*$ are
compatible with (i)--(ii). Using all the notation of Theorem
\ref{mi4thm2}, we find by the same argument that `$\,\cdot\,$' is
compatible with (iii) provided
\begin{gather*}
\bigl[\bigl(\bigl[\bigl((X_i\t Y_j)\t_{\al_{ij}\t\be_{ij},Z_{ij}\t
Z_{ij},\pi_{ij}}(Z_{ij}\t K_{ij})\bigr)/G\t
H_{ij}\bigr],\rho\ci\pi_{\R_i}\bigr)\bigr]=
\\
\sum_{Q\in\cQ(G,T^G)}E(G,T^G,Q)
\begin{aligned}[t] \bigl[\bigl(\bigl[\bigl((X_i\!\t\! Y_j)
&\t_{\al_{ij}\!\t\!\be_{ij},Z_{ij}\!\t\!Z_{ij},\pi_{ij}}
(Z_{ij}\!\t\!K_{ij})\bigr)\\
&/Q\!\t\!H_{ij}\bigr],
\rho\!\ci\!\pi_{\R_i}\!\ci\!\pi_{ij}\bigr)\bigr]
\end{aligned}
\end{gather*}
in $\uoSF(\F,\Up,\La)$. This holds because by (iii), both sides are
equal to
\begin{equation*}
\sum_{\substack{Q\in\cQ(G,T^G)\\ Q'\in\cQ(H_{ij},T^{H_{ij}})}}
\!\!\!\!\!\!\!
\begin{aligned}[t]
E(G,T^G,Q)E(H_{ij},T^{H_{ij}},Q')&\bigl[\bigl(\bigl[\bigl((X_i\!\t\!
Y_j)\t_{\al_{ij}\!\t\!\be_{ij},Z_{ij}\!\t\!Z_{ij},\pi_{ij}}\\
&(Z_{ij}\!\t\!K_{ij})\bigr)
/Q\!\t\!Q'\bigr],\rho\!\ci\!\pi_{\R_i}\!\ci\!\pi_{ij}\bigr)\bigr].
\end{aligned}
\end{equation*}
Here we use the facts that $\cQ(G\!\t\!H_{ij},T^G\!\t\!T^{H_{ij}})
\!=\!\cQ(G,T^G)\!\t\!\cQ(H_{ij},T^{H_{ij}})$, and
$E(G\!\t\!H_{ij},T^G\!\t\!T^{H_{ij}},Q\t Q')=E(G,T^G,Q)
E(H_{ij},T^{H_{ij}},Q')$, and each $Q\!\in\!\cQ(G,T^G)$ is a torus,
so $\cQ(Q,Q)\!=\!\{Q\}$ and $E(Q,Q,Q)=1$. Therefore `$\,\cdot\,$' is
compatible with (iii) and is well-defined. Modifying the argument of
Theorem \ref{mi4thm2} in the same way, $\phi^*$ is well-defined.

Compatibility of $\Pi^\mu$ with (i)--(ii) above is easy. To show
$\Pi^\mu$ is compatible with (iii) we must show it takes both sides
of \eq{mi5eq18} to the same thing in $\uoSF(\F,\Up,\La)$. That is,
we must prove that
\e
\begin{gathered}
\sum_{\substack{\text{$P\in\cP(X,T^G)$, $Q\in\cQ(G,T^G)$ and}\\
\text{$R\!\in\!\cR(X,G,T^G):R\!\subseteq\!P\cap Q$,
$M^X_G(P,Q,R)\!\ne\!0$}}}
\begin{aligned}[t]
&M^X_G(P,Q,R)\mu(R)\,\cdot\\
&\bigl[\bigl([X^P/C_G(Q)],\rho\!\ci\!\io^{P\cap Q}\bigr)\bigr]=
\end{aligned}
\\
\sum_{\substack{\text{$Q'\!\in\!\cQ(G,T^G)$,
$P'\!\in\!\cP(X,Q')$, and}\\
\text{$R\!\in\!\cR(X,Q',Q'):R\!\subseteq\!P'$,
$M^X_{Q'}(P',Q',R)\!\ne\!0$}}}
\!\!\!\!\!\!\!\!\!\!\!\!\!\!\!\!\!\!\!\!\!
\begin{aligned}[t]
E(G,T^G,Q')\,\cdot\, & M^X_{Q'}(P',Q',R)\mu(R)\,\cdot\\
&\bigl[\bigl([X^{P'}/Q'],\rho\!\ci\!\io^{P'\cap Q'}\bigr)\bigr]
\end{aligned}
\end{gathered}
\label{mi5eq19}
\e
in $\uoSF(\F,\Up,\La)$, using $\cQ(Q',Q')=\{Q'\}$ and
$C_{Q'}(Q')=Q'$ in the bottom line.

We rewrite the top line of \eq{mi5eq19} using \eq{mi5eq18}. Since
$\cQ(C_G(Q),T^G)=\{Q'\in\cQ(G,T^G):Q\subseteq Q'\}$ for $Q\subseteq
Q'$, this gives
\e
\begin{gathered}
\sum_{\substack{\text{$P\!\in\!\cP(X,T^G)$, $Q\!\in\!\cQ(G,T^G)$}\\
\text{and $R\!\in\!\cR(X,G,T^G):$}\\ \text{$R\subseteq P\cap Q$,
$M^X_G(P,Q,R)\ne 0$}}} \!\!\!\!\!\!\!\!\!\!\!\!\!\!\!
\!\!\!\!\!\!\!\!\!\!\!\!\ M^X_G(P,Q,R)\mu(R)\cdot \!\!\!\!\!\!\!\!
\sum_{\substack{Q'\in\cQ(G,T^G):\\ Q\subseteq Q'}}\!
\begin{aligned}[t]
& E(C_G(Q),T^G,Q')\,\cdot \\
&\bigl[\bigl([X^P/Q'],\rho\!\ci\!\io^{P\cap Q'}\bigr)\bigr].
\end{aligned}
\end{gathered}
\label{mi5eq20}
\e

We claim that the term in the bottom line of \eq{mi5eq19} with fixed
$P',Q',R$ agrees with the sum of terms in \eq{mi5eq20} with fixed
$P,Q',R$, where $P'=P\cap Q'$. To explain the relation between $P$
and $P'$, note that for $P,Q,R,Q'$ in \eq{mi5eq20} we have
$M^X_G(P,Q,R)\ne 0$, so $P$ is the smallest element of $\cP(X,T^G)$
containing $P\cap Q$ by Lemma \ref{mi5lem3}, and $X^{P\cap Q}=X^P$.
But $Q\subseteq Q'$, so $P\cap Q\subseteq P\cap Q'=P'\subseteq P$,
which shows that $X^{P'}=X^P$. Note too that $P$ is the smallest
element of $\cP(X,T^G)$ containing $P'$, so $P$ and $P'$ determine
each other uniquely given $Q'$, and fixing $P',Q',R$ in \eq{mi5eq19}
is equivalent to fixing $P,Q',R$ in~\eq{mi5eq20}.

Thus $\mu(R)\cdot[([X^P/Q'],\rho\ci\io^{P\cap Q'})]$ are common
terms in \eq{mi5eq19} and \eq{mi5eq20}, and the sums of coefficients
of these for fixed $P,P',Q',R$ are equal provided
\e
E(G,T^G,Q')M^X_{Q'}(P',Q',R)=E(C_G(Q),T^G,Q') \!\!\!\!\!\!\!\!\!\!
\sum_{Q\in\cQ(G,T^G):Q\subseteq Q'} \!\!\!\!\!\!\!\!\!\! M^X_G(P,Q,R).
\label{mi5eq21}
\e
Now for $Q,Q''\in\cQ(G,T^G)$ with $Q\subseteq Q'\subseteq Q''$ we
have
\begin{equation*}
\left\vert\frac{N_G(T^G)}{C_G(Q'')\!\cap\! N_G(T^G)}\right\vert=
\left\vert\frac{N_{C_G(Q)}(T^G)}{C_{C_G(Q)}(Q'')\!\cap\!
N_{C_G(Q)}(T^G)}\right\vert\cdot
\left\vert\frac{N_G(T^G)}{C_G(Q)\!\cap\!N_G(T^G)}\right\vert,
\end{equation*}
noting that $C_{C_G(Q)}(Q'')=C_G(Q'')$ as $Q\subseteq Q''$ and
intersecting top and bottom of the l.h.s.\ with $C_G(Q)$. Thus from
\eq{mi5eq17} we deduce that
\e
E(C_G(Q),T^G,Q')=\left\vert\frac{N_G(T^G)}{C_G(Q)\cap
N_G(T^G)}\right\vert\cdot E(G,T^G,Q').
\label{mi5eq22}
\e

Combining this with \eq{mi5eq5} we see that \eq{mi5eq21} is
equivalent to
\e
m_{Q'}^X(R,P')=
\sum_{\substack{Q\in\cQ(G,T^G):\\ Q\subseteq Q'}}\,\,
\sum_{\begin{subarray}{l} \text{$\hat P\in\cP(X,T^G)$,
$\hat Q\in\cQ(G,T^G)$:}\\
\text{$\hat P\subseteq P$, $\hat Q\subseteq Q$, $R=\hat P\cap\hat
Q$}
\end{subarray}\!\!\!\!\!\!\!\!\!\!\!\!\!\!\!\!\!\!\!\!\!\!\!\! }
\!\!\!\!\!\!\!\!\!\!\!\!\!\!\!\! m^X_{T^G}(\hat P,P)n^G_{T^G}(\hat
Q,Q),
\label{mi5eq23}
\e
where the l.h.s.\ is $M^X_{Q'}(P',Q',R)$. Fixing $\hat Q\subseteq
Q'$ in the r.h.s.\ of \eq{mi5eq23} and summing over $Q$, we find
$\sum_{Q\in\cQ(G,T^G):\hat Q\subseteq Q\subseteq Q'}n^G_{T^G}(\hat
Q,Q)$ is 1 if $\hat Q=Q'$ and 0 otherwise. So the r.h.s.\ becomes
$\sum_{\hat P\in\cP(X,T^G):\hat P\subseteq P,\; R=\hat P\cap
Q'}m^X_{T^G}(\hat P,P)$, which eventually reduces to
$m_{Q'}^X(R,P')$ as $P$ is the smallest element of $\cP(X,T^G)$
containing $P'=P\cap Q'$. This proves \eq{mi5eq23}, and hence
\eq{mi5eq21} and \eq{mi5eq19}, which shows $\Pi^\mu$ is compatible
with (iii) and is well-defined. Also $\Pi^\vi_n$ is a special case
of $\Pi^\mu$, and the changes to the proof above for
$\hat\Pi^\nu_\F$ are straightforward.
\end{proof}

The analogue of Corollary \ref{mi4cor} is immediate.

\begin{cor} The projections $\bar\Pi^{\Up,\La}_*,\Pi^{\Up,\La}_*$
commute with\/ `$\,\cdot\,$'$,\ab\phi_*,\ab\phi^*,\ab \ot$ on
$\uSF(*),\uoSF,\uSF(*,\Up,\La)$, and\/ $\bar\Pi^{\Up,\La}_*$
commutes with\/ $\Pi^\mu,\Pi^\vi_n,\hat\Pi^\nu_\F$. The analogues of
Theorems \ref{mi3thm1} and \ref{mi5thm2} and Proposition
\ref{mi5prop2} hold for~$\uoSF,\oSF(*,\Up,\La)$.
\label{mi5cor}
\end{cor}

Here is a useful way of representing elements
of~$\uoSF,\oSF(\F,\Up,\La)$.

\begin{prop} $\uoSF(\F,\Up,\La)$ and\/ $\oSF(\F,\Up,\La)$ are
generated over $\La$ by elements $[(U\t[\Spec\K/T],\rho)]$, for\/
$U$ a quasiprojective $\K$-variety and\/ $T$ an algebraic $\K$-group
isomorphic to $\bG_m^k\t K$ for $k\ge 0$ and\/ $K$ finite abelian.
\label{mi5prop3}
\end{prop}

\begin{proof} As in the proof of Theorem \ref{mi4thm1},
$\uoSF,\oSF(\F,\Up,\La)$ are generated over $\La$ by elements
$[([X/G],\rho)]$ for $X$ a quasiprojective $\K$-variety and $G$ an
affine algebraic $\K$-group which we can take to be $\GL(m,\K)$, so
in particular for $G$ {\it very special}. Definition
\ref{mi5def9}(iii) then implies $\uoSF,\oSF(\F,\Up,\La)$ are
generated over $\La$ by $[([X/Q],\rho\ci\io^Q)]$ for $X$ a
quasiprojective $\K$-variety and $Q$ a torus.

Given such $X,Q$ there is a finite collection of closed
$\K$-subgroups $T_i$ in $Q$ for $i\in I$ occurring as $\Stab_Q(x)$
for $x\in X$, and the set of such $x$ is a locally closed
$\K$-subvariety $X_i$ of $X$ with $X=\coprod_{i\in I}X_i$. Here
$T_i\cong\bG_m^{k_i}\t K_i$ for $K_i$ finite abelian, as $Q$ is a
torus. Then $Q/T_i$ acts freely on $X_i$, and $X_i/(Q/T_i)$ is an
algebraic $\K$-space which may be written as a disjoint union of
finitely many quasiprojective $\K$-subvarieties $U_{ij}$ for $j\in
J_i$. Thus $[X/Q]$ is a disjoint union of $\K$-substacks
1-isomorphic to $U_{ij}\t[\Spec\K/T_i]$. Definition \ref{mi5def9}(i)
gives $[([X/Q],\rho\ci\io^Q)]=\sum_{i,j}[(U_{ij}\t[\Spec\K/T_i],
\rho_{ij})]$, so $\uoSF,\oSF(\F,\Up,\La)$ are generated over $\La$
by such $[(U_{ij}\t[\Spec\K/T_i],\rho_{ij})]$, as we want.
\end{proof}

Such $[(U\t[\Spec\K/T],\rho)]$ are linearly independent for
nonisomorphic~$T$.

\begin{prop} Suppose $\sum_{i\in I}c_i[(U_i\t[\Spec\K/T_i],\rho_i)]=0$
in $\uoSF(\F,\Up,\La)$, where $I$ is finite set, $c_i\in\La$, $U_i$
a quasiprojective $\K$-variety and\/ $T_i$ an algebraic $\K$-group
isomorphic to $\bG_m^{k_i}\t K_i$ for $k_i\ge 0$ and\/ $K_i$ finite
abelian, with\/ $T_i\not\cong T_j$ for $i\ne j$. Then
$c_j[(U_j\t[\Spec\K/T_j],\rho_j)]=0$ for all\/~$j\in I$.
\label{mi5prop4}
\end{prop}

\begin{proof} Let $j\in I$, and define $\mu$ in Definition
\ref{mi5def5} by $\mu(T)=1$ if $T\cong T_j$ and $\mu(T)=0$
otherwise. Then $\Pi^\mu$ is well-defined on $\uoSF(\F,\Up,\La)$ by
Theorem \ref{mi5thm3}, and on $[(U_i\t[\Spec\K/T_i],\rho_j)]$ it is
the identity if $i=j$ and 0 otherwise, since $T_i\not\cong T_j$ for
$i\ne j$. The result follows by applying $\Pi^\mu$ to~$\sum_{i\in
I}\cdots=0$.
\end{proof}

We identify $\uoSF,\oSF(\Spec\K,\Up,\La)$, as in
Proposition~\ref{mi4prop2}.

\begin{prop} Define a commutative $\La$-algebra $\bar\La$ with\/
$\La$-basis isomorphism classes $[T]$ of\/ $\K$-groups $T$ of the
form $\bG_m^k\t K$, for $k\ge 0$ and\/ $K$ finite abelian, with
$\La$-bilinear multiplication given by $[T][T']=[T\t T']$ on basis
elements. Define $\bar\imath_\La:\bar\La\ra\uoSF(\Spec\K,\Up,\La)$
by $\bar\imath_\La:\sum_ic_i[T_i]\mapsto\sum_ic_i[[\Spec\K/T_i]]$.
Then $\bar\imath_\La$ is an algebra isomorphism. It restricts to an
isomorphism from the subalgebra $\La[\{1\}]\cong\La$ in $\bar\La$
to~$\oSF(\Spec\K,\Up,\La)$.
\label{mi5prop5}
\end{prop}

\begin{proof} Since $[[\Spec\K/T]]\cdot[[\Spec\K/T']]=[[\Spec\K/T\t
T']]$ in $\uoSF(\Spec\K,\Up,\La)$, $\bar\imath_\La$ is an {\it
algebra morphism}. By Proposition \ref{mi5prop3},
$\uoSF(\Spec\K,\Up,\La)$ is generated over $\La$ by
$[U\t[\Spec\K/T]]$ with $U$ a quasiprojective $\K$-variety and
$T\cong\bG_m^k\t K$. But $[U\t[\Spec\K/T]]=\Up([U])[[\Spec\K/T]]$ by
Definition \ref{mi5def9}(ii), so $\uoSF(\Spec\K,\Up,\La)$ is
generated over $\La$ by such $[[\Spec\K/T]]$, and $\bar\imath_\La$
is {\it surjective}.

Elements of $\bar\La$ may be written as $\sum_{i\in I}c_i[T_i]$ for
$I$ finite, $c_i\in\La$ and $T_i\not\cong T_j$ for $i\ne j$. Suppose
$\bar\imath_\La\bigl(\sum_{i\in I}c_i[T_i]\bigr)=0$. Then
$\sum_{i\in I}c_i[[\Spec\K/T_i]]=0$, so Proposition \ref{mi5prop4}
gives $c_i[[\Spec\K/T_i]]=0$ for each $i\in I$. Applying
$\Pi^{\Up,\La}_{\Spec\K}$ from Definition \ref{mi5def9} and
$i_\La^{-1}$ from Proposition \ref{mi4prop2} gives
\begin{equation*}
0=i_\La^{-1}\ci\Pi^{\Up,\La}_{\Spec\K}\bigl(c_i[[\Spec\K/T_i]]
\bigr)=c_i i_\La^{-1}\bigl([[\Spec\K/T_i]]\bigr)=c_i(\ell-1)^{-\dim
T_i}
\end{equation*}
in $\La$. So $c_i=0$ for all $i\in I$, and $\sum_{i\in
I}c_i[T_i]=0$. Thus $\bar\imath_\La$ is {\it injective}, and so an
isomorphism. Finally, as $\rho:[\Spec\K/T]\ra\Spec\K$ is
representable if and only if $T\cong\{1\}$, $\oSF(\Spec\K,\Up,\La)$
is the image under $\bar\imath_\La$ of~$\La[\{1\}]$.
\end{proof}

Since every finite abelian group $K$ is isomorphic to a product of
cyclic groups $\Z_{p^k}$ of prime power order, $\bar\La$ is the free
commutative $\La$-algebra generated by $[\bG_m]$ and $[\Z_{p^k}]$
for $p$ prime and $k\ge 1$. The proof of Proposition \ref{mi4prop3}
gives:

\begin{prop} The following maps are $\La$-linear and injective:
\begin{align*}
(\bar\Pi^{\Up,\La}_\F\ci\io_\F)\ot_\Q\id_\La&:
\CF(\F)\ot_\Q\La\longra\oSF(\F,\Up,\La),\\
\mu\ci\bigl((\bar\Pi^{\Up,\La}_\F\ci\io_\F)\ot_\Q
\bar\imath_\La\bigr)&:\CF(\F)\ot_\Q\bar\La\longra\uoSF(\F,\Up,\La),
\end{align*}
where in the second line $\mu:\uoSF(\Spec\K,\Up,\La)\ot_\Q
\uoSF(\F,\Up,\La)\ra\uoSF(\F,\Up,\La)$ is the combination of the
tensor product\/ $\ot:\uoSF(\Spec\K,\Up,\La)\t
\uoSF(\F,\Up,\La)\ra\uoSF(\Spec\K\t\F,\Up,\La)$ of Definition
\ref{mi5def10} with the isomorphism~$\Spec\K\t\F\cong\F$.
\label{mi5prop6}
\end{prop}

Again, this shows that the spaces $\uoSF,\oSF(\F,\Up,\La)$ are quite
large, though not as large as $\uSF,\SF(\F)\ot_\Q\La$, and therefore
that the relations Definition \ref{mi5def9}(i)--(iii) have some kind
of consistency about them.

Given a generator $[(\R,\rho)]$, for each $r\in\oR(\K)$ with
$\rho_*(r)=x\in\oF(\K)$ we have a $\K$-group morphism
$\rho_*:\Aut_\K(r)\ra\Aut_\K(x)$. Roughly speaking, the difference
between the spaces $\uSF(\F,\Up,\La)$ of \S\ref{mi43} and the
$\uoSF,\oSF(\F,\Up,\La)$ above is that the $\uoSF,\oSF(\F,\Up,\La)$
keep track of the restriction of $\rho_*$ to a maximal torus of
$\Aut_\K(r)$, but $\uSF(\F,\Up,\La)$ loses this information.

Proposition \ref{mi5prop3} shows that in $\uoSF,\oSF(\F,\Up,\La)$ we
can always reduce to $[(\R,\rho)]$ with all stabilizer groups
$\Aut_\K(r)$ for $r\in\oR(\K)$ of the form $\bG_m^k\t K$ for $K$
finite abelian. That is, $\uoSF,\oSF(\F,\Up,\La)$ {\it abelianize
stabilizer groups}. We can regard the $\Pi^\vi_n$ of \S\ref{mi52} as
doing the same job: $[[\Spec\K/G]]$ is the sum of components
$\Pi^\vi_n([[\Spec\K/G]])$ which behave like multiples
of~$[[\Spec\K/\bG_m^n]]$.

In the applications of \cite{Joyc3,Joyc4,Joyc5}, the concept of {\it
virtual rank\/} given by the $\Pi^\vi_n$ is more useful than the
{\it real rank\/} given by the $\Pi^\re_n$ of \S\ref{mi51}. One
reason for this is that the author can prove there are no relations
analogous to Definition \ref{mi5def9}(i)--(iii) which are compatible
with the $\Pi^\re_n$ in the way that these are compatible with the
$\Pi^\vi_n$, and which are also compatible with multiplication
`$\,\cdot\,$' and pullbacks $\phi^*$. This suggests there is some
kind of consistency between the $\Pi^\mu,\Pi^\vi_n,\hat\Pi^\nu_\F$
and `$\,\cdot\,$'$,\phi^*$ which the author does not yet understand.

\section{Extension to the case $\ell=1$}
\label{mi6}

We now extend the constructions of \S\ref{mi4}--\S\ref{mi5} to the
case when $\ell-1$ is not invertible in $\La$, and in particular to
the case $\ell=1$, which includes Euler characteristics $\chi$. We
do this in \S\ref{mi61} by supposing the algebra $\La$ of
\S\ref{mi41} has a subalgebra $\La^\ci$ containing $\Up([X])$ for
varieties $X$ and some rational functions of $\ell$, but not
$(\ell-1)^{-1}$, and that we are given a surjective algebra morphism
$\pi:\La^\ci\ra\Om$ with $\pi(\ell)=1$. Then $\Th=\pi\ci\Up$ is the
motivic invariant we are interested in, which takes values in $\Om$.
This can be done in all our examples.

Section \ref{mi62} shows that the coefficients $E(G,T^G,Q)$ of
\eq{mi5eq17} actually lie in $\La^\ci$ (this is not obvious), and
computes them when $G=\GL(m,\K)$. Therefore the relations Definition
\ref{mi5def9}(i)--(iii) for $\uSF,\SF(\F,\Up,\La)$ make sense with
coefficients in $\La^\ci$ rather than in $\La$, and applying $\pi$
they also make sense with coefficients in $\Om$. So in \S\ref{mi63}
we define new spaces $\uoSF,\oSF(\F,\Up,\La^\ci)$ and
$\uoSF,\oSF(\F,\Th,\Om)$ with these relations, with the usual
operations `$\,\cdot\,$'$,\phi_*,\phi^*,\ot,\Pi^\mu,\Pi^\vi_n,
\hat\Pi^\nu_\F$. These will be important in \cite{Joyc3,Joyc4,Joyc5}
for defining invariants counting coherent sheaves on Calabi--Yau
3-folds. When $\Om=\Q$ and $\Th=\chi$ we also define smaller spaces
$\uhSF,\hSF(\F,\chi,\Q)$ exploiting special properties of $\chi$ on
fibrations.

\subsection{Initial assumptions and examples}
\label{mi61}

For our next constructions we need more data than in
Assumption~\ref{mi4ass}.

\begin{ass} Suppose Assumption \ref{mi4ass} holds, and
$\La^\ci$ is a $\Q$-subalgebra of $\La$ containing the image of
$\Up$ and the elements $\ell^{-1}$ and
$(\ell^k+\ell^{k-1}+\cdots+1)^{-1}$ for $k=1,2,\ldots$, but {\it
not\/} containing $(\ell-1)^{-1}$. Let $\Om$ be a commutative
$\Q$-algebra, and $\pi:\La^\ci\ra\Om$ a surjective $\Q$-algebra
morphism, such that $\pi(\ell)=1$. Define
\begin{equation*}
\Th:\{\text{isomorphism classes $[X]$ of quasiprojective
$\K$-varieties $X$}\}\longra\Om
\end{equation*}
by $\Th=\pi\ci\Up$. Then~$\Th([\bA^1])=1$.
\label{mi6ass}
\end{ass}

Note that given $\La^\ci$ satisfying the above conditions, there is
a natural choice for $\Om,\pi$: by assumption $(\ell-1)\La^\ci$ is
an {\it ideal\/} in $\La^\ci$, not containing 1, so we may take
$\Om=\La^\ci/(\ell-1)\La^\ci$ to be the quotient algebra, with
projection $\pi:\La^\ci\ra\Om$. We can satisfy Assumption
\ref{mi6ass} in all the examples of~\S\ref{mi41}.

\begin{ex} In Example \ref{mi4ex1}, let $\La^\ci_\Ho$ be the
subalgebra of $P(x,y)/Q(x,y)$ in $\La_\Ho$ for which $xy-1$ does not
divide $Q(x,y)$. Set $\Om_\Ho=\Q(x)$, the $\Q$-algebra of rational
functions in $x$, and define $\pi_\Ho:\La^\ci_\Ho\ra\Om_\Ho$ by
$\pi_\Ho:P(x,y)/Q(x,y)\mapsto P(x,x^{-1})/Q(x,x^{-1})$. Then
Assumption \ref{mi6ass} holds.
\label{mi6ex1}
\end{ex}

\begin{ex} In Example \ref{mi4ex2}, let $\La^\ci_\Po$ be the
subalgebra of $P(z)/Q(z)$ in $\La_\Po$ for which $z\pm 1$ do not
divide $Q(z)$. Here are three possibilities for~$\Om_\Po,\pi_\Po$:
\begin{itemize}
\setlength{\itemsep}{0pt}
\setlength{\parsep}{0pt}
\item[(a)] Set $\Om_\Po=\Q$ and $\pi_\Po:f(z)\mapsto f(-1)$. Then
$\Th_\Po([X])=\pi_\Po\ci\Up_\Po([X])$ is the {\it Euler
characteristic} of~$X$.
\item[(b)] Set $\Om_\Po=\Q$ and $\pi_\Po:f(z)\mapsto f(1)$. Then
$\Th_\Po([X])=\pi_\Po\ci\Up_\Po([X])$ is the {\it sum of the virtual
Betti numbers} of~$X$.
\item[(c)] Set $\Om_\Po=\Q\op\Q$, a product of algebras, and
$\pi_\Po:f(z)\mapsto(f(-1),f(1))$. This combines (a) and~(b).
\end{itemize}
Assumption \ref{mi6ass} holds in each case.
\label{mi6ex2}
\end{ex}

\begin{ex} In Example \ref{mi4ex3}, let $\La^\ci_\uni$ be the
subalgebra of $\La_\uni$ generated by elements $\ell^{-1}$,
$(\ell^k+\ell^{k-1}+\cdots+1)^{-1}$ for $k=1,2,\ldots$, and $[X]$
for quasiprojective $\K$-varieties $X$. Define $\Om_\uni$ to be the
quotient algebra $\La^\ci_\uni/(\ell-1)\La^\ci_\uni$, with
projection $\pi_\uni:\La^\ci_\uni\ra\Om_\uni$. We have a morphism
$\La_\uni\ra\La_\Po$ taking $\La^\ci_\uni\ra\La^\ci_\Po$ and
$\ell\mapsto\ell$, so $(\ell-1)^{-1}\notin\La^\ci_\uni$ as
$(\ell-1)^{-1}\notin\La^\ci_\Po$ in Example \ref{mi6ex2}, and
Assumption \ref{mi6ass} holds.
\label{mi6ex3}
\end{ex}

Here are some useful facts about~$\Up([G])$.

\begin{lem} Let Assumptions \ref{mi4ass} and \ref{mi6ass} hold,
and\/ $G$ be a special algebraic $\K$-group of rank\/ $k$. Then
$\Up([G])\in(\ell-1)^k\La^\ci$ and\/~$\Up([G])^{-1}\in
(\ell-1)^{-k}\La^\ci$.
\label{mi6lem1}
\end{lem}

\begin{proof} Let $T^G$ be a maximal torus in $G$, so that
$T^G\cong\bG_m^k$. Then $G\ra G/T^G$ is a principal bundle with
fibre $T^G$, which is special. So by Definition \ref{mi2def1} it is
a Zariski locally trivial fibration, and
$\Up([G])=\Up([T^G])\Up([G/T^G])$ by Lemma \ref{mi4lem1}. Hence
$\Up([G])\in (\ell-1)^k\La^\ci$, as~$\Up([T^G])=(\ell-1)^k$.

As $G$ is special we can embed it in $\GL(m,\K)$ with
$\GL(m,\K)\ra\GL(m,\K)/G$ a Zariski locally trivial fibration, so
$\Up([\GL(m,\K)])=\Up([G])\Up([\GL(m,\K)/G])$. Applying Lemma
\ref{mi4lem2} yields
\e
\Up([G])^{-1}\!=\!(\ell\!-\!1)^{-m}\cdot\ell^{-m(m-1)/2}
\ts\prod_{k=1}^m(\ell^{k-1}\!+\!\cdots\!+\!1)^{-1}\cdot
\Up([\GL(m,\K)/G]).
\label{mi6eq1}
\e
Now the diagonal matrices $\bG_m^m$ in $\GL(m,\K)$ act on
$\GL(m,\K)/G$, and the stabilizer of each point is isomorphic to
$\bG_m^j\t K$ for $j\ge 0$ and finite abelian $K$. Since this is
conjugate to a subgroup of $G$, we have $j\le k$. Hence each orbit
of $\bG_m^m$ is isomorphic $\bG_m^{m-j}$ for some $0\le j\le k$.

Thus we may write $\GL(m,\K)/G$ as a finite disjoint union of
$\bG_m^m$-invariant $\K$-subvarieties $X_i$, such that the
$\bG_m^m$-orbits make $X_i$ into a fibre bundle with fibre
$\bG_m^{m-j_i}$ for $0\le j_i\le k$. Now $\bG_m^n$-bundles are
Zariski locally trivial fibrations as $\bG_m^n$ is special. So
refining the decomposition if necessary we can suppose $X_i\cong
\bG_m^{m-j_i}\t Y_i\cong\bG_m^{m-k}\t\bG_m^{k-j_i}\t Y_i$ for some
quasiprojective $\K$-varieties $Y_i$. Thus
$\Up([X_i])=(\ell-1)^{m-k}\Up([\bG_m^{k-j_i}\t Y_i])$, and
$\Up([\GL(m,\K)/G])=\sum_i\Up([X_i])\in(\ell-1)^{m-k}\La^\ci$.
Combining this with \eq{mi6eq1} and using Assumption \ref{mi6ass}
shows~$\Up([G])^{-1}\in(\ell-1)^{-k}\La^\ci$.
\end{proof}

\subsection{Properties of the $E(G,T^G,Q)$}
\label{mi62}

We shall show $E(G,T^G,Q)$ in \eq{mi5eq17} lies in $\La^\ci$. This
is far from obvious, as by Lemma \ref{mi6lem1} each term in
\eq{mi5eq17} lies in $(\ell-1)^{\dim Q-\rk G}\La^\ci$. Effectively,
in the sum over $Q'$ in \eq{mi5eq17} the terms in $(\ell-1)^{-n}$
for $0<n\le\rk\,G-\dim Q$ all cancel.

\begin{thm} Let Assumptions \ref{mi4ass} and \ref{mi6ass} hold
and\/ $G$ be a very special algebraic $\K$-group with maximal torus
$T^G$. Then $E(G,T^G,Q)$ in \eq{mi5eq17} lies in $\La^\ci$ for
all\/~$Q\in\cQ(G,T^G)$.
\label{mi6thm1}
\end{thm}

\begin{proof} Let $X$ be a quasiprojective $\K$-variety acted on by
a very special $\K$-group $G$, so that $[[X/G]]\in\uSF(\Spec\K)$.
For $\mu$ a {\it weight function}, Definition \ref{mi5def5} defines
$\Pi^\mu:\uSF(\Spec\K)\ra\uSF(\Spec\K)$. Applying $\Up'$ of Theorem
\ref{mi4thm1} to \eq{mi5eq7} and noting that $C_G(Q)$ is special for
all $Q\in\cQ(G,T^G)$ as $G$ is very special yields
\e
\begin{gathered}
\Up'\ci\Pi^\mu\bigl([X/G]\bigr)\!=
\!\!\!\!\!\!\!\!\!\!\!\!\!\!\!\!\!\!\!\!\!\!\!
\sum_{\substack{\text{$P\in\cP(X,T^G)$, $Q\in\cQ(G,T^G)$ and}\\
\text{$R\!\in\!\cR(X,G,T^G):R\!\subseteq\!P\cap Q$,
$M^X_G(P,Q,R)\!\ne\!0$}}}
\begin{aligned}[t]
&M^X_G(P,Q,R)\mu(R)\,\cdot\\
&\Up([X^P])\Up([C_G(Q)])^{-1}.
\end{aligned}
\end{gathered}
\label{mi6eq2}
\e
Substituting in \eq{mi5eq5} gives a sum over $P,Q,R,P',Q'$ with
$R=P'\cap Q'$. Comparing this with \eq{mi5eq17} we see that the sum
over $Q$ is proportional to that defining $E(G,T^G,Q')$, so
\eq{mi6eq2} becomes
\e
\begin{gathered}
\Up'\ci\Pi^\mu\bigl([X/G]\bigr)= \!\!\!\!\!\!
\sum_{\substack{\text{$P'\subseteq P$ in $\cP(X,T^G)$,}\\
\text{$Q'$ in $\cQ(G,T^G)$}}} \!\!\!\!\!\!
\begin{aligned}[t]
m_{T^G}^X(P',P)\Up([X^P])\mu(P'\cap Q')\,\cdot &\\
\Up([Q'])^{-1}E(G,T^G,Q')&.
\end{aligned}
\end{gathered}
\label{mi6eq3}
\e

If $P$ is a closed $\K$-subgroup of $T^G$, write $X_{T^G}^P=\{x\in
X:\Stab_{T^G}(x)=P\}$, a subvariety of $X$. It is easy to see that
if $X_{T^G}^P\ne\emptyset$ then $P\in\cP(X,T^G)$, and for
$P\in\cP(X,T^G)$ we have $X^P=\coprod_{P'\in\cP(X,T^G):P'\subseteq
P}X^{P'}_{T^G}$. Therefore $\Up([X^P])=\sum_{P'\in\cP(X,T^G):P'
\subseteq P}\Up([X^{P'}_{T^G}])$. Inverting this combinatorially
using properties of the $m_{T^G}^X(P',P)$ yields
$\Up([X^{P'}_{T^G}])=\sum_{P\in\cP(X,T^G):P'\subseteq
P}m_{T^G}^X(P',P)\Up([X^P])$. Comparing this with \eq{mi6eq3} we see
that
\e
\Up'\ci\Pi^\mu\bigl([X/G]\bigr)= \!\!\!\!\!\!\!\!\!
\sum_{\text{$P'\in\cP(X,T^G)$, $Q'\in\cQ(G,T^G)$}
\!\!\!\!\!\!\!\!\!\!\!\!\!\!\!\!\!\!\!\!\!\!}\!\!\!\!\!\!\!\!\!
\Up([X_{T^G}^{P'}])\mu(P'\cap Q') \Up([Q'])^{-1}E(G,T^G,Q').
\label{mi6eq4}
\e

Now $E(G,T^G,Q)\in(\ell-1)^{\dim Q-\rk\,G}\La^\ci$ as above, proving
the theorem for $Q=T^G$. So suppose $\dim Q<\dim T^G=\rk\,G$. Choose
a $\K$-subgroup $T\subseteq T^G$ with $T\cong\bG_m^{\rk\,G-\dim Q}$
such that $K=T\cap Q$ is finite, and if $Q\ne Q'\in \cQ(G,T^G)$ then
$T\cap Q'\not\cong K$. This is possible if $\dim Q>0$, as there are
infinitely many $T$. But it may not be if $\dim Q=0$, as $T=T^G$ is
the only choice.

Define $\mu$ in Definition \ref{mi5def5} by $\mu([H])=1$ if $H\cong
K$ and $\mu([H])=0$ otherwise. Set $X=G/T$. Then we have
1-isomorphisms $[X/G]\cong[\Spec\K/T]\cong[(T^G/T)/T^G]$. Since
$T\not\cong K$ as $\dim T>0$, we find from Definition \ref{mi5def5}
that $\Pi^\mu([(T^G/T)/T^G])=0$, so $\Pi^\mu([X/G])=0$ by Theorem
\ref{mi5thm1}, and \eq{mi6eq4} is zero. Suppose some $P',Q'$ give a
nonzero term on the r.h.s.\ of \eq{mi6eq4}. Then $P'$ is conjugate
in $G$ to a subgroup of $T$ as $X^{P'}_{T^G}\ne\emptyset$, and
$P'\cap Q'\cong K$ as $\mu(P'\cap Q')\ne 0$. Hence $\dim P'\le\dim
T$, and $\dim P'+\dim Q'\le\rk\,G$ as~$\dim P'\cap Q'=0$.

If $\dim P'=\dim T$ then $P'$ is conjugate to $T$ as $T$ is
connected, giving $P'=\ga T$ for $\ga\in W(G,T^G)$. Then $P'\cap
Q'\cong K$ and the choice of $T$ imply $Q'=\ga Q$. Rearranging
\eq{mi6eq4} to put terms $P',Q'=\ga T,\ga Q$ on the left gives
\begin{gather}
\frac{\md{W(G,T^G)}}{\md{\{\ga\in W(G,T^G):\ga T=T,\;\ga Q=Q\}}}
\Up\bigl([{\ts X_{T^G}^T}\,]\bigr)\Up([Q])^{-1}E(G,T^G,Q)=
\nonumber\\
- \!\!\!\!\!\!\!
\sum_{\substack{P'\in\cP(X,T^G),\; Q'\in\cQ(G,T^G):\;\dim P'<\dim T,\\
\dim P'+\dim Q'\le\rk\,G,\; P'\cap Q'\cong K}}
\!\!\!\!\!\!\!\!\!\!\!\!\!\!\!\!\!\!\!\!\! \Up\bigl([{\ts
X_{T^G}^{P'}}\,]\bigr)\Up([Q'])^{-1}E(G,T^G,Q').
\label{mi6eq5}
\end{gather}

Since $X=G/T$ we find that $X_{T^G}^T=N_G(T)/T$, so that
\begin{equation*}
\Up([{\ts X_{T^G}^T}\,])\!=\!\Up([N_G(T)])/\Up([T])\!=\!
(\ell\!-\!1)^{\dim Q-\rk\,G}\md{N_G(T)/C_G(T)}\Up([C_G(T)]).
\end{equation*}
As $G$ is very special $C_G(T)$ is special, and
$\rk\,C_G(T)=\rk\,G$, so
$\Up([C_G(T)])^{-1}\in(\ell-1)^{-\rk\,G}\La^\ci$ by Lemma
\ref{mi6lem1}. The orbits of $T^G$ on $X_{T^G}^{P'}$ are all
isomorphic to $T^G/P'\cong\bG_m^{\rk\,G-\dim P'}$, so the argument
of Lemma \ref{mi6lem1} shows that
$\Up([X_{T^G}^{P'}])\in(\ell-1)^{\rk\,G-\dim P'}\La^\ci$. Combining
these with \eq{mi6eq5} shows that
\e
\begin{gathered}
E(G,T^G,Q)\!=\!\!\!\!\!\!\!\!\!
\sum_{\begin{subarray}{l} P'\in\cP(X,T^G),\; Q'\in\cQ(G,T^G):
\dim P'<\rk G-\dim Q,\\ \dim P'+\dim Q'\le\rk\,G,\; P'\cap Q'\cong
K\end{subarray}
\!\!\!\!\!\!\!\!\!\!\!\!\!\!\!\!\!\!\!\!\!\!\!\!\!\!\!\!\!\!
\!\!\!\!\!\!\!\!\!\!\!\!\!\!\! }
\!\!\!\!\!\!\!\!\!\!\!\!\!\!\!\!\!\!\!\!\!\!\!
\begin{aligned}[t]
(\text{term in $(\ell-1)^{\rk\,G-\dim P'-\dim Q'}\La^\ci$})\,\cdot&\\
E(G,T^G,Q')&.
\end{aligned}
\end{gathered}
\label{mi6eq6}
\e

Let $k=1,\ldots,\rk\,G$ be given, and suppose by induction that
\begin{itemize}
\item[$(*_k)$] If $Q\in\cQ(G,T^G)$ then $E(G,T^G,Q)$ lies in $\La^\ci$
when $\dim Q\ge k$, and in $(\ell-1)^{\dim Q-k}\La^\ci$ otherwise.
\end{itemize}
When $k=\rk\,G$ this is immediate, from above. Supposing $(*_k)$
holds we can use \eq{mi6eq6} to prove $(*_{k-1})$, by applying
$(*_k)$ to $E(G,T^G,Q')$ in \eq{mi6eq6} and thinking carefully about
the effect of the inequalities $\dim P'<\rk G-\dim Q$, $\dim P'+\dim
Q'\le\rk\,G$. So $(*_0)$ holds by induction. This completes the
proof, except for the problem in choosing $T$ above when $\dim Q=0$.
We can solve this by a similar argument involving $[X/G]$ for
$X=G/Q$ and $Q\in\cQ(G,T^G)$ finite.
\end{proof}

So we may define:

\begin{dfn} Let Assumptions \ref{mi4ass} and \ref{mi6ass} hold and
$G$ be a very special algebraic $\K$-group with maximal torus $T^G$.
For all $Q\in\cQ(G,T^G)$ set $F(G,T^G,Q)=\pi(E(G,T^G,Q))$. This is
well-defined by Theorem~\ref{mi6thm1}.
\label{mi6def1}
\end{dfn}

We now continue Example \ref{mi5ex}, and calculate the
$E,F(G,T^G,Q)$ when $G=\GL(m,\K)$. Let $G,T^G,m,n,\phi$ be as in
Example \ref{mi5ex}, and define $Q$ in $\cQ(G,T^G)$ by \eq{mi5eq2}.
Write $m_k=\md{\phi^{-1}(\{k\})}$ for $k=1,\ldots,n$, so that
$m=m_1+\cdots+m_n$. From \eq{mi5eq22} with $Q'=Q$ we see that
\begin{equation*}
E(G,T^G,Q)=\frac{\md{\{w\in W(G,T^G):w\vert_Q=
\id_Q\}}}{\md{W(G,T^G)}}E(C_G(Q),T^G,Q).
\end{equation*}
Now $C_G(Q)\cong\prod_{k=1}^n\GL\bigl(m_k,\K\bigr)$ by \eq{mi5eq1}
with $Q\cong\prod_{k=1}^n\bG_m\cdot\id_{m_k}$, where $\id_{m_k}$ is
the identity matrix in $\GL(m_k,\K)$. As
\begin{equation*}
E(G\t H,T^G\t T^H,Q_G\t Q_H)=E(G,T^G,Q_G)\cdot E(H,T^H,Q_H),
\end{equation*}
we deduce that
\e
E\bigl(\GL(m,\K),\bG_m^m,Q\bigr)=\ts\frac{1}{m!}
\prod_{k=1}^nm_k!E(m_k),
\label{mi6eq7}
\e
where $E(m)=E\bigl(\GL(m,\K),\bG_m^m,\bG_m\cdot\id_m\bigr)$.
Applying $\pi$ gives
\e
F\bigl(\GL(m,\K),\bG_m^m,Q\bigr)=\ts\frac{1}{m!}
\prod_{k=1}^nm_k!F(m_k),
\label{mi6eq8}
\e
where~$F(m)=F\bigl(\GL(m,\K),\bG_m^m,\bG_m\cdot\id_m\bigr)$.

So it is enough to compute $E(m),F(m)$. For small values of $m$ we
can do this directly using \eq{mi5eq17}, Example \ref{mi5ex} and
Lemma \ref{mi4lem2}, giving
\begin{gather*}
E(1)=1,\quad F(1)=1,\quad
E(2)=(\ell+1)^{-1}\bigl(-\ell^{-1}-\ha\bigr), \quad
F(2)=-\ts\frac{3}{4},
\\
E(3)=(\ell^2+\ell+1)^{-1}\bigl(\ell^{-3}+\ell^{-2}+\ell^{-1}
+\ts\frac{1}{3}\bigr), \quad F(3)=\ts\frac{10}{9},\;\ldots.
\end{gather*}
For larger values of $m$ it is helpful to have an inductive formula
for $E(m),F(m)$. By writing $[[\mathbb{KP}^m/\GL(m+1,\K)]]=
[[\Spec\K/\GL(m,\K)\ltimes\bA^m]]$ in two different ways in
$\uoSF(\Spec\K,\Up,\La)$, after some calculation we find that
\e
\begin{split}
&\sum_{n=1}^{m+1}\frac{1}{n!}\sum_{m+1=m_1+\cdots+m_n,\; m_k\ge 1}
\prod_{k=1}^n\frac{\ell^{m_k}-1}{\ell-1}E(m_k)=\\
&\ell^{-m}\sum_{n=1}^m\frac{(-1)^n}{n!}\sum_{m=m_1+\cdots+m_n,\;
m_k\ge 1} \prod_{k=1}^n\frac{\ell^{m_k}-1}{\ell-1}E(m_k).
\end{split}
\label{mi6eq9}
\e
Applying $\pi$ then gives
\e
\begin{split}
&\sum_{n=1}^{m+1}\frac{1}{n!}\sum_{m+1=m_1+\cdots+m_n,\; m_k\ge 1}
\prod_{k=1}^nm_kF(m_k)=\\
&\sum_{n=1}^m\frac{(-1)^n}{n!}\sum_{m=m_1+\cdots+m_n,\; m_k\ge 1}
\prod_{k=1}^nm_kF(m_k).
\end{split}
\label{mi6eq10}
\e
Here \eq{mi6eq9} and \eq{mi6eq10} contain $(\ell^m+\cdots+1)E(m+1)$
and $(m+1)F(m+1)$ on the top lines with $n=1$, and all other terms
involve $E(m')$ and $F(m')$ for $m'\le m$. So we can use
\eq{mi6eq9}--\eq{mi6eq10} to find $E(m),F(m)$ inductively.

\subsection{Spaces $\uoSF,\oSF(\F,\Up,\La^\ci)$,
$\uoSF,\oSF(\F,\Th,\Om)$, $\uhSF,\hSF(\F,\chi,\Q)$}
\label{mi63}

We restrict the spaces $\oSF(*,\Up,\La)$ in \S\ref{mi53} to
$\La^\ci$, and then project to~$\Om$.

\begin{dfn} Let Assumptions \ref{mi4ass} and \ref{mi6ass} hold,
and $\F$ be an algebraic $\K$-stack with affine geometric
stabilizers. Consider pairs $(\R,\rho)$, where $\R$ is a finite type
algebraic $\K$-stack with affine geometric stabilizers and
$\rho:\R\ra\F$ is a 1-morphism, with {\it equivalence} of pairs as
in Definition \ref{mi3def1}. Define $\uoSF,\oSF(\F,\Up,\La^\ci)$ to
be the $\La^\ci$-modules generated by equivalence classes
$[(\R,\rho)]$ as above, with $\rho$ representable for
$\oSF(\F,\Up,\La^\ci)$, and with relations Definition
\ref{mi5def9}(i)--(iii). These make sense with $\La^\ci$ in place of
$\La$ since (i)--(iii) only involve multiplying by elements of
$\La^\ci$ in $\La$. In particular, $\Up([U])\in\La^\ci$ in (ii), and
$E(G,T^G,Q)\in\La^\ci$ in equation \eq{mi5eq18} of (iii) by
Theorem~\ref{mi6thm1}.

Define $\uoSF,\oSF(\F,\Th,\Om)$ to be the $\Om$-modules generated by
equivalence classes $[(\R,\rho)]$ as above, with $\rho$
representable for $\oSF(\F,\Th,\Om)$, and with relations Definition
\ref{mi5def9}(i)--(iii) projected to $\Om$ using $\pi$ in the
obvious way. That is, in (ii) we have $[(\R\t U,\rho\ci\pi_\R)]
=\Th([U])[(\R,\rho)]$, and \eq{mi5eq18} becomes
\begin{equation*}
[(\R,\rho)]=\ts\sum_{Q\in\cQ(G,T^G)}F(G,T^G,Q)
\bigl[\bigl([X/Q],\rho\ci\io^Q\bigr)\bigr].
\end{equation*}
Since $\pi:\La^\ci\ra\Om$ is supposed {\it surjective} we have
isomorphisms
\e
\begin{split}
\uoSF(\F,\Th,\Om)&\cong\uoSF(\F,\Up,\La^\ci)/
(\Ker\pi\cdot\uoSF(\F,\Up,\La^\ci)),\\
\oSF(\F,\Th,\Om)&\cong\oSF(\F,\Up,\La^\ci)/
(\Ker\pi\cdot\oSF(\F,\Up,\La^\ci)),
\end{split}
\label{mi6eq11}
\e
where $\Ker\pi$ is an ideal in $\La^\ci$. Define projections
\begin{align*}
&\bar\Pi^{\Up,\La^\ci}_\F:\uSF(\F)\longra\uoSF(\F,\Up,\La^\ci), &
&\bar\Pi^{\Up,\La}_\F:\uoSF(\F,\Up,\La^\ci)\longra\uoSF(\F,\Up,\La),\\
&\bar\Pi^{\Th,\Om}_\F:\uSF(\F)\longra\uoSF(\F,\Th,\Om), &
&\bar\Pi^{\Th,\Om}_\F:\uoSF(\F,\Up,\La^\ci)\longra\uoSF(\F,\Th,\Om)
\end{align*}
by \eq{mi4eq3}, replacing $c_i$ by $\pi(c_i)$ on the r.h.s.\ for
$\bar\Pi^{\Th,\Om}_\F:\uSF(\F,\Up,\La^\ci)\ra\uoSF(\F,\Th,\Om)$.
These are well-defined since they map relations in the domain to
relations in the target, as in Definitions \ref{mi4def1} and
\ref{mi5def9}. Define multiplication `$\,\cdot\,$', pushforwards
$\phi_*$, pullbacks $\phi^*$, tensor products $\ot$ and operators
$\Pi^\mu,\Pi^\vi_n,\hat\Pi^\nu_*$ on the spaces
$\uoSF,\oSF(*,\Up,\La^\ci)$ and $\uoSF,\oSF(*,\Th,\Om)$ exactly as
in Definition~\ref{mi5def10}.
\label{mi6def2}
\end{dfn}

From the proofs of Theorems \ref{mi4thm2} and \ref{mi5thm3} we
deduce the analogous result for the $\uoSF,\oSF(*,\Up,\La^\ci)$.
This is nearly immediate, as the relations in $\uoSF,\oSF(*,\Up,
\La^\ci)$ are the same as in $\uoSF,\oSF(*,\Up,\La)$. We know that
under the operations `$\,\cdot\,$'$,\ldots,\hat\Pi^\nu_*$ relations
are taken to linear combinations of relations with coefficients in
$\La$, and we must check these coefficients may be chosen in
$\La^\ci$, which is fortunately obvious. Projecting coefficients
from $\La^\ci$ to $\Om$ using $\pi$ proves the same thing for the
$\uoSF,\oSF(*,\Th,\Om)$, giving:

\begin{thm} These operations `$\,\cdot\,$'$,\phi_*,\phi^*,\ot,\Pi^\mu,
\Pi^\vi_n,\hat\Pi^\nu_*$ on $\uoSF,\oSF(*,\Up,\La^\ci)$ and\/
$\uoSF,\oSF(*,\Th,\Om)$ are compatible with the relations, and so
are well-defined.
\label{mi6thm2}
\end{thm}

As for Corollaries \ref{mi4cor} and \ref{mi5cor}, we deduce:

\begin{cor} The projections $\bar\Pi^{\Up,\La^\ci}_*,
\bar\Pi^{\Up,\La}_*,\bar\Pi^{\Th,\Om}_*$ commute with the operations
`$\,\cdot\,$'$,\phi_*,\phi^*,\ot,\Pi^\mu,\Pi^\vi_n,\hat\Pi^\nu_*$ on
$\uSF(*),\uoSF(*,\Up,\La^\ci),\uoSF(*,\Up,\La),\uoSF(*,\Th,\Om)$.
The analogues of Theorems \ref{mi3thm1} and \ref{mi5thm2} and
Proposition \ref{mi5prop2} hold for the spaces
$\uoSF,\oSF(*,\Up,\La^\ci)$ and\/~$\uoSF,\oSF(*,\Th,\Om)$.
\label{mi6cor}
\end{cor}

The analogues of Propositions \ref{mi5prop3} and \ref{mi5prop4}
apply for the $\uoSF,\oSF(*,\Up,\La^\ci)$ and
$\uoSF,\oSF(*,\Th,\Om)$, replacing $\La$ by $\La^\ci,\Om$. The
analogue of Proposition \ref{mi5prop5} is:

\begin{prop} Define commutative $\La^\ci$- and\/ $\Om$-algebras
$\bar\La^\ci,\bar\Om$ with\/ $\La^\ci$- and\/ $\Om$-bases
isomorphism classes $[T]$ of\/ $\K$-groups $T$ of the form
$\bG_m^k\t K$, for $k\ge 0$ and\/ $K$ finite abelian, with
multiplication $[T][T']=[T\t T']$. Define
$\bar\imath_{\La^\ci}:\bar\La^\ci\ra\uoSF(\Spec\K,\Up,\La^\ci)$
and\/ $\bar\imath_\Om:\bar\Om\ra\uoSF(\Spec\K,\Th,\Om)$ by
$\sum_ic_i[T_i]\mapsto\sum_ic_i[[\Spec\K/T_i]]$. Then
$\bar\imath_{\La^\ci},\bar\imath_\Om$ are algebra isomorphisms. They
restrict to isomorphisms from $\La^\ci[\{1\}],\Om[\{1\}]$ to
$\oSF(\Spec\K,\Up,\La^\ci),\oSF(\Spec\K,\Th,\Om)$.
\label{mi6prop1}
\end{prop}

\begin{proof} For the $\La^\ci$ case we follow Proposition
\ref{mi5prop5} replacing $\La$ by $\La^\ci$ throughout, {\it except
that\/} to deduce injectivity we apply $i_\La^{-1}\ci
\Pi^{\Up,\La}_{\Spec\K}\ci\bar\Pi^{\Up,\La}_{\Spec\K}$ to project to
$\La$, not $\La^\ci$. The $\Om$ case then follows using
\eq{mi6eq11}, since~$\bar\Om\cong\bar\La^\ci/
(\Ker\pi\cdot\bar\La^\ci)$.
\end{proof}

If $\F$ is a finite type algebraic $\K$-stack with affine geometric
stabilizers then $[\F]\in\uoSF(\Spec\K,\Up,\La^\ci)$ or
$\uoSF(\Spec\K,\Th,\Om)$ and so $\bar\imath_{\La^\ci}^{-1}([\F]),
\bar\imath_\Om^{-1}([\F])$ lie in $\bar\La^\ci,\bar\Om$. We can
regard these as generalizations of $\Up'$ in Theorem \ref{mi4thm1},
which work even when $\ell=1$. In particular, when $\Om=\Q$ and
$\Th=\chi$ as in Example \ref{mi6ex2}(a), $\bar\imath_\Q^{-1}([\F])
\in\bar\Q$ is a kind of {\it generalized Euler characteristic}
of~$\F$.

As for Proposition \ref{mi5prop6} we have:

\begin{prop} The following maps are\/ $\La^\ci$- or\/ $\Om$-linear
and injective:
\begin{align*}
(\bar\Pi^{\Up,\La^\ci}_\F\ci\io_\F)\ot_\Q\id_{\La^\ci}&:
\CF(\F)\ot_\Q\La^\ci\longra\oSF(\F,\Up,\La^\ci),\\
(\bar\Pi^{\Up,\Om}_\F\ci\io_\F)\ot_\Q\id_\Om&:
\CF(\F)\ot_\Q\Om\,\,\longra\oSF(\F,\Up,\Om),\\
\mu\ci\bigl((\bar\Pi^{\Up,\La^\ci}_\F\ci\io_\F)\ot_\Q
\bar\imath_{\La^\ci}\bigr)&:
\CF(\F)\ot_\Q\bar\La^\ci\longra\uoSF(\F,\Up,\La^\ci),\\
\mu\ci\bigl((\bar\Pi^{\Up,\Om}_\F\ci\io_\F)\ot_\Q
\bar\imath_\Om\bigr)&:
\CF(\F)\ot_\Q\bar\Om\,\,\longra\uoSF(\F,\Up,\Om),
\end{align*}
defining $\mu$ as in Proposition~\ref{mi5prop6}.
\label{mi6prop2}
\end{prop}

We generalize the $\pi^\stk_\F$ of \eq{mi3eq3} to
$\oSF(\F,\Up,\La^\ci)$ and~$\oSF(\F,\Th,\Om)$.

\begin{dfn} Let Assumptions \ref{mi4ass} and \ref{mi6ass} hold
with $\K$ of characteristic zero. Suppose ${\rm X}:\La^\ci\ra\Q$ or
${\rm X}:\Om\ra\Q$ is an algebra morphism with ${\rm X}\ci\Up([U])=
\chi([U])$ or ${\rm X}\ci\Th([U])=\chi([U])$ for all quasiprojective
$\K$-varieties $U$, where $\chi$ is the Euler characteristic. Such
morphisms $\rm X$ exist in all of Examples
\ref{mi6ex1}--\ref{mi6ex3}. Let $\F$ be an algebraic $\K$-stack with
affine geometric stabilizers. Define
$\bar\pi^\stk_\F:\oSF(\F,\Up,\La^\ci)\ra\CF(\F)$ or
$\bar\pi^\stk_\F:\oSF(\F,\Th,\Om)\ra\CF(\F)$ by
\e
\bar\pi_\F^\stk\bigl(\ts\sum_{i=1}^nc_i[(\R_i,\rho_i)]\bigr)=
\ts\sum_{i=1}^n{\rm X}(c_i)\CF^\stk(\rho_i)1_{\R_i},
\label{mi6eq12}
\e
following \eq{mi3eq3}. By a complicated proof similar to Theorems
\ref{mi4thm2} and \ref{mi5thm3} we can show that $\bar\pi^\stk_\F$
is compatible with the relations defining $\oSF(\F,\Up,\La^\ci)$ and
$\oSF(\F,\Th,\Om)$, and so is well-defined. The analogues of
Propositions \ref{mi3prop1} and \ref{mi3prop2} and Theorem
\ref{mi3thm2} then hold, by the same proofs as in~\S\ref{mi3}.
\label{mi6def3}
\end{dfn}

In the situation of Examples \ref{mi4ex2} and \ref{mi6ex2}(a) we
have $\Om=\Q$ and $\Th=\chi$, the Euler characteristic, so we have
defined spaces $\uoSF,\oSF(\F,\chi,\Q)$ which are very like the
constructible functions $\CF(\F)$ of \S\ref{mi23}, in that
pushforwards $\phi_*$ `integrate' along the fibres of $\phi$ using
$\chi$. Now for $\K$ of characteristic zero, if $\phi:X\ra Y$ is a
fibration of quasiprojective $\K$-varieties with fibre $F$ then
$\chi(X)=\chi(F)\chi(Y)$, even if $\phi$ is {\it not\/} a Zariski
locally trivial fibration. This is a special property of the Euler
characteristic which does not hold for other motivic invariants such
as virtual Poincar\'e polynomials, and lies behind the proof of
\eq{mi2eq3}. We modify the relations in $\uoSF,\oSF(\F,\chi,\Q)$ to
include this.

\begin{dfn} Let $\K$ have characteristic zero, and $\Up,\La,\ldots$
be as in Examples \ref{mi4ex2} and \ref{mi6ex2}(a), with $\Om=\Q$
and $\Th=\chi$. Let $\F$ be an algebraic $\K$-stack with affine
geometric stabilizers. Define spaces $\uhSF(\F,\chi,\Q)$ and
$\hSF(\F,\chi,\Q)$ exactly as for $\uoSF,\oSF(\F,\chi,\Q)$ in
Definition \ref{mi6def2}, but replacing relation (ii) by
\begin{itemize}
\setlength{\itemsep}{0pt}
\setlength{\parsep}{0pt}
\item[(ii$'$)] Let $\mathfrak{Q},\R$ be finite type algebraic
$\K$-stacks with affine geometric stabilizers, $\rho:\R\ra\F$ a
1-morphism, $n\in\Z$ and $\pi_\R:\mathfrak{Q}\ra\R$ a representable
1-morphism such that $U_x=\mathfrak{Q}\t_{\pi_\R,\R,x}\Spec\K$ is a
quasiprojective $\K$-variety with $\chi([U_x])=n$ for all
$x:\Spec\K\ra\R$. Then~$[(\mathfrak{Q},\rho\ci\pi_\R)]
=n[(\R,\rho)]$.
\end{itemize}
Taking $\mathfrak{Q}=\R\t U$ recovers (ii), so this strengthens the
relations. Thus there are natural surjective
projections~$\hat\Pi^{\chi,\Q}_\F:\uoSF,\oSF(\F,\chi,\Q)\ra
\uhSF,\hSF(\F,\chi,\Q)$.
\label{mi6def4}
\end{dfn}

One can then prove that all the material above on operations
`$\,\cdot\,$'$,\phi_*,\ab\phi^*,\ab\ot,\ab\Pi^\mu,\ab\Pi^\vi_n,\ab
\hat\Pi^\nu_\F,\ab\pi^\stk_\F$ and properties of
$\uoSF,\oSF(*,\chi,\Q)$ also works for $\uhSF,\hSF(*,\chi,\Q)$.
Suppose $\F$ is a $\K$-scheme or algebraic $\K$-space, so that its
stabilizer groups are trivial. Proposition \ref{mi5prop3} implies
that $\hSF(\F,\chi,\Q)$ is spanned over $\Q$ by elements
$[(U,\rho)]$ for $U$ a quasiprojective $\K$-variety. Using (ii$'$)
it is then easy to show $\hat\pi^\stk_\F:\hSF(\F,\chi,\Q)\ra\CF(\F)$
is an {\it isomorphism}. Therefore the $\hSF(\F,\chi,\Q)$ coincide
with $\CF(\F)$ for schemes and algebraic spaces.

\medskip

\noindent{\small\sc The Mathematical Institute, 24-29 St. Giles,
Oxford, OX1 3LB, U.K.}

\noindent{\small\sc E-mail: \tt joyce@maths.ox.ac.uk}


\begin{thebibliography}{99}

\bibitem{Bitt} F. Bittner, {\it The universal Euler characteristic
for varieties of characteristic zero}, Comp. Math. 140 (2004),
1011--1032. math.AG/0111062.

\bibitem{Bore} A. Borel, {\it Linear Algebraic Groups}, second
edition, Graduate Texts in Math. 126, Springer-Verlag, New York,
1991.

\bibitem{Chea} J. Cheah, {\it On the cohomology of Hilbert
schemes of points}, J. Algebraic Geometry 5 (1996), 479--511.

\bibitem{Chev} {\it Anneaux de Chow et applications}, S\'eminaire
C. Chevalley, 2e ann\'ee, \'Ecole Normale Sup\'erieure,
Secr\'etariat math\'ematique, Paris, 1958.

\bibitem{DaKh} V. Danilov and A. Khovanskii, {\it Newton polyhedra
and an algorithm for computing Hodge--Deligne numbers}, Math.\ USSR
Izvestiya 29 (1987), 279--298.

\bibitem{Deli} P. Deligne, {\it Th\'eorie de Hodge I}, pages 425--430
in Actes du Congr\`es Int. Math. (Nice, 1970), vol. 1,
Gauthier--Villars, Paris, 1971.

\bibitem{Gome} T.L. G\'omez, {\it Algebraic stacks}, Proc. Indian Acad.
Sci. Math. Sci. 111 (2001), 1--31. math.AG/9911199.

\bibitem{Joyc1} D.D. Joyce, {\it Constructible functions on
Artin stacks}, J. London Math. Soc. 74 (2006), 583--606.
math.AG/0403305.

\bibitem{Joyc2} D.D. Joyce, {\it Configurations in abelian
categories. I. Basic properties and moduli stacks}, Advances in
Math. 203 (2006), 194--255. math.AG/0312190.

\bibitem{Joyc3} D.D. Joyce, {\it Configurations in abelian
categories. II. Ringel--Hall algebras}, Advances in Math. 210
(2007), 635--706. math.AG/0503029.

\bibitem{Joyc4} D.D. Joyce, {\it Configurations in abelian
categories. III. Stability conditions and identities},
math.AG/0410267, 2004. To appear in Advances in Mathematics.

\bibitem{Joyc5} D.D. Joyce, {\it Configurations in abelian
categories. IV. Invariants and changing stability conditions},
math.AG/0410268, version 4, 2006.

\bibitem{Kres} A. Kresch, {\it Cycle groups for Artin stacks},
Invent. math. 138 (1999), 495--536. math.AG/9810166.

\bibitem{LaMo} G. Laumon and L. Moret-Bailly, {\it Champs alg\'ebriques},
Ergeb. der Math. und ihrer Grenzgebiete 39, Springer-Verlag, Berlin, 2000.

\bibitem{Looi} E. Looijenga, {\it Motivic measures}, Asterisque 276
(2002), 267--297. \hfil\break
math.AG/0006220.

\bibitem{Toen} B. Toen, {\it Anneaux de Grothendieck des $n$-champs
d'Artin}, \hfil\break math.AG/0509098, 2005.

\bibitem{Yasu} T. Yasuda, {\it Motivic integration over
Deligne--Mumford stacks}, Advances in Math. 207 (2006), 707--761.
math.AG/0312115.

\end{thebibliography}
\end{document}